\documentclass[10pt]{article}

\usepackage{mathpazo}   
\usepackage{comment}
\usepackage[title]{appendix}
\usepackage{wrapfig}
\usepackage{subcaption}
\usepackage{tikz}
\usepackage{pgfplots}
\pgfplotsset{compat=1.18}
\usetikzlibrary{patterns, pgfplots.fillbetween}
\graphicspath{{.}}
\usepackage[font=footnotesize]{caption}
\usepackage{nameref}
\usepackage{mathtools}
\usepackage{empheq}
\usepackage{comment}
\usepackage[shortlabels,inline]{enumitem}
\setlist[enumerate]{nosep}
\usepackage[colorlinks=true,
linkcolor=refkey,
urlcolor=lblue,
citecolor=red]{hyperref}
\usepackage[doc,wmm,hhb]{optional}
\usepackage[dvipsnames,svgnames,x11names]{xcolor}

\usepackage{float}
\usepackage{soul}
\usepackage{graphicx}
\definecolor{labelkey}{rgb}{0,0.08,0.45}
\definecolor{refkey}{rgb}{0,0.6,0.0}
\definecolor{Brown}{rgb}{0.45,0.0,0.05}
\definecolor{lime}{rgb}{0.00,0.8,0.0}
\definecolor{lblue}{rgb}{0.5,0.5,0.99}
\definecolor{OliveGreen}{rgb}{0,0.6,0}
\definecolor{tyrianpurple}{rgb}{0.4, 0.01, 0.24}

\colorlet{hlcyan}{cyan!30}

\usepackage{stmaryrd}
\usepackage{amssymb}

\makeatletter
\def\namedlabel#1#2{\begingroup
   \def\@currentlabel{#2}%
   \label{#1}\endgroup
}
\makeatother

\newcommand{\seppthree}{\setlength{\itemsep}{-3pt}}

\usepackage[margin=0.92in,footskip=0.25in]{geometry}
\newcommand{\X}{L}
\newcommand{\Y}{U}
\newcommand{\com}{\ensuremath{\operatorname{CE}}_{\ominus}}

\providecommand{\siff}{\Leftrightarrow}

\newcommand{\nnn}{\ensuremath{{n\in{\mathbb N}}}}

\newcommand{\menge}[2]{\big\{{#1}~\big |~{#2}\big\}}

\newcommand{\fenv}[1]%
{\ensuremath{\,\overrightarrow{\operatorname{env}}_{#1}}}
\newcommand{\benv}[1]%
{\ensuremath{\,\overleftarrow{\operatorname{env}}_{#1}}}

\newcommand{\scal}[2]{\left\langle{#1},{#2}  \right\rangle}

\newcommand{\RR}{\ensuremath{\mathbb R}}
\newcommand{\HH}{\ensuremath{\mathcal H}}

\newcommand{\RX}{\ensuremath{\,\left]-\infty,+\infty\right]}}
\newcommand{\RRX}{\ensuremath{\,\left[-\infty,+\infty\right]}}

\newcommand{\NN}{\ensuremath{\mathbb N}}

\newcommand{\dom}{\ensuremath{\operatorname{dom}}\,}

\newcommand{\inte}{\ensuremath{\operatorname{int}}}

\newcommand{\bd}{\ensuremath{\operatorname{bdry}}}

\newcommand{\minf}{\ensuremath{-\infty}}
\newcommand{\pinf}{\ensuremath{+\infty}}

{\begin{list}{}{%
\settowidth{\labelwidth}{\textrm{#1~}}%
\setlength{\leftmargin}{\labelwidth+\labelsep}}}
{\end{list}}
\usepackage{amsthm}
\usepackage{aliascnt}
\makeatletter
\def\th@plain{%
	\thm@notefont{}
	\itshape 
}
\def\th@definition{%
	\thm@notefont{}
	\normalfont 
}
\makeatother
\usepackage[capitalize,nameinlink]{cleveref}
\crefname{equation}{}{equations}

\crefname{item}{}{items}
\crefname{enumi}{}{}
\newtheorem{theorem}{Theorem}[section]

\newaliascnt{lemma}{theorem}
\newtheorem{lemma}[lemma]{Lemma}
\aliascntresetthe{lemma}
\crefname{lemma}{Lemma}{Lemmas}
\Crefname{lemma}{Lemma}{Lemmas}

\newaliascnt{lem}{theorem}

\aliascntresetthe{lem}
\crefname{lem}{Lemma}{Lemmas}
\Crefname{lem}{Lemma}{Lemmas}

\newaliascnt{corollary}{theorem}
\newtheorem{corollary}[corollary]{Corollary}
\aliascntresetthe{corollary}
\crefname{corollary}{Corollary}{Corollaries}
\Crefname{corollary}{Corollary}{Corollaries}

\newaliascnt{cor}{theorem}

\aliascntresetthe{cor}
\crefname{cor}{Corollary}{Corollaries}
\Crefname{cor}{Corollary}{Corollaries}

\newaliascnt{proposition}{theorem}
\newtheorem{proposition}[proposition]{Proposition}
\aliascntresetthe{proposition}
\crefname{proposition}{Proposition}{Propositions}
\Crefname{proposition}{Proposition}{Propositions}

\newaliascnt{prop}{theorem}

\aliascntresetthe{prop}
\crefname{prop}{Proposition}{Propositions}
\Crefname{prop}{Proposition}{Propositions}

\newaliascnt{definition}{theorem}
\newtheorem{definition}[definition]{Definition}
\aliascntresetthe{definition}
\crefname{definition}{Definition}{Definitions}
\Crefname{definition}{Definition}{Definitions}

\newaliascnt{defn}{theorem}
\newtheorem{defn}[defn]{Definition}
\aliascntresetthe{defn}
\crefname{defn}{Definition}{Definitions}
\Crefname{defn}{Definition}{Definitions}

\newaliascnt{thm}{theorem}

\aliascntresetthe{thm}
\crefname{thm}{Theorem}{Theorems}
\Crefname{thm}{Theorem}{Theorems}

\newaliascnt{example}{theorem}
\newtheorem{example}[example]{Example}
\aliascntresetthe{example}
\crefname{example}{Example}{Examples}
\Crefname{example}{Example}{Examples}

\newaliascnt{ex}{theorem}
\newtheorem{ex}[ex]{Example}
\aliascntresetthe{ex}
\crefname{ex}{Example}{Examples}
\Crefname{ex}{Example}{Examples}

\newaliascnt{fact}{theorem}
\newtheorem{fact}[fact]{Fact}
\aliascntresetthe{fact}
\crefname{fact}{Fact}{Facts}
\Crefname{fact}{Fact}{Facts}

\newaliascnt{remark}{theorem}
\newtheorem{remark}[remark]{Remark}
\aliascntresetthe{remark}
\crefname{remark}{Remark}{Remarks}
\Crefname{remark}{Remark}{Remarks}

\newaliascnt{rem}{theorem}

\aliascntresetthe{rem}
\crefname{rem}{Remark}{Remarks}
\Crefname{rem}{Remark}{Remarks}

\crefname{chapter}{Appendix}{chapters}

\providecommand{\LA}{\Leftarrow}
\providecommand{\RA}{\Rightarrow}

\providecommand{\RR}{\mathbb{R}}

\providecommand{\dom}{\operatorname{dom}}

\providecommand{\gra}{\operatorname{gra}}

\providecommand{\NN}{\mathbb{N}}

\providecommand{\U}{ {U}}

\providecommand{\RR}{\mathbb{R}}
\providecommand{\NN}{\mathbb{N}}

\definecolor{myblue}{rgb}{0.9,0.9,0.98}

  \newcommand*\mybluebox[1]{%
    \colorbox{myblue}{\hspace{1em}#1\hspace{1em}}}

\allowdisplaybreaks 

\usepackage{relsize}

\begin{document}

\setlength{\abovedisplayskip}{8pt}
\setlength{\belowdisplayskip}{8pt}	
\newsavebox\myboxA
\newsavebox\myboxB
\newlength\mylenA

\newcommand*\xoverline[2][0.75]{%
    \sbox{\myboxA}{$#2$}%
    \setbox\myboxB\null
    \ht\myboxB=\ht\myboxA%
    \dp\myboxB=\dp\myboxA%
    \wd\myboxB=#1\wd\myboxA
    \sbox\myboxB{$\overline{\copy\myboxB}$}
    \setlength\mylenA{\the\wd\myboxA}
    \addtolength\mylenA{-\the\wd\myboxB}%
    \ifdim\wd\myboxB<\wd\myboxA%
       \rlap{\hskip 0.5\mylenA\usebox\myboxB}{\usebox\myboxA}%
    \else
        \hskip -0.5\mylenA\rlap{\usebox\myboxA}{\hskip 0.5\mylenA\usebox\myboxB}%
    \fi}
\makeatother

\makeatletter
\renewcommand*\env@matrix[1][\arraystretch]{%
  \edef\arraystretch{#1}%
  \hskip -\arraycolsep
  \let\@ifnextchar\new@ifnextchar
  \array{*\c@MaxMatrixCols c}}
\makeatother

\providecommand{\wbar}{\xoverline[0.9]{w}}
\providecommand{\ubar}{\xoverline{u}}

\newcommand{\nn}[1]{\ensuremath{\textstyle\mathsmaller{({#1})}}}
\newcommand{\crefpart}[2]{%
  \hyperref[#2]{\namecref{#1}~\labelcref*{#1}~\ref*{#2}}%
}
\newcommand\bigzero{\makebox(0,0){\text{\LARGE0}}}
	

%

\author{
Sedi Bartz\thanks{
Mathematics and Statistics, UMass Lowell, MA 01854, USA. E-mail:
\texttt{sedi\_bartz@uml.edu}.},~
Heinz H.\ Bauschke\thanks{
Mathematics, University
of British Columbia,
Kelowna, B.C.\ V1V~1V7, Canada. E-mail:
\texttt{heinz.bauschke@ubc.ca}.}~~~and~
Yuan Gao\thanks{
Mathematics, University
of British Columbia,
Kelowna, B.C.\ V1V~1V7, Canada. 
 E-mail: \texttt{y.gao@ubc.ca}.}
}

\title{\textsf{ 
On the equivalence of $c$-potentiability and $c$-path boundedness\\ in the sense of 
Artstein-Avidan, Sadovsky, and Wyczesany
}
}

\date{August 17, 2026}

\maketitle

\begin{abstract}
A cornerstone of convex analysis, established by Rockafellar in 1966, asserts that a set has a potential if and only if it is cyclically monotone. This characterization was generalized to hold for any {\color{black} finite-valued} cost function $c$ and lies at the core structure of optimal transport plans. However, this equivalence fails to hold for costs that attain infinite values. In this paper, we explore potentiability for an infinite-valued cost $c$ under the assumption of $c$-path boundedness, a condition that was first introduced by Artstein-Avidan, Sadovsky and Wyczesany. 
{\color{black} 
This condition is necessary for potentiability and it is more restrictive than $c$-cyclic monotonicity; 
however, this condition alone does not imply potentiability.
}
We provide general settings and other conditions under which $c$-path boundedness is sufficient for potentiability, and therefore equivalent. We provide a general theorem for potentiability, requiring no topological assumptions on the spaces or the cost. We then provide sufficiency in separable metric spaces and costs that are continuous in their domain. Finally, we introduce the notion of a $c$-path bounded extension and use it to prove the existence of potentials for a special class of costs on $\mathbb{R}^2$. 
We illustrate our discussion and results with several examples. 
\end{abstract}
{ 


\section{Introduction}
\label{s:intro}
Optimal transport theory has rapidly evolved with vast and far-reaching applications in fields of mathematics, physics, economics and more (see, for example, \cite{Friesecke, Peyre, Santambrogio, Villani}). Given a cost $c(x,y)$ of transporting a mass unit from the location $x$ in the space $X$ to a location $y$ in the space $Y$, and given a probability mass distribution $\mu$ in $X$ and a probability distribution $\nu$ in $Y$, the optimal transport problem consists of finding a transport plan $\pi$ (a probability distribution in $\Pi(\mu,\nu)$, the set of probability distributions on $X\times Y$ with marginal distributions $\mu$ and $\nu$) such that the total cost of transportation is minimal, that is, one would like to find a minimizing plan $\pi$ to the optimal transport problem
$$
\inf_{\pi\in\Pi(\mu,\nu)}\int_{X\times Y}c(x,y)d\pi(x,y).
$$
Emerging from Brenier's work on the quadratic cost, and then generalized to arbitrary {\color{black} finite-valued 
 (also known as real-valued)} cost functions $c$ (see, for example, \cite{Villani} for an account), it is well known that, under suitable conditions, a transport plan $\pi$ is optimal if and only if it is concentrated in a subset $G$ of $X\times Y$ for which there is a potential function $\phi:X\to\RX$ such that $G\subseteq\partial_c\phi$, 
where $\partial_c\phi$ is the $c$-subdifferential of $\phi$ (see \cref{d:c-sub}). We shall call such a set $G$ a $c$-potentiable set.

At the heart of this fact and of modern proofs of the Kantorovich Duality in optimal transport is the characterization that a set $G$ is $c$-potentiable if and only if it is $c$-cyclically monotone (\cref{d:max cost}).   
In 1966, Rockafellar \cite{Rocky66} was the first to introduce cyclic monotonicity and proved that this notion is equivalent to the existence of a convex potential. In the mid-1980s and early 1990s, this celebrated result was generalized to hold for an arbitrary {\color{black} finite-valued} cost function $c: X \times Y \to \mathbb{R}$ by Rochet \cite{Rochet} and by R\"uschendorf \cite{Rueschendorf} who showed the equivalence between $c$-cyclic monotonicity (\cref{d:max cost}) and $c$-potentiability (\cref{d:c-potential}). 

However, if $c$ is a non-traditional cost, i.e., $c$ is allowed to take on the value $\pinf$, then $c$-cyclic monotonicity does not guarantee optimality of a transport plan (see \cite[Example~3.1]{Ambrosio}) and it does not imply $c$-potentiability (see \cite[Section~2.8]{Kasia}). 
Non-traditional costs are natural in various problems emerging from physics and other fields and have gained interest in applications and studies of optimal transport (see~\cite{Kasia, Pascale} for recent accounts). 

In 2022, Artstein-Avidan, Sadovsky and Wyczesany \cite{Kasia} introduced the concept of $c$-path boundedness of a set $G$, a necessary condition for $c$-potentiability and a more restrictive notion than $c$-cyclic monotonicity for non-traditional costs. They proved that $c$-path boundedness is sufficient for $c$-potentiability (and therefore equivalent) under the condition that $G$ is countable in~\cite{Kasia} and under the ``no infinite black hole" assumption on $G$ in a later publication \cite{Kasiafix}. 
In this paper, we study the sufficiency of $c$-path boundedness for $c$-potentiability in cases and scenarios where the ``no infinite black hole" condition is not satisfied. Our discussion does not invoke measure-theoretic aspects. However, one of the strongest motivations for our study stems from the fact that the set $G$ is where a transport plan is concentrated and from the fact that when the total cost of transportation is finite, then potentiability does imply optimality of the transport plan (see~\cite{Beigl}).
We establish sufficiency of $c$-path boundedness for $c$-potentiability in settings commonly used in discussions and studies of optimal transport. In particular, we establish sufficiency in settings that include non-traditional costs which are continuous on their domain in separable metric spaces. We elaborate on the main contributions of this paper further with more details at the end of this section. To this end, we will use the following settings and preliminary facts.

Throughout, we assume that 
\begin{empheq}[box=\mybluebox]{equation}
\text{$X$ and $Y$ are two nonempty sets
}
\end{empheq}
and 
\begin{empheq}[box=\mybluebox]{equation}
\text{$c\colon X\times Y \to \RR\cup\{+\infty\}$}
\end{empheq}
is a \emph{(non-traditional) cost} function with nonempty domain: 
\begin{empheq}[box=\mybluebox]{equation}
D := \dom c := \menge{(x,y)\in X\times Y}{c(x,y)<+\infty}\neq\varnothing.
\end{empheq}

\begin{defn}[$c$-conjugates and $c$-class functions]
\label{d:conjugate}
 We define the \emph{$c$-conjugate} of a function $\psi: X \to [-\infty, +\infty]$ by \footnote{{\color{black} We use the convention that $\inf \varnothing =+\infty$.}} 
    \begin{equation}\label{equ:c-transf}
    (\forall y\in Y)\quad \psi^c(y) := \inf\menge{c(x,y)-\psi(x)}{x\in X_y}, 
\;\;\text{where}\;\; 
X_y:=\menge{x\in X}{(x,y) \in D}.
    \end{equation}
 Similarly, we define the \emph{$c$-conjugate} of a function $\varphi: Y \to [-\infty, +\infty]$ by 
    \begin{equation}\label{equ:c-transg}
(\forall x\in X) \quad   \varphi^{\overline{c}}(x):= 
\inf\menge{c(x,y) - \varphi(y)}{y\in Y_x}, 
\;\;\text{where}\;\;
    Y_x := \menge{y\in Y}{(x,y)\in D}. 
    \end{equation}
One says that $\psi:  X\to [-\infty, +\infty]$ is a 
    \emph{$c$-class function} if there is a function $\varphi: Y \to [-\infty, +\infty]$ such that $\psi=\varphi^{\overline{c}}$. 
\end{defn}

\begin{defn}[$c$-subdifferential]\label{d:c-sub}
Let $\psi: X \to [-\infty, +\infty]$ and suppose that 
$x\in X$ satisfies $\psi(x) \in \mathbb{R}$. 
Then the \emph{$c$-subdifferential} of $\psi$ at $x$ is
    \begin{equation}\label{equ:c-subdiff}
        \partial_c \psi(x) :=\menge{ y\in Y }{\psi(x)+\psi^c(y) = c(x,y)<+\infty}. 
    \end{equation}
\end{defn}

\begin{remark}
\label{r:c-sub}
    An alternative characterization of the $c$-subdifferential 
$\partial_c\psi(x)$ from \cref{equ:c-subdiff} is 
\begin{equation}
\partial_c \psi(x) 
:= \menge{y\in Y}{c(x,y)<+\infty \text{\ \ and\ \ } (\forall z\in X)\;\;c(x,y)-c(z,y) \leq \psi(x)-\psi(z)}.
\end{equation}
\end{remark}

\begin{remark}[convex analysis setting]
\label{r:classic}
Suppose that $X$ is a real Banach space, $Y=X^*$, and 
$c(x,y)=-\scal{x}{y}$ is the negative pairing of $X\times X^*$. 
Let $f \colon X\to\left]-\infty,+\infty\right]$ be proper, i.e., $f\not\equiv +\infty$. 
Then  $(-f)^c=-f^*$ 
is the negative Fenchel conjugate of $f$, 
and $\partial_c(-f)(x)=\partial f(x)$ is the subdifferential 
of $f$ at $x\in\dom f$ in the sense of convex analysis.  
\end{remark}

Note that the graph of 
$\partial_c\psi$ satisfies 
\begin{equation}
\gra \partial_c\psi := \menge{(x,y)\in X\times Y}{y\in\partial_c\psi(x)}
\subseteq D.
\end{equation}
More generally, we are interested in 
the following properties of subsets of $D$:

\begin{defn}[path, condensation graph, and strongly connected components]
\label{d:path}
Let $G$ be a nonempty subset of $D$, and 
let $(x_s,y_s)$ and $(x_e, y_e)$ be two points in $G$. 
A \emph{path in $G$, with starting point $(x_s,y_s)$ and endpoint $(x_e,y_e)$,} 
is a finite sequence $((x_k, y_k))_{k=1}^{N+1} \in G^{N+1}$ such that 
the following hold: 
\begin{enumerate}
\item $(x_1, y_1) = (x_s,y_s)$,
\item $(x_{N+1},y_{N+1}) = (x_e, y_e)$,
\item $c(x_{k+1}, y_k) < +\infty$ for every $k\in\{1,\ldots, N\}$.
\end{enumerate}
We denote the set of all paths in $G$, with 
starting point $(x_s,y_s)$ and endpoint $(x_e,y_e)$, 
by 
$P^G_{(x_s,y_s)\to (x_e,y_e)}$.
If there is a path from 
$(x_s,y_s)$ to $(x_e,y_e)$ as well as a path
from $(x_e,y_e)$ back to $(x_s,y_s)$, then 
we say that 
$(x_s,y_s)$ and $(x_e,y_e)$ are \emph{mutually reachable} or 
\emph{strongly connected}.
This induces an equivalence relation on $G$, 
and the quotient set is the \emph{condensation graph} of $G$
and its elements are the \emph{strongly connected components}. 
If there is only one strongly connected component, i.e., 
any two points in $G$ are strongly connected, then 
$G$ is called \emph{strongly connected}. 
\end{defn}

\begin{defn}[$c$-path boundedness and $c$-cyclic monotonicity]
\label{d:max cost}
Let $G$ be a nonempty subset of $D$,
and suppose that $(x_s,y_s), (x_e,y_e)$ are points in $G$. 
We define the \emph{maximal inner variation} 
from $(x_s,y_s)$ to $(x_e,y_e)$ 
within $G$ 
by 
\begin{equation}
        \label{e:Ranti}
        F_G: G \times G \to [-\infty, +\infty]:\big((x_s, y_s),(x_e,y_e)\big) \mapsto \sup_{\left((x_k,y_k)\right)_{k=1}^{N+1}\in P^G_{(x_s,y_s)\to (x_e,y_e)}} \sum_{k=1}^N \big(c(x_k,y_k)-c(x_{k+1}, y_k)\big), 
    \end{equation}
    where we recall that 
    $P^G_{(x_s,y_s)\to (x_e,y_e)}$ is the set of paths in $G$ from $(x_s,y_s)$ to $(x_e,y_e)$. 
Note that there is a path from $(x_s,y_s)$ to $(x_e,y_e)$ 
precisely when 
$\minf < F_G((x_s,y_s),(x_e,y_e))\leq\pinf$ because 
$F_G((x_s,y_s),(x_e,y_e))=\minf$ 
$\Leftrightarrow$ there is no path from $(x_s,y_s)$ to 
$(x_e,y_e)$ using the usual convention that 
$\sup\varnothing=\minf$. 
Using the maximal inner variation, one can elegantly 
describe $c$-path boundedness, which was recently introduced 
in \cite[Definition~1.1]{Kasia} and the classical $c$-cyclic 
monotonicity \cite{Rocky66}:
\begin{align}
\text{$G$ is \emph{$c$-path bounded}}
\;\; &:\Leftrightarrow \;\;
(\forall (x_s,y_s)\in G)(\forall(x_e,y_e) \in G)
\;\;
F_G\big((x_s,y_s),(x_e,y_e)\big) <+\infty; 
\label{e:c-pathbound} \\
\label{e:cyclical-mono}
\text{$G$ is \emph{$c$-cyclically monotone}}
\;\; &:\Leftrightarrow \;\;
(\forall (x,y)\in G)
\;\;
F_G\big((x,y),(x,y)\big)\leq 0. 
\end{align}
Note that if $G=\varnothing$, then it is 
$c$-path bounded and $c$-cyclically monotone 
(because $F_G\equiv\minf$).
 Finally, we obviously have 
\begin{equation}
\label{r:inclusion}
G\subseteq H \subseteq D \;\;\Rightarrow\;\;
F_G\leq F_H. 
\end{equation}
\end{defn}

\begin{definition}[$c$-potentiability]
\label{d:c-potential}
Let $G$ be a nonempty subset of $D$. 
Then $G$ \emph{has a $c$-potential {\color{black}(or $G$ is $c$-potentiable)}} 
if there exists a $c$-class function $\psi\colon X\to\RRX$ such that 
$G \subseteq \gra \partial_c\psi$. 
\end{definition}


For a given nonempty subset $G$ of $D$, 
the following general 
relationships between these notions were provided 
by Artstein-Avidan et~al.: 

\begin{fact}[$c$-potentiability $\Rightarrow$ $c$-path boundedness]
\label{f:0707a}
{\color{black} Let $G$ be a nonempty subset of $D$ that is $c$-potentiable. Then $G$ is $c$-path bounded.}
\end{fact}
\begin{proof}
See \cite[Proof of Theorem~1 on page~16]{Kasia}.
\end{proof}

\begin{fact}[$c$-path boundedness $\Rightarrow$ $c$-cyclic monotonicity]
\label{f:0707b}
Let $G$ be a nonempty subset of $D$ that is $c$-path bounded. 
Then $G$ is $c$-cyclically monotone. 
\end{fact}
\begin{proof}
This was observed on \cite[page~11]{Kasia}. 
\end{proof}

The central question is whether or not the two implications 
of \cref{f:0707a} and \cref{f:0707b} are actually 
equivalences. 
In the context of classical convex analysis 
and of traditional cost functions, 
the answer is affirmative: 

\begin{remark}[Rockafellar-Rochet-R\"uschendorf] 
Suppose that we are in the classical convex analysis 
setting ($Y=X^*$ and $c=-\scal{\cdot}{\cdot}$, see \cref{r:classic}). 
Rockafellar's well-known antiderivative construction (see \cite{Rocky66}) 
shows that for a nonempty subset $G$ of $D=X\times X^*$, we have 
the equivalences 
\begin{equation}
\text{$G$ is potentiable}
\;\;\Leftrightarrow\;\;
\text{$G$ is path-bounded}
\;\;\Leftrightarrow\;\;
\text{$G$ is cyclically monotone.}
\end{equation}
More generally, if $X$ and $Y$ are general sets and $c$ is a 
\emph{traditional cost}, i.e., $D=\dom c= X\times Y$, then 
variations of Rockafellar's construction by 
Rochet \cite{Rochet} and by R\"uschendorf \cite{Rueschendorf} 
give exactly the same result:
\begin{equation}
\text{$G$ is $c$-potentiable}
\;\;\Leftrightarrow\;\;
\text{$G$ is $c$-path bounded}
\;\;\Leftrightarrow\;\;
\text{$G$ is $c$-cyclically monotone.}
\end{equation}
\end{remark}

Unlike the traditional cost case, it is not possible to 
obtain $c$-potentiability from $c$-cyclic monotonicity 
as Artstein-Avidan et al.\ recently showed: 

\begin{example}[$c$-path boundedness $\not\hspace{-1mm}\Leftarrow$ 
$c$-cyclic monotonicity for the polar cost] 
\label{ex:polar}
Suppose that $X=Y=\mathbb{R}$ and that 
$c$ is the \emph{polar cost}, i.e., 
\begin{equation}
c(x,y) := \begin{cases}
-\ln(xy-1), &\text{if $xy>1$;}\\
\pinf, &\text{if $xy\leq 1$.}
\end{cases}
\end{equation}
Then $D = \menge{(x,y)\in X\times Y}{xy>1}$ and 
\begin{equation}
G := \menge{(x,3-2x)}{\tfrac{3}{4}\leq x <1}
\cup \big\{(\tfrac{3}{2},\tfrac{3}{4}) \big\}
\end{equation}
is a nonempty 
subset of $D$ that is
$c$-cyclically monotone but $G$ is not $c$-path bounded and thus 
not $c$-potentiable. 
\end{example}
\begin{proof}
See \cite[page~10]{Kasia}. 
\end{proof}

\begin{remark}[the two fundamental questions]
\cref{f:0707a}, \cref{f:0707b}, and \cref{ex:polar} 
raise the following two fundamental questions:
{\rm 
\begin{enumerate}
\item[\textbf{Q1.}] 
What assumptions guarantee that 
$c$-potentiability $\Leftarrow$ $c$-path boundedness?\\[-3mm]
\item[\textbf{Q2.}] 
What assumptions guarantee
that $c$-path boundedness $\Leftarrow$ $c$-cyclic monotonicity?
\end{enumerate}
}
\end{remark}

The state-of-the-art for \textbf{Q1} is the following very recent result 
by Artstein-Avidan et al.:

\begin{fact}[Artstein-Avidan et al., 2025]
\label{f:fix}
Let $G$ be a nonempty subset of $D$ such that $G$ 
is $c$-path bounded. 
Suppose that one of the following conditions holds:
\begin{enumerate}
\item $G$ is countable.
\item $G$ is uncountable and there is \emph{no infinite black hole}, i.e., 
for every nonempty proper infinite subset $G_0$ of $G$, there exists 
$(x_0,y_0)\in G_0$ and $(x_1,y_1)\in G\smallsetminus G_0$ 
such that  $c(x_1,y_0)<\pinf$. 
\end{enumerate}
Then $G$ is $c$-potentiable. 
\end{fact}
\begin{proof}
See \cite[Theorem~1 (Corrected) on page~2]{Kasiafix}. 
\end{proof}

\begin{remark}
\label{r:fix}
Let $G$ be a nonempty 
strongly connected\footnote{Recall \cref{d:path}.} 
subset of $D$. 
Then $G$ has no infinite black hole\footnote{
Indeed, suppose to the contrary that there is an infinite black hole. 
Then there exists a nonempty proper infinite subset $G_0$ of $G$ 
such that for all $(x_0,y_0)\in G_0$ and all $(x_1,y_1)\in G\smallsetminus G_0$, 
we have $c(x_1,y_0)=\pinf$. 
Letting $(x_s,y_s)\in G_0$ and 
$(x_e,y_e)\in G\smallsetminus G_0$, 
we see that there is no path in $G$ from $(x_s,y_s)$ 
to $(x_e,y_e)$ because at some point we have to switch
from $G_0$ to $G\smallsetminus G_0$ and then obtain a contradiction.
}. 
\end{remark}

One now obtains the following pleasant partial answer to \textbf{Q1}:

\begin{corollary}
\label{c:fix}
Let $G$ be a nonempty strongly connected 
subset of $D$. 
Then 
$G$ is $c$-potentiable
$\Leftrightarrow$
$G$ is $c$-path bounded.
\end{corollary}
\begin{proof}
``$\Rightarrow$'': \cref{f:0707a}.
``$\Leftarrow$'': If $G$ is countable, then $G$ is potentiable 
by \cref{f:fix}.(i). Otherwise, $G$ is uncountable 
and without infinite black holes (see \cref{r:fix}) and 
so \cref{f:fix}.(ii) applies.
\end{proof}

{\color{black} One might wonder if the answer for \textbf{Q1} is ``None", i.e., no extra condition is needed to guarantee that $c$-potentiability $\Leftarrow$ $c$-path boundedness. 
\cref{ex:SUPER} presented in \cref{s:important counter} below will illustrate
that \emph{some} additional condition is required.}

Let us turn to \textbf{Q2}. It follows from \cref{ex:polar} that 
\emph{some} assumptions are needed. 
The following positive results are known:

\begin{fact}
\label{f:0707c}
Let $G$ be a nonempty strongly connected subset of $D$ that is $c$-cyclically monotone. 
Then $G$ is $c$-path bounded and $c$-potentiable.
\end{fact}
\begin{proof}
See \cite[page~17]{Kasia} for the proof of 
$c$-path boundedness. 
Combined with \cref{c:fix}, we obtain $c$-potentiability.
See also \cite[Proposition~3.2]{Beigl} for related 
work when $X$ and $Y$ are Polish spaces.
\end{proof}

\begin{fact}
\label{f:0707d}
Suppose that $X$ and $Y$ are separable metric spaces 
and that the cost function $c: X\times Y\to \left]-\infty, +\infty\right]$ is continuous \footnote{{\color{black} The continuity of $c$ is in the sense of \cref{d:continuous} with $G:=X\times Y$.}}.
Let $G$ be a nonempty subset of $D$ that is $c$-cyclically monotone. 
Suppose furthermore that 
$G$ is path connected in the sense of topology and that 
$\sup c(G)<\pinf$. 
Then $G$ is $c$-path bounded. 
\end{fact}
\begin{proof}
See \cite[Corollary~4.1]{Kasia}. 
\end{proof}

\begin{fact}
Suppose that $X$ and $Y$ are {\color{black} separable}
metric spaces, 
and that the cost function $c: X\times Y\to \left]-\infty, +\infty\right]$ is continuous.
{\color{black} Let $G$ be a nonempty subset 
of $D$ that is $c$-cyclically monotone. Suppose furthermore that $\overline{G}$ is compact in $X\times Y$ and $\sup c(G)<+\infty$.} 
Then $G$ is $c$-path bounded. 
\end{fact}
\begin{proof}
See\footnote{The hypothesis in \cite{Kasia} requires that $X$ and $Y$ are separable metric spaces, $G$ is bounded, and $\sup c(G)<\pinf$. However, 
{\color{black} these assumptions (separability, boundedness and $\sup c(G)<+\infty$)} do not guarantee the compactness of $\overline{G}$ which is required in their proof.} \cite[Proposition~4.2]{Kasia}. 
\end{proof}

We are now ready to highlight the main results of this paper.

\begin{enumerate}
\item [\textbf{R1:}]
\textit{Let $G$ be a nonempty subset of $D$ that is $c$-path bounded. 
Suppose that at least one of the following holds:} 
\begin{enumerate}
\item \textit{There is a nonempty finite subset $\Omega$ of $G$
such that for every $(x,y)\in G$, there exists a path in $G$ from 
$(x,y)$ to some point $(x_e,y_e)\in \Omega$ \footnote{Recall \cref{d:path}.};}

\item \textit{There is a nonempty finite subset $A$ of $G$ such that for every $(x,y)\in G$, there exists a point $(x_s,y_s)\in A$ and a path in $G$ from $(x_s,y_s)$ to $(x,y)$ \footnote{Recall \cref{d:path}.}.}
\end{enumerate}
\textit{Then $G$ has a {\color{black} $c$-potential}.}
\end{enumerate}

The result \textbf{R1} (see \cref{t:exist1} and \cref{t:exist2} below) applies in specific scenarios with an infinite black hole 
where previous results are not applicable (see \cref{r:Coulomb}).  

\begin{enumerate}
\item [\textbf{R2:}]
\textit{Suppose that $X$ and $Y$ are separable metric spaces, and 
that the domain $D$ of the cost function $c$ is open and that $c|_D$ is continuous. 
Let $G$ be a nonempty subset of $D$. 
Then: 
\begin{align*}
{\color{black} \text{$G$ has an upper semicontinuous $c$-potential}}
\Leftrightarrow \text{$G$ is $c$-path bounded.}
\end{align*}
}
\end{enumerate}

The result \textbf{R2} (see \cref{t:continuous} below) establishes the equivalence between $c$-potentiability and $c$-path boundedness under mild topological assumptions satisfied in a great number of scenarios in optimal transport.

\begin{enumerate}
\item [\textbf{R3:}]
\textit{Suppose that $X=Y=\mathbb{R}$, and that $c$ and $G$ satisfy some regularity assumptions (for precise statements, see \cref{t:cpathextension}). Then $G$ has a strongly connected $c$-path bounded extension and thus has a {\color{black} $c$-potential}.}
\end{enumerate}

The result \textbf{R3} provides another way to show $c$-potentiability by the 
new technique of constructing a $c$-path bounded extension. 
\textbf{R3} is applicable in specific contexts where existing results fall short, as illustrated in \cref{ex:5.1}.

The remainder of the paper is organized as follows. 

In \cref{s:aux}, we present several auxiliary results used in later proofs 
{\color{black} 
and an example where $c$-path boundedness is not sufficient 
for $c$-potentiability.
}
We investigate the existence of {\color{black} $c$-potentials} 
for costs defined on the product set $X \times Y$ without assuming any topology 
in \cref{s:general}.
In \cref{s:metric}, we study the metric case  and we introduce the notion of ball-chain connectedness, which is a less restrictive condition than topological connectedness. It allows us to show that $c$-cyclical monotonicity together with topological connectedness ensures the existence of a $c$-antiderivative and hence $c$-potentiability (\cref{t:0713c}). Furthermore, we prove that, for costs which are continuous on their open domains in separable metric spaces, {\color{black} $c$-path boundedness is equivalent to the existence of an upper semicontinuous $c$-potential (\cref{t:continuous})}.
Finally, in \cref{s:extension}, we introduce the notion of a $c$-path bounded 
extension.  It is used to prove $c$-potentiability for certain costs on $\mathbb{R} \times \mathbb{R}$ for which none of the preceding results are applicable.

\section{Auxiliary results}

\label{s:aux}

\subsection{$c$-antiderivatives vs $c$-potentials}

Potentiability is more conveniently understood 
using the following ideas developed 
in the work of Artstein-Avidan et al.\ 
(see \cite[Section~3]{Kasia}) which we summarize in this section. 

\begin{defn}[$c$-antiderivative]
\label{d:c-anti}
Let $G$ be a nonempty subset of $D$.  
We say $f\colon  G \to \RR$ 
is a \emph{$c$-antiderivative} of $G$ if 
\begin{equation}
\label{e:canti}
\big(\forall(x,y)\in G\big)
\big(\forall(u,v)\in G\big)
\quad 
c(x,y)-c(u,y)\leq f(x,y)-f(u,v). 
\end{equation}
\end{defn}

\begin{remark}
\label{r:easyRemark}
From \cref{d:c-anti}, we know that if $f$ is a $c$-antiderivative of a nonempty set $G$, then $f|_\Omega$ is a $c$-antiderivative of any nonempty subset $\Omega$ of $G$. 
\end{remark}

\begin{remark}
\label{r:keyremark}
Let $G$ be a nonempty subset of $D$,
and let $f\colon G\to\RR$ be a 
$c$-antiderivative of $G$. 
If $(x,y_1),(x,y_2)$ belong to $G$, 
then $f(x,y_1)=f(x,y_2)$ 
(see \cite[Remark~3.2]{Kasia}). 
Setting $P_X(G) := \menge{x\in X}{\{x\}\times Y \cap G\neq\varnothing}$ and letting $y_x$ be any selection 
such that $(x,y_x)\in G$ for $x\in P_X(G)$, we see 
that the function  
\begin{equation}
f_G \colon P_X(G)\to\RR
\colon x\mapsto f(x,y_x)
\end{equation}
is well-defined, 
and it satisfies
\begin{equation}
\label{e:0710a}
\big(\forall(x,y)\in G\big)
\big(\forall(u,v)\in G\big)
\quad 
c(x,y)-c(u,y)\leq f_G(x)-f_G(u). 
\end{equation}
\end{remark}

These notions lead to useful characterizations of 
$c$-potentiability:

\begin{fact}[characterization of $c$-potentiability by systems of inequalities]
\label{f:et}
Let $G$ be a nonempty subset of $D$. 
Then the following are equivalent:
\begin{enumerate}
\item $G$ has a $c$-potential.
\item There exists a function $\psi\colon P_X(G)\to\RR$ such that 
\begin{equation}
\label{e:0710b}
\big(\forall(x,y)\in G\big)
\big(\forall(u,v)\in G\big)
\quad 
c(x,y)-c(u,y)\leq \psi(x)-\psi(u). 
\end{equation}
\end{enumerate}
If this happens, then 
\begin{equation}
\label{e:inf}
\Psi : X \to [-\infty, +\infty] \colon x\mapsto \inf_{(u,v)\in G} \big(c(x,v)-c(u,v)+\psi(u)\big) 
\end{equation}
is a $c$-potential for $G$ and $\Psi$
extends $\psi$, i.e., $\Psi|_{P_X(G)}=\psi$.
\end{fact}
\begin{proof}
See \cite[Theorem~3.1]{Kasia} and
\cite[Theorem~3.2]{Bartz} for the case when 
$D=X\times Y$. See also \cite[Theorem~3.5]{Bartz2} in relation to optimal transport. 
\end{proof}

\begin{remark}[$c$-antiderivative vs $c$-potential]
\label{r:c-anti vs c-potent}
Let $G$ be a nonempty subset of $D$. 
\begin{enumerate}
\item If $G$ has a $c$-antiderivative $f$, 
then $f_G$ (see \cref{r:keyremark}) is a solution 
to \cref{e:0710b}; hence, by \cref{f:et}, there exists 
a $c$-potential for $G$.
\item Conversely, if $G$ has a $c$-potential, then 
\cref{f:et} guarantees the existence of a function $\psi\colon P_X(G)\to\RR$ 
such that \cref{e:0710b} holds; consequently, the function 
$f\colon G\to\RR\colon (x,y)\mapsto \psi(x)$
is a $c$-antiderivative of $G$. 
\end{enumerate}
To sum up, 
\begin{equation}
\label{e:opera}
\text{$G$ has a $c$-antiderivative}
\;\Leftrightarrow\;
\text{$G$ has a $c$-potential.}
\end{equation}
\end{remark}
Because of equivalence \cref{e:opera}, for the rest of the paper, 
\begin{empheq}[box=\mybluebox]{equation}
\text{we will say that a set ``has a $c$-potential" or ``has a $c$-antiderivative" interchangeably.}
\end{empheq}

\subsection{Paths and maximal inner variations}

We started introducing paths in \cref{d:path}.
In this section, we record additional relevant and related results.

\begin{defn}[types of connectedness]
\label{d:connected}
Let $G$ be a nonempty subset of $D$. 
Then\footnote{All paths in this definition lie in $G$.}
$(x_s,y_s),(x_e,y_e)$ in $G$ are \ldots
\begin{enumerate}
\item \emph{strongly connected} or \emph{mutually reachable}, if there is a path from $(x_s,y_s)$ to $(x_e,y_e)$ and a path from $(x_e,y_e)$ to $(x_s,y_s)$;
\item
\emph{unilaterally connected}, if there is either a path from $(x_s, y_s)$ to $(x_e,y_e)$ or a path from $(x_e,y_e)$ to $(x_s,y_s)$ but not both;
\item
\emph{disconnected}, if there is no path from $(x_s,y_s)$ to $(x_e,y_e)$ and no path from $(x_e,y_e)$ to $(x_s,y_s)$.
\end{enumerate}
As pointed out in \cref{d:path}, strong connectedness
induces an equivalence relation which gives rise to the condensation graph and its strongly connected components. 
\end{defn}

The next result resembles a ``reverse triangle inequality'' for the maximal 
inner variation:

\begin{proposition}
\label{p:connectedanti}    
Let $G$ be a nonempty subset of $D$ that is 
$c$-path bounded. 
Then for all $(x_s,y_s),(x_m,y_m),(x_e,y_e)$ in $G$, we have
\begin{equation}
\label{e:keyineq}
     F_G\big((x_s,y_s), (x_m, y_m)\big)+F_G\big((x_m,y_m),(x_e,y_e)\big)\leq F_G\big((x_s,y_s),(x_e,y_e)\big).
\end{equation}
\end{proposition}
\begin{proof}
By $c$-path boundedness of $G$, we have 
$F_G((x_s,y_s),(x_m,y_m))<+\infty$ and 
$F_G((x_m,y_m),(x_e,y_e))<+\infty$. 
It is thus enough to consider the case  
when $\minf < F_G((x_s,y_s),(x_m,y_m))$ 
and $\minf < F_G((x_m,y_m),(x_e,y_e))$.
Fix $\varepsilon>0$, and take two paths 
$(u_i, v_i)_{i=1}^{M+1}$ in $P^G_{(x_s,y_s)\to (x_m,y_m)}$ and $(a_k,b_k)_{k=1}^{N+1}$ in $P^G_{(x_m,y_m)\to (x_e,y_e)}$ such that
\begin{subequations}
\label{e:keyineq1}
    \begin{align}
    F_G\big((x_s,y_s),(x_m, y_m)\big)-\frac{\varepsilon}{2}
    &\leq 
    \sum_{i=1}^M \big(c(u_i,v_i)-c(u_{i+1},v_i)\big), \\ 
        F_G\big((x_m,y_m),(x_e, y_e)\big)-\frac{\varepsilon}{2}
        &\leq  \sum_{k=1}^N \big(c(a_k,b_k)-c(a_{k+1},b_k)\big).  
    \end{align}
\end{subequations}
    Notice that $$\big((u_1,v_1)=(x_s,y_s),\ldots, (u_{M+1},v_{M+1})=(x_m,y_m)=(a_1,b_1), (a_2,b_2),\ldots, (a_{N+1},b_{N+1})=(x_e,y_e)\big)
    $$
is a path in 
    $P^G_{(x_s,y_s)\to (x_e,y_e)}$; thus, 
    \begin{equation}
    \label{e:keyinequ2}
        \sum_{i=1}^M \big(c(u_i, v_i)-c(u_{i+1},v_i)\big)
        +\sum_{k=1}^{N}\big(c(a_k,b_k)-c(a_{k+1},b_k)\big)
        \leq F_G\big((x_s,y_s),(x_e,y_e)\big).
    \end{equation}
    Combining \cref{e:keyineq1} and \cref{e:keyinequ2}, 
    we obtain 
    \begin{displaymath}
        F_G\big((x_s,y_s),(x_m,y_m)\big)+F_G\big((x_m,y_m),(x_e, y_e)\big)-\varepsilon\leq F_G\big((x_s,y_s),(x_e,y_e)\big),
    \end{displaymath}
    and this implies \cref{e:keyineq}.
\end{proof}

\begin{remark}
\label{r:connectedanti}
Let $(x_s,y_s), (x_m,y_m), (x_e,y_e)\in G$. Then 
\cref{d:max cost} implies that 
$c(x_s,y_s)-c(x_m, y_s)\leq F_G((x_s,y_s), (x_m, y_m))$ and 
$c(x_m,y_m)-c(x_e,y_m)\leq F_G((x_m,y_m),(x_e,y_e))$.
Consequently, using \cref{p:connectedanti}, we obtain 
\begin{equation}
\label{e:2.7-1}
c(x_s,y_s)-c(x_m, y_s)+F_G\big((x_m,y_m), (x_e,y_e)\big)\leq F_G\big((x_s,y_s),(x_e,y_e)\big)
\end{equation}
and 
\begin{equation}
\label{e:2.7-2}
F_G\big((x_s,y_s), (x_m,y_m)\big)+c(x_m,y_m)-c(x_e,y_m)\leq F_G\big((x_s,y_s),(x_e,y_e)\big).
\end{equation}
\end{remark}

\subsection{$c$-antiderivative for costs with disconnected domain}
\label{s:2.3}

\begin{proposition}
Suppose that there is a (possibly infinite) partition of $D$, 
say $(D_i)_{i\in I}$ such that for $i,j\in I$ with $i\neq j$, {\color{black} $(x,y)$ and $(u,v)$ are disconnected \footnote{This means there is neither a path in $D$ from $(x,y)$ to $(u,v)$ nor a path in $D$ from $(u,v)$ to $(x,y)$, i.e., $P^D_{(x,y)\to (u,v)}=P^D_{(u,v)\to (x,y)}=\varnothing$.}} for all $(x,y)\in D_i, (u,v)\in D_j$. For $i\in I$, define
\begin{align*}
c_i\colon X\times Y \to \left]-\infty, +\infty\right]\colon (x,y)\mapsto\begin{cases}
c(x,y), &\text{if $(x,y)\in D_i$;}\\
\pinf, &\text{otherwise.}
\end{cases}
\end{align*}
Let $G$ be a nonempty subset of $D$. Set $\widetilde{I}:=\{i\in I \ |\ G\cap D_i\neq \varnothing\}$.  Then 
\begin{equation}
\label{e:2.3}
\text{for every $i\in \widetilde{I}$, there is a $c_i$-antiderivative of $G\cap D_i$} \Leftrightarrow \text{there is a $c$-antiderivative of $G$.} 
\end{equation}
\end{proposition}
\begin{proof}
"$\Rightarrow$": Suppose $f_i: G\cap D_i \to \mathbb{R}$ is a $c_i$-antiderivative of $G\cap D_i$. 
Then for $(x,y)\in G=\bigcup_{i\in \widetilde{I}} (G\cap D_i)$, we set $f(x,y):=f_i(x,y)$. Let $(x,y), (u,v)\in G$.
\begin{enumerate}
\item [Case 1:] Suppose $(x,y), (u,v)\in G\cap D_i$. Then 
\begin{equation}
\label{e:2.3e1}
f(x,y)-f(u,v)=f_i(x,y)-f_i(u,v)\geq c_i (x,y)-c_i(u,y)=c(x,y)-c_i(u,y).
\end{equation}
Note that $(u,y)\notin D_j$ for $j\neq i$ otherwise $\big((x,y), (u,y)\big)$ is a path from $(x,y)\in D_i$ to $(u,y)\in D_j$, which contradicts the assumption of $D_i, D_j$. Hence, $(u,y)\in D_i$ or $(u,y)\notin D$. If $(u,y)\in D_i$, then $c_i(u,y)=c(u,y)$ and the right-hand side of \cref{e:2.3e1} is $c(x,y)-c(u,y)$, so $f$ satisfies \cref{e:canti}. If $(u,y)\notin D$, then $c(x,y)-c(u,y)=-\infty$ and $f$ satisfies \cref{e:canti}. 

\item[Case 2:] Suppose $(x,y)\in G\cap D_i$ and $(u,v)\in G\cap D_j$ where $i\neq j$. Then $(u,y)\notin D$ (i.e., $c(u,y)=+\infty$) by the assumption of $D_i,D_j$. Hence,
\begin{align*}
f(x,y)-f(u,v)=f_i(x,y)-f_j(u,v)>c(x,y)-c(u,y)=-\infty.
\end{align*}
\end{enumerate}

\par "$\Leftarrow$": Suppose $f: G\to \mathbb{R}$ is a $c$-antiderivative of $G$. Then for $i\in \widetilde{I}$, and $(x,y), (u,v)\in G\cap D_i$
\begin{equation}
\label{e:2.3e2}
f|_{G\cap D_i}(x,y)-f|_{G\cap D_i} (u,v)=f(x,y)-f(u,v)\geq c(x,y)-c(u,y).
\end{equation}
Note that $(u,y)\notin D_j$ for $j\neq i$, because otherwise we contradict the disconnectedness assumption for $D_i,D_j$. Hence, $(u,y)\in D_i$ or $(u,y)\notin D$, and $c(u,y)=c_i(u,y)$. Moreover, $c(x,y)=c_i(x,y)$ because $(x,y)\in D_i$. From \cref{e:2.3e2}, we have
\begin{align*}
f|_{G\cap D_i}(x,y)-f|_{G\cap D_i} (u,v)\geq c_i(x,y)-c_i(u,y).
\end{align*}
It follows that $f|_{G\cap D_i}$ is a $c_i$-antiderivative of $G\cap D_i$. 
\end{proof}

\subsection{A simple direct proof of \cref{f:fix}.(i)}
\label{sec:countable}

It's worthwhile to record a direct proof of \cref{f:fix}.(i), 
which was proved in \cite[Appendix]{Kasia} with  
more involved arguments 
using tools from graph theory.

The proof below does follow the original Zorn's-Lemma approach 
in the (wrong-in-general) \cite[Proof of Theorem~3.3]{Kasia}.
We note that the countability assumption is essential
in one step of the proof as the ingenious example 
presented in \cite[Section~2]{Kasiafix} shows.

\begin{proposition}
    \label{p:countableset}
Let $G$ be a nonempty subset of $D$ such that $G$ is $c$-path bounded 
and $G$ is \emph{countable}. Then $G$ has a $c$-antiderivative\footnote{The proof actually offers an inductive 
construction of $f$.}.
\end{proposition}
\begin{proof}
By assumption, there is a sequence $(x_n,y_n)_\nnn$ in $G$ 
such that $G = \{(x_n,y_n)\}_\nnn$. 
{\color{black} We assume first that $G$ is infinite and 
WLOG that $n\mapsto (x_n,y_n)$ is injective.}
For $n\in\NN$, we abbreviate $I_n := \NN\cap[0,n]$, 
$z_n := (x_n,y_n)$, and $G_n := \{z_0,\ldots,z_n\}$.
Using induction, 
we shall prove that there exists a sequence of functions $(f_n)_\nnn$
such that $(\forall\nnn)$ $f_n\colon G_n\to\RR$, 
$f_{n+1}|_{G_n}=f_n$, and
\begin{equation}
\label{e:0713d}
(\forall i\in I_n)(\forall j\in I_n)
\quad 
F_G(z_i,z_j)\leq f_n(z_i)-f_n(z_j),
\end{equation}
where the maximal inner variation $F_G$ is defined in \cref{e:Ranti}. 

For the base case, we define $f_0\colon G_0=\{z_0\}\to\RR\colon 
z_0\mapsto 0$ which clearly satisfies \cref{e:0713d} for $n=0$.

Now we tackle the inductive step and assume that \cref{e:0713d}
holds for $n\in\NN$. 
Because $G$ is $c$-path bounded, we have $(\forall j\in I_n)$  
$F_G(z_{n+1},z_j)<\pinf$ and $F_G(z_j,z_{n+1})<\pinf$. 
Using \cref{p:connectedanti}, we have 
$(\forall i\in I_n)(\forall j\in I_n)$ 
$F_G(z_j,z_{n+1})+F_G(z_{n+1},z_i)\leq F_G(z_j,z_i)$; 
moreover, \cref{e:0713d} yields
$F_G(z_j,z_i) \leq f_n(z_j)-f_n(z_i)$. 
Altogether, 
\begin{equation}
(\forall i\in I_n)(\forall j\in I_n)\quad 
F_G(z_{n+1}, z_i)+f_n(z_i)\leq f_n(z_j)-F_G(z_j,z_{n+1}). 
\end{equation}
Hence 
\begin{align}
\label{e:0713e}
\minf \leq \alpha &:= 
\max_{i\in I_n} \big(F_G(z_{n+1}, z_i)+f_n(z_i)\big) 
\leq \min_{j\in I_n} \big (f_n(z_j)-F_G(z_j,z_{n+1})\big) =:\beta \leq \pinf.
\end{align}
Note that $\alpha<\pinf$ and $\minf < \beta$. 
Pick  $\gamma\in [\alpha,\beta]\cap\RR\neq\varnothing$, 
and define $f_{n+1}$ by
\begin{equation}
f_{n+1}\colon G_{n+1}\to\RR\colon 
z_i \mapsto \begin{cases}
f_n(z_i), &\text{if $i\in I_n$;}\\
\gamma, &\text{if $i=n+1$.}
\end{cases}
\end{equation}
It's clear that $f_{n+1}$ is {\color{black} finite-valued} and that it extends $f_n$.
It remains to verify 
\begin{equation}
\label{e:0713f}
(\forall i\in I_{n+1})(\forall j\in I_{n+1})
\quad 
F_G(z_i,z_j)\leq f_{n+1}(z_i)-f_{n+1}(z_j).
\end{equation}
By the extension property,
these are clear when both $i$ and $j$ are in $I_n$
(by \cref{e:0713d}), or
when $i=j=n+1$ (because $G$ is $c$-path bounded, hence $c$-cyclically 
monotone). 
If $i\in I_n$ and $j=n+1$, 
$f_{n+1}(z_{n+1}) = \gamma
\leq \beta \leq f_n(z_i)-F_G(z_i,z_{n+1})$ and 
thus $F_G(z_i,z_j)\leq f_{n+1}(z_i)-f_{n+1}(z_j)$ as desired.
And if $i=n+1$ and $j\in I_n$, then 
$f_{n+1}(z_{n+1}) = \gamma
\geq \alpha \geq F_G(z_{n+1},z_j)+f_n(z_j)$ and so 
$F_G(z_i,z_j) \leq f_{n+1}(z_i)-f_{n+1}(z_j)$ as well. 
This completes the proof by induction.

Finally, we define 
\begin{equation}
f\colon G\to \RR\colon z_n = (x_n,y_n)\mapsto f_n(z_n). 
\end{equation}
Consider two points $(x_i,y_i)$ and $(x_j,y_j)$ in $G$. 
By definition of $F_G$, we have 
$c(x_i,y_i)-c(x_j,y_i) \leq F_G((x_i,y_i),(x_j,y_j)) 
= F_G(z_i,z_j)$. 
Set $n=\max\{i,j\}$. By \cref{e:0713d}, we have
on the other hand, 
$F_G(z_i,z_j)\leq f_n(z_i)-f_n(z_j) = f(z_i)-f(z_j)
= f(x_i,y_i)-f(x_j,y_j)$. 
Altogether, 
$c(x_i,y_i)-c(x_j,y_i) \leq f(x_i,y_i)-f(x_j,y_j)$, 
which verifies \cref{e:canti} and thus $f$ is an antiderivative of $G$. 
{\color{black} The case when $G$ is finite is treated similarly 
by working with a finite sequence instead.} 
\end{proof} 

{\color{black}
\subsection{A $c$-path bounded set that is not $c$-potentiable}
\label{s:important counter}
The following example shows that the implication ``$c$-potentiability $\Leftarrow$ $c$-path boundedness" does not hold without extra assumptions. The idea is from \cite[Section~2]{Kasiafix}. We simplified the proof with the help of \texttt{ChatGPT}.
\begin{example}[$c$-potentiability $\not\Leftarrow$ $c$-path boundedness]
\label{ex:SUPER}
Set $\mathcal{S}:=\RR^\mathbb{N}$, 
and suppose that 
$X=Y=\mathbb{N}\cup \mathcal{S}$ and 
\begin{align*}
c\colon X\times Y \to\left]-\infty, +\infty\right]\colon 
(x,y) \mapsto \begin{cases}
0, &\text{if $x=y$;}\\
-a_k, &\text{if $x=k\in \NN$ and $y=(a_n)_\nnn\in \mathcal{S}$;}\\
+\infty, &\text{otherwise.}
\end{cases}
\end{align*}
Furthermore, set 
\begin{equation*}
G:=\{(u,u)\ |\ u\in X=Y\}.
\end{equation*}
Then $G$ is $c$-path bounded but not $c$-potentiable. 
\end{example}
\begin{proof}
For $(x_s,y_s), (x_e,y_e)\in G$, we have
\begin{align*}
F_G\big((x_s,y_s), (x_e, y_e)\big)= \begin{cases}
0, &\text{if $(x_s,y_s)=(x_e, y_e)$;}\\
a_k, &\text{if $x_s=y_s=(a_n)_\nnn$ and $x_e=y_e=k\in \NN$;}\\
-\infty, &\text{otherwise.}
\end{cases}
\end{align*}
This implies that $G$ is $c$-path bounded. 

\par It remains to show that $G$ is not $c$-potentiable. 
Suppose to the contrary that $G$ has a $c$-potential. 
From \cref{r:c-anti vs c-potent}.(ii), we know that $G$ has a $c$-antiderivative $f: G\to \mathbb{R}$. In particular, 
\begin{equation}
\label{e:section2.5e1}
(\forall (x,x)\in \mathcal{S}\times \mathcal{S})(\forall (n,n)\in \NN\times \NN)\quad c(x,x)-c(n,x)\leq f(x,x)-f(n,n). 
\end{equation}
Now set 
\begin{align*}
(\forall n\in\NN)\quad
\widetilde{a}_n:=n-f(n,n).
\end{align*}
Then $(\widetilde{a}_n)_\nnn\in \mathcal{S}$. Applying \cref{e:section2.5e1} with $x:=(\widetilde{a}_n)_\nnn$, we have
\begin{align*}
(\forall n\in \NN) \quad c(x,x)-c(n,x)=\widetilde{a}_n=n-f(n,n)\leq f(x,x)-f(n,n).
\end{align*}
This implies $f(x,x)\geq n$ for every $n\in \NN$, i.e., $f(x,x)=+\infty$, which contradicts the fact that $f$ is finite-valued. 
\end{proof}
}

\section{$c$-antiderivatives in the general setting}

\label{s:general}

The following simple yet powerful result gives a recipe for 
constructing $c$-antiderivatives.

\begin{lemma}
\label{l:lemma1}
Let $G$ be a nonempty subset of $D$ that is $c$-path bounded. 
Let $\Omega$ be a nonempty subset of $G$, and define\footnote{Recall \cref{d:max cost}.} 
    \begin{align}
        \label{e:sup}
f\colon G\to\RRX\colon (x,y) &\mapsto     \sup_{(x_e,y_e)\in \Omega} F_G\big((x,y), (x_e,y_e)\big).
    \end{align}
If $f$ is {\color{black} finite-valued}, then 
$f$ is a $c$-antiderivative of $G$.
\end{lemma}
\begin{proof}
It is enough to check \cref{e:canti}.
To this end, 
let $(x,y), (u,v)$ be in $G$ and assume that $c(x,y)-c(u,y)>-\infty$ 
(otherwise, the inequality \cref{e:canti} trivially holds). 
Applying \cref{e:2.7-1} in \cref{r:connectedanti} with $(x_s,y_s)=(x,y)$, $(x_m,y_m)=(u,v)$ and $(x_e,y_e)\in \Omega$, 
we have 
\begin{equation}
\label{e:0711a}
          c(x,y)-c(u,y)+F_G\big((u,v),(x_e,y_e)\big)\leq F_G\big((x,y),(x_e,y_e)\big).
\end{equation}
Taking the supremum 
over $(x_e,y_e)\in \Omega$ in \cref{e:0711a} yields 
      \begin{align*}
          &\ \ \quad c(x,y)-c(u,y)+\sup_{(x_e,y_e)\in \Omega} F_G\big((u,v),(x_e,y_e)\big)\leq \sup_{(x_e,y_e)\in \Omega} F_G\big((x,y),(x_e,y_e)\big), 
          \end{align*}
          i.e., 
$c(x,y)-c(u,y)+f(u,v)\leq f(x,y)$. 
Because $f$ is {\color{black} finite-valued}, this is equivalent to 
$c(x,y)-c(u,y)\leq f(x,y)-f(u,v)$ and we are done. 
\end{proof}

We also have the following counterpart to \cref{l:lemma1}: 

\begin{lemma}
\label{l:lemma2}
Let $G$ be a nonempty subset of $D$ that is $c$-path bounded. 
Let $A$ be a nonempty subset of $G$, and define\footnote{Recall \cref{d:max cost}.} 
    \begin{align}
        \label{e:sup2}
f\colon G\to\RRX\colon (x,y) &\mapsto     -\sup_{(x_s,y_s)\in A} F_G\big((x_s,y_s), (x,y)\big). 
    \end{align}
If $f$ is {\color{black} finite-valued}, then 
$f$ is a $c$-antiderivative of $G$.
\end{lemma}
\begin{proof}
The proof is similar to that of \cref{l:lemma1}: indeed, 
the counterpart to \cref{e:0711a}, which follows from \cref{e:2.7-2} in 
\cref{r:connectedanti} with $(x_m,y_m)=(x,y)$ and $(x_e,y_e)=(u,v)$, is 
\begin{equation}
\label{e:0711c}
          F_G\big((x_s,y_s),(x,y)\big)+c(x,y)-c(u,y)\leq F_G\big((x_s,y_s),(u,v)\big).
\end{equation}
The rest of the proof is analogous. 
\end{proof}

\begin{theorem}
\label{t:exist1}
Let $G$ be a nonempty subset of $D$ that is $c$-path bounded. 
Suppose that there is a nonempty finite subset $\Omega$ of $G$
such that for every $(x,y)\in G$, there exists a path from 
$(x,y)$ to some point $(x_e,y_e)\in \Omega$. 
Then the function $f$ from 
\cref{e:sup} is a $c$-antiderivative of $G$. 
\end{theorem}
\begin{proof}
Let $(x,y)\in G$. 
In view of \cref{l:lemma1}, it suffices to show that 
$f(x,y)\in\RR$. 
Because of $c$-path boundedness of $G$, we know that 
$F_G((x,y),(x_e,y_e))<+\infty$ for all $(x_e,y_e)\in \Omega$. 
On the other hand, the assumption implies that there is 
some $(x_e,y_e)\in \Omega$ that can be reached through a path in $G$ 
starting at $(x,y)$; thus, 
$\minf < F_G((x,y),(x_e,y_e))$. 
Altogether, because $\Omega$ is finite, we deduce that  
$f(x,y)=\max_{(x_e,y_e)\in \Omega} F_G((x,y),(x_e,y_e)) \in\RR$.
\end{proof}

\begin{theorem}
\label{t:exist2}
Let $G$ be a nonempty subset of $D$ that is $c$-path bounded. 
Suppose that there is a nonempty finite subset $A$ of $G$
such that for every $(x,y)\in G$, there exists a point 
$(x_s,y_s)\in A$ and a path from 
$(x_s,y_s)$ to $(x,y)$. 
Then the function $f$ from 
\cref{e:sup2} is a $c$-antiderivative of $G$. 
\end{theorem}
\begin{proof}
The proof is analogous to that of \cref{t:exist1}.
\end{proof}

\begin{corollary}
\label{c:finiteequi}
Let $G$ be a nonempty subset of $D$ and assume $G$ is $c$-path bounded. 
Suppose that the condensation graph\footnote{See \cref{d:path}.} of $G$  
has only finitely many strongly connected components. 
Then $G$ has a $c$-antiderivative.
\end{corollary}
\begin{proof}
Denote the finitely many strongly connected components of $G$
by $C_1,\ldots,C_M$. Pick representatives 
$(x_i,y_i)$ for each $C_i$. 
If $(x,y)\in G$, then $(x,y)$ belongs to exactly one $C_i$ and 
thus there is a path from $(x,y)$ to $(x_i,y_i)$. 
Now apply \cref{t:exist1} with 
{\color{black}
$\Omega:=\{(x_1,y_1),\ldots,
(x_M,y_M)\}$. 
}
\end{proof}

\begin{corollary}
\label{c:3.5}
Let $G$ be a nonempty subset of $D$ and assume $G$ is $c$-path bounded. Suppose that there are $A\subseteq X$ and $B\subseteq Y$ such that $G\subseteq A\times B \subseteq D$. Then $G$ has only one strongly connected component and has a $c$-antiderivative. 
\end{corollary}
\begin{proof}
Let $(x,y), (u,v)\in G$. Then $x,u\in A$ and $y,v\in B$, hence $(x,v), (u,y)\in A\times B \subseteq D$. This means $\big((x,y), (u,v)\big)$ is a path from $(x,y)$ to $(u,v)$ and $\big((u,v), (x,y)\big)$ is a path from $(u,v)$ to $(x,y)$. Therefore, $(x,y), (u,v)$ are strongly connected. Since $(x,y), (u,v)$ are arbitrary in $G$, the condensation graph of $G$ has only one strongly connected component. By \cref{c:finiteequi}, we conclude that $G$ has a $c$-antiderivative.   
\end{proof}

As an application of \cref{c:3.5}, we prove the existence of $c$-potential of the Bregman cost that was studied in \cite{Kainth}.

\begin{example}[Bregman cost]
\label{ex:Bregman}
Suppose that $X=Y=\mathbb{R}^n$ and that $f:X \to \left]-\infty, +\infty\right]$ is proper with domain $\dom(f)$. Assume $f$ is differentiable on $\inte \dom(f)$. Then the \emph{Bregman cost} is given by
\begin{equation}
c\colon X\times Y \to\RX \colon 
(x,y)\mapsto 
\begin{cases}
f(x)-f(y)-\langle x-y, \nabla f(y)\rangle, &\text{if $(x,y)\in \dom(f)\times \inte \dom(f)$}\\
\pinf, &\text{otherwise,}
\end{cases}
\end{equation}
whose domain is $D=\dom(f)\times \inte \dom(f)$. Let $G$ be a nonempty subset of $D$. Combining \cref{f:0707c}, \cref{f:0707b}, \cref{c:3.5} and \cref{e:opera}, we obtain the equivalence
\begin{equation}
\text{$G$ has a $c$-potential}
\;\Leftrightarrow\;
\text{$G$ has a $c$-antiderivative}
\;\Leftrightarrow\;
\text{$G$ is $c$-path bounded}
\;\Leftrightarrow\;
\text{$G$ is $c$-cyclically monotone.}
\end{equation}
\end{example}

\begin{corollary}
\label{c:Coulomb}
Let $G$ be a nonempty subset of $D$ that is $c$-path bounded. 
Suppose that there exists $(u,v)\in G$ such that 
the set 
\begin{equation}
\label{e:0711d}
\menge{y\in Y}{\exists x\in X, \ (x,y)\in G, \ c(u,y)=\pinf}
\end{equation}
contains at most finitely many elements. 
Then $G$ has a $c$-antiderivative.
\end{corollary}
\begin{proof}
Denote the set in \cref{e:0711d} by $S$.
For every $y\in S$, pick $x_y\in X$ such that 
$(x_y,y)\in G$ and $c(u,y)=\pinf$. 
Then the set $\Omega := \{(u,v)\}\cup \{(x_y,y)\}_{y\in S}$ is finite
by assumption.
Now let $(x,y)\in G$. 
If $c(u,y)<\pinf$, then $((x,y),(u,v))$ 
is a one-step path from $(x,y)$ to $(u,v)\in \Omega$. 
Now assume that $c(u,y)=\pinf$. 
Then $y\in S$ and $((x,y),(x_y,y))$ is a one-step path 
from $(x,y)$ to $(x_y,y)\in\Omega$ because 
$(x_y,y)\in G\subseteq D$ and so 
$c(x_y,y)<\pinf$. 
The conclusion thus follows from \cref{t:exist1}.
\end{proof}

We now demonstrate how 
\cref{c:Coulomb} is applicable to the Coulomb cost, 
 which has important applications in quantum physics 
 and density {\color{black} functional} theory (see, e.g., \cite{Marino}):
 
\begin{example}[Coulomb cost]
\label{ex:Coulomb}
Suppose that $X=Y=\RR^n$ {\color{black} 
is equipped with the Euclidean norm $\|\cdot\|$} and that $c$ is the 
\emph{Coulomb cost}, given by
\begin{equation}
c\colon X\times Y \to\RX \colon 
(x,y)\mapsto 
\begin{cases}
\frac{1}{\|x-y\|}, &\text{if $x\neq y$;}\\
\pinf, &\text{if $x=y$,}
\end{cases}
\end{equation}
and with $D = \dom c = \menge{(x,y)\in X\times Y}{x\neq y}$. 
Let $G$ be a nonempty subset of $D$, and 
let $(u,v)\in G$.
{\color{black} If $c(u,y)=\pinf$, then necessarily $y=u$.
It follows that the set in \cref{e:0711d} contains at most one element!}
Combining \cref{c:Coulomb} with \cref{f:0707a} and \cref{e:opera}, 
we obtain the equivalence
\begin{equation}
\text{$G$ has a $c$-potential}
\;\Leftrightarrow\;
\text{$G$ has a $c$-antiderivative}
\;\Leftrightarrow\;
\text{$G$ is $c$-path bounded.}
\end{equation}
\end{example}

\begin{remark}
\label{r:Coulomb}
With the results in this section, we can show that a $c$-potential exists in some senarios where the ``no infinite black hole" \footnote{Recall \cref{f:fix}.} condition does not hold. For instance, let $X=Y=\mathbb{R}$ and $c: X\times Y \to \left]-\infty, +\infty\right]$ be the Coulomb cost in \cref{ex:Coulomb}. Moreover, set $$G:=\{(2,1)\}\cup ([3,4]\times \{2\}).$$ 
We first note that 
\begin{equation*}
\text{$[3,4]\times \{2\}$ is an infinite black hole.}
\end{equation*}
Indeed, for every $(x_0,y_0)\in [3,4]\times \{2\}$ and $(x_1, y_1)\in \{(2,1)\}$, we have $(x_1, y_0)=(2,2)$, so $c(x_1, y_0)=+\infty$. 
Furthermore, 
\begin{equation*}
\text{$G$ is $c$-path bounded,}
\end{equation*}
which we verify by cases on the starting point and endpoint of the path:
\begin{enumerate}
\item [Case 1:] Suppose $(x_s,y_s)=(x_e,y_e)\in \{(2,1)\}$ and $((x_n,y_n))_{n=1}^{N+1}\in P^G_{(x_s,y_s)\to (x_e,y_e)}$. Then $(x_n, y_n)=(x_s, y_s)=(x_e,y_e)$ for all $n\in \{1, \ldots, N+1\}$ and 
\begin{align*}
\sum_{n=1}^N \big(c(x_n, y_n)-c(x_{n+1}, y_n)\big)=0. 
\end{align*}
Hence, $F_G((x_s,y_s), (x_e,y_e))=0<+\infty$.

\item [Case 2:] Suppose $(x_s,y_s), (x_e,y_e)\in [3,4]\times \{2\}$ and $((x_n,y_n))_{n=1}^{N+1}\in P^G_{(x_s,y_s)\to (x_e,y_e)}$. Then $(x_n,y_n)\in [3,4]\times \{2\}$ for all $n\in \{1, \ldots, N+1\}$ and 
\begin{align*}
\sum_{n=1}^N \big(c(x_n, y_n)-c(x_{n+1}, y_n)\big)
=\sum_{n=1}^N \big(c(x_n, 2)-c(x_{n+1}, 2)\big)=c(x_s, 2)-c(x_e, 2)<+\infty.
\end{align*}
This means $F_G ((x_s,y_s), (x_e,y_e))=c(x_s,2)-c(x_e, 2)<+\infty$.

\item [Case 3:] Suppose $(x_s, y_s)\in \{(2,1)\}$, $(x_e,y_e)\in [3,4]\times \{2\}$ and $((x_n,y_n))_{n=1}^{N+1}\in P^G_{(x_s,y_s)\to (x_e,y_e)}$. Then there is $k\in \{1, \ldots, N\}$ such that $(x_n, y_n)\in \{(2,1)\}$ for $n\leq k$ and $(x_n, y_n)\in [3,4]\times \{2\}$ for $n>k$. Hence,
\begin{align*}
\sum_{n=1}^N \big(c(x_n, y_n)-c(x_{n+1}, y_n)\big)
&=\sum_{n=1}^k \big(c(x_n, y_n)-c(x_{n+1}, y_n)\big)
+\sum_{n=k+1}^N \big(c(x_n, y_n)-c(x_{n+1}, y_n)\big)\\
&=c(x_k, y_k)-c(x_{k+1}, y_k)
+\sum_{n=k+1}^N \big(c(x_n, y_n)-c(x_{n+1}, y_n)\big)\\
&=c(2, 1)-c(x_{k+1}, 1)+c(x_{k+1}, 2)-c(x_e, 2)\\
&=\frac{1}{2-1}-\frac{1}{x_{k+1}-1}+\frac{1}{x_{k+1}-2}-\frac{1}{x_e-2}\\
&=1-\frac{1}{x_e-2}+\frac{1}{(x_{k+1}-2)(x_{k+1}-1)}\\
&\leq \frac{3}{2}-\frac{1}{x_e-2}<+\infty.
\end{align*}
Therefore, $F_G((x_s,y_s), (x_e,y_e))<+\infty$. 

\item[Case 4:] Suppose $(x_s,y_s)\in [3,4]\times \{2\}$ and $(x_e,y_e)\in \{(2,1)\}$. Then $P^G_{(x_s,y_s)\to (x_e,y_e)}=\varnothing$ and thus $F_G((x_s,y_s), (x_e,y_e))=-\infty<+\infty$. 
\end{enumerate}
Because $G$ has an infinite black hole, \cref{f:fix} is not applicable to prove $c$-potentiability of $G$. However, as shown in \cref{ex:Coulomb}, 
where \cref{c:Coulomb} is applied, $G$ has a $c$-antiderivative, i.e., $G$ is $c$-potentiable.  
\end{remark}

We conclude this section with another example of a graph that 
has a $c$-antiderivative: 
   
\begin{example}
\label{ex:not unique}
Suppose that $X=Y=\RR$ and that $c$ is given by 
\begin{equation}
c\colon \RR\times \RR \to\RX \colon 
(x,y)\mapsto 
\begin{cases}
{-|x-y| = y-x}, &\text{if $y\leq x$;}\\
\pinf, &\text{if $x<y$.}
\end{cases}
\end{equation}
Now set {\color{black}$$G:= \RR(1,1)=\{t(1,1) \ |\ t\in \mathbb{R}\}=\{(t,t)\ |\ t\in \mathbb{R}\}.$$}
Given $(x_s,x_s),(x_e,x_e)$ in $G$, 
there is a path from $(x_s,x_s)$ to $(x_e,x_e)$
$\siff$ 
$c(x_e,x_s)<\pinf$
$\siff$ 
$x_s\leq x_e$, in which case $c(x_e,x_s)=x_s-x_e$. 
It follows that 
\begin{equation}
(\forall x_s\in \RR)(\forall x_e\in\RR)
\quad 
F_G\big((x_s,x_s),(x_e,x_e)\big)= 
\begin{cases}
|x_s-x_e|=x_e-x_s, &\text{if $x_s\leq x_e$;}\\
\minf, &\text{if $x_s>x_e$.}
\end{cases}
\end{equation}
It follows that $G$ is $c$-path bounded and 
so $c$-cyclically monotone (see \cref{f:0707b}). 
The condensation graph\footnote{Recall~\cref{d:path}.} is the same as 
$G$ and every strongly connected component is a singleton; 
in particular, there are infinitely many strongly connected components. 
Is there a subset {\color{black}$\Omega = S(1,1)=\{(t,t) \ |\ t\in S\subseteq\mathbb{R}\}$} of $G$ satisfying
the hypothesis of \cref{t:exist1}? 
To be reachable from any point in $G$, 
the set $S$ must necessarily satisfy 
$\sup S=\pinf$; consequently, it cannot be finite.
Now consider
{\color{black}$\Omega = \NN(1,1)=\{(n,n) \ |\ n\in \NN\}$}, which is reachable from all points in $G$. 
Then the function $f$ from \cref{e:sup} satisfies
$f\equiv\pinf$; in particular, it cannot be a $c$-antiderivative of $G$.
We learn that the finiteness assumption of $\Omega$ in 
\cref{t:exist1} is crucial! 
However, a function $f$ defined on $G$, written as 
$f(x,x)=\phi(x)$, where $\phi\colon\RR\to\RR$, 
is a $c$-antiderivative of $G$ if and only if $x\mapsto x+\phi(x)$ is decreasing {\color{black} (also known as ``nonincreasing'')}%
\footnote{Indeed, $f$ is a $c$-antiderivative of $G$
$\Leftrightarrow$ $(\forall x\in\RR)(\forall y\in\RR)$ 
$c(x,x)-c(y,x)\leq \phi(x)-\phi(y)$ $\siff$ 
$(\forall x\leq y)$ $y-x\leq \phi(x)-\phi(y)$ $\siff$
$x\mapsto x+\phi(x)$ is decreasing.}. 
So $f(x,x)=-x$ is a $c$-antiderivative of $G$. 
This illustrates that the finiteness assumptions in 
\cref{t:exist1} and \cref{c:finiteequi} are sufficient, but 
not necessary. 
\end{example}

\section{$c$-antiderivatives in the metric setting with open domain}

\label{s:metric}

In this section, we assume that our general sets 
\begin{empheq}[box=\mybluebox]{equation}
\label{e:metricXY}
\text{$X$ and $Y$ are metric spaces.
}
\end{empheq}
We equip the product space $X\times Y$ with the metric 
$d_{X\times Y}((x,y), (u,v)):=\scriptstyle\sqrt{\textstyle d^2_X(x,u)+d^2_Y(y,v)}$. 
The \emph{open ball} centered at $(x,y)$ with radius 
$\rho \in \left]0,\pinf\right]$ is denoted by 
$\mathcal{B}_{\rho} (x,y)$:
\begin{equation}
\mathcal{B}_{\rho} (x,y) := 
\menge{(u,v)\in X\times Y}{d_{X \times Y}\big((x,y), (u,v)\big)<\rho}. 
\end{equation}
We additionally assume that 
\begin{empheq}[box=\mybluebox]{equation}
\label{e:openD}
\text{$D = \dom c = \menge{(x,y)\in X\times Y}{c(x,y)<\pinf}$ is open, }
\end{empheq}
which allows us to introduce the following notion:

\begin{definition}[maximal ball]
\label{d:maxball}
Set 
\begin{equation}
\label{d:maxR}
\big(\forall (x,y)\in D\big)\quad    
\rho(x,y):=\sup\menge{\rho>0}{\mathcal{B}_{\rho}(x,y)\subseteq D} 
\;\;\text{and}\;\;
\mathcal{B}(x,y) := \mathcal{B}_{\rho(x,y)}(x,y). 
\end{equation}
Because $D$ is open, we have that $\mathcal{B}(x,y)\subseteq D$ 
for every $(x,y)\in D$. 
We call $\mathcal{B}(x,y)$ the \emph{maximal (open) ball} centered at $(x,y)$ contained in $D$. 
\end{definition}

\begin{remark}
\label{r:oneball}
Let $(x,y)\in D$, and let $(u,v)\in \mathcal{B}(x,y)$. 
Then $d((x,y),(u,v))={\scriptstyle\sqrt{\textstyle d_X^2(x,u)+d_Y^2(y,v)}} < \rho(x,y)$. 
It follows that 
$d((x,y),(u,y)) = d(x,u) <\rho(x,y)$
and
$d((x,y),(x,v))= d(y,v)<\rho(x,y)$. 
Hence $(u,y),(x,v)$ both belong to $\mathcal{B}(x,y)\subseteq D$ 
and so $c(u,y)<\pinf$ and $c(x,v)<\pinf$.
Therefore, 
$((x,y),(u,v))$ is a one-step path from $(x,y)$ to $(u,v)$ 
while $((u,v),(x,y))$ is a one-step path from $(u,v)$ to $(x,y)$. 
In turn, we deduce that 
\begin{equation}
\label{e:oneball}
\text{every two points $(u_1,v_1),(u_2,v_2)$ in 
$\mathcal{B}(x,y)$ are strongly connected.
}
\end{equation}
\end{remark}

\begin{defn}[ball chain connectedness] 
\label{d:openballs}
Let $G$ be a nonempty subset of $D$, and 
let $(x_s,y_s)$ and $(x_e, y_e)$ be two points in $G$. 
A \emph{ball chain (of length $N$) in $G$, with starting point $(x_s,y_s)$ and endpoint $(x_e,y_e)$,} 
is a finite sequence $((x_k, y_k))_{k=1}^{N+1} \in G^{N+1}$ such that 
the following hold: 
\begin{enumerate}
\item $(x_1, y_1) = (x_s,y_s)$,
\item $(x_{N+1},y_{N+1}) = (x_e, y_e)$,
\item $\mathcal{B}(x_k,y_k)\cap \mathcal{B}(x_{k+1},y_{k+1})\cap G\neq
\varnothing$ for every $k\in\{1,\ldots,N\}$. 
\end{enumerate}
This induces an equivalence relation on $G$, giving rise 
to ball chain connectedness via (maximal) open balls. 
The resulting quotient set contains the \emph{ball chain connected components}. If there is only one ball chain connected component,
then we say that $G$ is \emph{ball chain connected}. 
Note that 
\cref{r:oneball} implies that if there is a ball chain from 
$(x_s,y_s)$ to $(x_e,y_e)$, then $(x_s,y_s)$ and $(x_e,y_e)$ 
are strongly connected. Consequently,
\begin{equation}
\label{e:0713b}
\text{ball chain connectedness $\Rightarrow$ strong connectedness.}
\end{equation}
\end{defn}

We now link ball chain connectedness to classical connectedness from topology: 

\begin{proposition}[topological connectedness $\Rightarrow$ ball chain connectedness]
\label{p:0713a}
Assume that $D$ is open. Let $G$ be a nonempty subset of $D$ that is topologically connected. 
Then $G$ is ball chain connected. 
\end{proposition}
\begin{proof} 
We prove the contrapositive and assume that 
$G$ is not ball chain connected. 
Then there exists a nonempty ball chain connected component $E$ of $G$ such that $G\smallsetminus E\neq\varnothing$. 
Hence $(\cup_{(x,y)\in E} \mathcal{B}(x,y))\cap G$ 
and $(\cup_{(u,v)\in G\smallsetminus E} \mathcal{B}(u,v))\cap G$ 
are nonempty open sets in the relative topology of $G$. 
Clearly, $G\subseteq (\cup_{(x,y)\in E} \mathcal{B}(x,y))
\cup (\cup_{(u,v)\in G\smallsetminus E} \mathcal{B}(u,v))$.
If 
$G \cap (\cup_{(x,y)\in E} \mathcal{B}(x,y)) \cap 
(\cup_{(u,v)\in G\smallsetminus E} \mathcal{B}(u,v)) \neq
\varnothing$, then there exist $(x,y)\in E$ and 
$(u,v)\in G\smallsetminus E$ such that 
$\mathcal{B}(x,y)\cap\mathcal{B}(u,v)\cap G\neq \varnothing$, 
which would imply the absurdity that $(x,y)$ and $(u,v)$ 
are ball chain connected despite the fact that $(u,v)$ does not 
belong to $E$.
Hence $G \cap (\cup_{(x,y)\in E} \mathcal{B}(x,y)) \cap 
(\cup_{(u,v)\in G\smallsetminus E} \mathcal{B}(u,v)) =
\varnothing$, and therefore $G$ is not topologically connected.
\end{proof}

\begin{example}[topological connectedness\; $\not\hspace{-1mm}\Leftarrow$ ball chain connectedness]
Suppose that $X=Y=\RR$ with the usual metric. 
Now suppose that $D$ is the open unit ball, and set 
$G := \{(0,0), (0.1,0)\}$. 
Then $G$ is ball chain connected, but not 
topologically connected.
\end{example}

\begin{proposition}
\label{t:0713c}
Suppose that $c$ is a cost with open domain (i.e., $D$ is open). Let $G$ be a nonempty subset of $D$ that is topologically connected, 
or ball chain connected or strongly connected. 
If $G$ is $c$-cyclically monotone, then $G$ is $c$-path bounded, 
$c$-potentiable, and it has a $c$-antiderivative.
\end{proposition}
\begin{proof}
If $G$ is topologically connected, then 
\cref{p:0713a} implies that $G$ is ball chain connected.
And if $G$ is ball chain connected, then $G$ is strongly connected by 
\cref{e:0713b}.
Now recall \cref{f:0707c} and \cref{e:opera}. 
\end{proof}



\begin{remark}
We note that \cref{t:0713c} is a considerable 
generalization of \cref{f:0707d}: indeed, in the latter result, we assume that
$c$ is continuous (which implies that $D$ is open) as well 
as the separability of the spaces $X$ and $Y$ and $\sup c(G)<\pinf$ 
which are not needed in \cref{t:0713c}. 
\end{remark}

Furthermore, if the cost satisfies some continuity assumptions, we obtain equivalence between $c$-path boundedness and $c$-potentiability. 

\begin{defn}[continuous extended real-valued function]
\label{d:continuous}
Let $S$ be a nonempty subset of $X \times Y$. A function $h\colon S\to \left]-\infty, +\infty\right]$ is 
\emph{continuous}\footnote{Technically, this is sequential continuity, but this is the same as continuity in our 
metric space setting.} if for every $(x,y)\in S$ and every sequence $(x_n,y_n)_\nnn$ in $S$ that converges to $(x,y)$, we have $h(x_n,y_n)\to h(x,y)$. 
\end{defn}

\begin{remark}
\label{r:continuous}
Note that if a function $h\colon S\to \left]-\infty, +\infty\right]$ is continuous, then $h|_\Omega$ is continuous for every nonempty 
subset $\Omega$ of $S$. 
\end{remark}

\begin{proposition}
\label{l:continuous}
Suppose that $X$ and $Y$ are metric spaces, and 
that the cost function $c$ has open domain (i.e., $D$ is open) and $c|_{D}$ is continuous. Let $\Omega$ be a nonempty subset of $D$ and assume that $f$ is a $c$-antiderivative 
of $\Omega$. Then there exists exactly one continuous function $\widetilde{f}\colon \overline{\Omega}\cap D\to \RR$ that is a $c$-antiderivative of $\overline{\Omega}\cap D$ and $\widetilde{f}|_\Omega = f$. 
\end{proposition}
\begin{proof}
Take $(x,y)\in (\overline{\Omega}\cap D)\smallsetminus\Omega$, and get a sequence 
$(x_n,y_n)_\nnn$ in $\Omega$ such that $(x_n,y_n)\to (x,y)$. 
Let $\varepsilon>0$.
By continuity of $c|_{D}$ and openness of $D$, there exists $N\in\NN$ such that 
for all $m,n\geq N$, we have 
$|c(x_m,y_m)-c(x_n,y_m)|=|c|_D(x_m,y_m)-c|_D(x_n,y_m)|<\varepsilon$ and 
$|c(x_n,y_n)-c(x_m,y_n)|=|c|_D (x_n,y_n)-c|_D (x_m,y_n)|<\varepsilon$. 
Let $m,n\geq N$.
Now $f$ is a $c$-antiderivative of $\Omega$; thus, 
\begin{align*}
c(x_m,y_m)-c(x_n,y_m)&\leq f(x_m,y_m)-f(x_n,y_n),\\
c(x_n,y_n)-c(x_m,y_n)&\leq f(x_n,y_n)-f(x_m,y_m).
\end{align*}
This implies that 
$|f(x_m,y_m)-f(x_n,y_n)|\leq \max\{|c(x_m,y_m)-c(x_n,y_m)|, 
|c(x_n,y_n)-c(x_m,y_n)|\}<\varepsilon$. 
We've shown that $(f(x_n,y_n))_\nnn$ is Cauchy and hence convergent. 
Moreover, this limit is \emph{independent}%
\footnote{
Indeed, 
if two sequences in $\Omega$ satisfy 
$(x_{1,n},y_{1,n})\to (x,y)$ and 
$(x_{2,n},y_{2,n})\to (x,y)$, then 
there exist real numbers $\ell_1,\ell_2$ such that 
$f(x_{1,n},y_{1,n})\to\ell_1$ and
$f(x_{2,n},y_{2,n})\to\ell_2$.
Now consider the interlaced sequence 
$(x_n,y_n)_\nnn = ((x_{1,0},y_{1,0}),(x_{2,0},y_{2,0}),(x_{1,1},y_{1,1}),(x_{2,1},y_{2,1}), \ldots)_\nnn$ 
which converges to $(x,y)$ and hence $f(x_n,y_n)\to\ell$ for some 
$\ell\in\RR$. Then $f(x_{2n},y_{2n})\to \ell_1$ and $f(x_{2n+1},y_{2n+1})
\to \ell_2$; thus $\ell_1=\ell=\ell_2$.}
of the sequence $(x_n,y_n)_\nnn$. 
This gives rise to the extension
\begin{align}
\label{e:fwtilde}
\widetilde{f}\colon \overline{\Omega}\cap D \to\RR\colon 
(x,y) \mapsto \begin{cases}
f(x,y), &\text{if $(x,y)\in \Omega$;}\\
\lim_n f(x_n,y_n), &\text{if $(x,y)\in (\overline{\Omega}\cap D)\smallsetminus \Omega$.}
\end{cases}
\end{align}
We claim that $\widetilde{f}$ is a $c$-antiderivative of $\overline{\Omega}\cap D$.
To this end, let $(x,y),(u,v)$ be in $\overline{\Omega}\cap D$. 
We are done as soon as we can show that 
\begin{equation}
\label{e:0713g}
c(x,y)-c(u,y) \leq \widetilde{f}(x,y)-\widetilde{f}(u,v). 
\end{equation}
By the definition of $\widetilde{f}$ in \cref{e:fwtilde}, 
there exist sequences $(x_n,y_n)_\nnn$ and $(u_n,v_n)_\nnn$ in 
$\Omega$ such that $(x_n,y_n)\to (x,y)$, 
$(u_n,v_n)\to (u,v)$, $f(x_n,y_n)\to \widetilde{f}(x,y)$, 
and $f(u_n,v_n)\to \widetilde{f}(u,v)$: 
indeed, if $(x,y)\in\Omega$, we pick $(x_n,y_n)_{\nnn} = (x,y)$; 
if $(u,v)\in\Omega$, we pick $(u_n,v_n)_{\nnn} = (u,v)$;
in the remaining cases, existence is guaranteed by density of $\Omega$ in $\overline{\Omega}\cap D$. 

Because $f$ is a $c$-antiderivative of $\Omega$, we have 
\begin{equation}
(\forall\nnn)\quad 
c(x_n,y_n)-c(u_n,y_n)\leq f(x_n,y_n)-f(u_n,v_n). 
\end{equation}
We now take the limit as $n\to\pinf$; indeed, using 
the continuity of $c$ and \cref{e:fwtilde}, we deduce that 
$c(x,y)-c(u,y)\leq \widetilde{f}(x,y)-\widetilde{f}(u,v)$, i.e., 
\cref{e:0713g} holds. Now we show the continuity of $\widetilde{f}$. Let $(x,y)\in \overline{\Omega}\cap D$ and $(x_n,y_n)_{\nnn}$ be a sequence in $\overline{\Omega}\cap D$ such that $(x_n,y_n)\to (x,y)$. Since $\widetilde{f}$ is a $c$-antiderivative of $\overline{\Omega}\cap D$, we have 
\begin{align*}
c(x_n,y_n)-c(x,y_n)&\leq \widetilde{f}(x_n,y_n)-\widetilde{f}(x,y),\\
c(x,y)-c(x_n,y)&\leq \widetilde{f}(x,y)-\widetilde{f}(x_n,y_n).
\end{align*}
Hence, for $n$ large enough, $|\widetilde{f}(x_n,y_n)-\widetilde{f}(x,y)|\leq \max\{|c(x_n,y_n)-c(x,y_n)|, 
|c(x,y)-c(x_n,y)|\}=\max\{|c|_D(x_n,y_n)-c|_D(x,y_n)|, 
|c|_D(x,y)-c|_D(x_n,y)|\}\to 0$, which implies that $\widetilde{f}$ is continuous on $\overline{\Omega}\cap D$. 
In particular, $f$ is continuous. 
Finally, assume that $g$ has the same 
properties as $\widetilde{f}$: $g$ is a $c$-antiderivative 
of $\overline{\Omega}\cap D$, $g$ is continuous, and 
$g|_\Omega = f$. 
Clearly, 
$\widetilde{f}|_\Omega = f = g|_\Omega$. 
Now let $(x,y)\in (\overline{\Omega}\cap D)\smallsetminus 
\Omega$ and let $(x_n,y_n)_{\nnn}$ be a sequence in $\Omega$ such that $(x_n,y_n)\to (x,y)$. 
Then 
$\widetilde{f}(x,y)\leftarrow f(x_n,y_n) \to g(x,y)$ 
because $\widetilde{f}$ and $g$ are continuous at $(x,y)$ 
and again 
$\widetilde{f}|_\Omega = f = g|_\Omega$. 
\end{proof}

{\color{black}
\begin{theorem}
\label{t:continuous}
Suppose that $X$ and $Y$ are separable metric spaces, and 
that the cost function $c$ has open domain and $c|_D$ is continuous. 
Let $G$ be a nonempty subset of $D$. 
Then: 
\begin{align*}
\text{$G$ has an upper semicontinuous $c$-potential} \Leftrightarrow \text{$G$ has a continuous $c$-antiderivative} 
\Leftrightarrow \text{$G$ is $c$-path bounded.}
\end{align*}
\end{theorem}

\begin{proof}
``$G$ has a continuous $c$-antiderivative $\Leftarrow$ $G$ is $c$-path bounded":
The separability assumption on $X$ and $Y$ 
yields the separability of $X \times Y$. 
Hence $G$ is separable as well.
Now take a countable subset $\Omega$ of $G$ that is dense in $G$.
It follows from \cref{r:inclusion} that 
$F_\Omega \leq F_G$; thus, $\Omega$ is $c$-path bounded (because 
$G$ is). 
By \cref{p:countableset}, $\Omega$ has a $c$-antiderivative, 
say $f$. 
Hence $\overline{\Omega}\cap D$ has a 
continuous $c$-antiderivative, say $\widetilde{f}$, 
according to \cref{l:continuous}. Since $G\subseteq (\overline{\Omega}\cap D)$, $\widetilde{f}|_G$ is a continuous $c$-antiderivative of $G$ by 
\cref{r:easyRemark} and \cref{r:continuous}. 

\par ``$G$ has an upper semicontinuous $c$-potential $\Leftarrow$ $G$ has a continuous $c$-antiderivative": Assume that  $f: G \to \mathbb{R}$ is a continuous $c$-antiderivative of $G$ and $f_G$ is defined in \cref{r:keyremark}. Since $c$ has open domain and $c|_D$ is continuous, it follows that for every $(u,v)\in G$, 
the function $x\mapsto c(x, v)-c(u,v)+f_G(u)$ is 
upper semicontinuous 
 on $X$.
Hence, 
\begin{align*}
    \text{$\Psi: X \to [-\infty, +\infty]: x\mapsto \inf_{(u,v)\in G} \big(c(x,v)-c(u,v)+f_G(u)\big)$ is upper semicontinuous.}
\end{align*}
On the other hand, $\Psi$ is a $c$-potential of $G$ by \cref{f:et}. Therefore, $\Psi$ is an upper semicontinuous $c$-potential of $G$. 

\par ``$G$ has an upper semicontinuous $c$-potential $\Rightarrow$ $G$ is $c$-path bounded": This follows from \cref{f:0707a}. 
\end{proof}
}

Let us now present applications of the results in this section to some well-known non-traditional costs.
We start with a generalization of the polar cost considered 
earlier in \cref{ex:polar}: 

\begin{example}[polar cost revisited]
Suppose that $X=Y$ is a real separable Hilbert space, with inner product $\langle \cdot, \cdot\rangle$. 
Following Artstein-Avidan, Sadovsky and Wyczesany \cite{Kasia2}, we assume 
that the cost function $c$ is the polar cost 
\begin{align*}
c\colon X \times X \to\RX \colon 
(x,y)\mapsto 
\begin{cases}
-\ln\big(\langle x, y\rangle-1\big), &\text{if $\langle x, y\rangle >1$;}\\
\pinf, &\text{if $\langle x, y \rangle \leq 1$.}
\end{cases}
\end{align*}
Then the domain 
$D=\{(x,y)\in X\times X \ |\ \langle x, y\rangle>1\}$ is open and {\color{black}$c|_D$} is continuous. 
Therefore, \cref{t:0713c} and \cref{t:continuous} imply the 
following for any nonempty subset $G$ of $D$: 
\begin{equation}
\text{$G$ is $c$-cyclically monotone and (top.) connected}
\;\Rightarrow\;
\text{$G$ is $c$-path bounded}
\;\Leftrightarrow\;
{\color{black}\text{$G$ has an upper semicontinuous $c$-potential}}.
\end{equation}
\end{example}

\begin{example}[a cost related to the Hellinger-Kantorovich distance]
Suppose that $X=Y$ is a separable metric space with metric $d$, 
and that the cost function $c$ is given by 
\begin{align*}
c\colon X \times X \to\RX \colon 
(x,y)\mapsto 
\begin{cases}
-\ln\big(\cos^2 d(x,y)\big), &\text{if $d(x,y)<\frac{\pi}{2}$;}\\
\pinf, &\text{if $d(x,y)\geq \frac{\pi}{2}$.}
\end{cases}
\end{align*}
This cost function is known to induce the 
so-called 
\emph{Hellinger-Kantorovich distance}; see, e.g., 
\cite{Liero} for details.
Then  $D=\{(x,y)\in X\times X \ |\ d(x,y)<\frac{\pi}{2}\}$ 
is open and {\color{black}$c|_D$} is continuous.
Again, \cref{t:0713c} and \cref{t:continuous} imply the 
following for any nonempty subset $G$ of $D$: 
\begin{equation}
\text{$G$ is $c$-cyclically monotone and (top.) connected}
\;\Rightarrow\;
\text{$G$ is $c$-path bounded}
\;\Leftrightarrow\;
{\color{black}\text{$G$ has upper semicontinuous $c$-potential}}.
\end{equation}
\end{example}

Finally, we revisit the Coulomb cost from \cref{ex:Coulomb}. 

\begin{example}[Coulomb cost revisited]
Suppose that $X=Y=\RR^n$ {\color{black} is 
equipped with the Euclidean norm $\|\cdot\|$} and that $c$ is the 
\emph{Coulomb cost}, given by
\begin{equation}
c\colon X\times Y \to\RX \colon 
(x,y)\mapsto 
\begin{cases}
\frac{1}{\|x-y\|}, &\text{if $x\neq y$;}\\
\pinf, &\text{if $x=y$,}
\end{cases}
\end{equation}
and with $D = \dom c = \menge{(x,y)\in X\times Y}{x\neq y}$. 
Clearly, $D$ is open and 
{\color{black}$c|_D$} is continuous. 
Again, \cref{t:0713c} and \cref{t:continuous} imply the 
following for any nonempty subset $G$ of $D$: 
\begin{equation}
\text{$G$ is $c$-cyclically monotone and (top.) connected}
\;\Rightarrow\;
\text{$G$ is $c$-path bounded}
\;\Leftrightarrow\;
{\color{black}\text{$G$ has upper semicontinuous $c$-potential}}.
\end{equation}
Compared to \cref{ex:Coulomb}, we not only get the equivalence 
between $c$-path boundedness and $c$-potentiability, but also {\color{black}upper semicontinuity of $c$-potential} and  
the useful sufficient condition of 
$c$-cyclic monotonicity and (topological) connectedness.
\end{example}

\section{Monotone costs and $c$-path bounded extensions in $\RR^2$}
\label{s:extension}
{\color{black} In this section, we discuss $c$-potentiability in scenarios where none of the previous results is applicable.
In order to establish the existence of 
such a $c$-potential,  
we introduce the notion of a $c$-path bounded extension in \cref{d:c-path bounded ex}. We use it to prove $c$-potentiability for a class of costs in $\mathbb{R}^2$, which we call monotone costs (\cref{d:monocost}). The main result of this section is \cref{t:cpathextension}, with which, finally, in \cref{ex:5.20}, we provide $c$-potentiability of the set 
$G$, presented in the following  \cref{ex:5.1}.}

\begin{example}
\label{ex:5.1}
Suppose that $X=Y=\RR$ and that 
\begin{equation}
c\colon \mathbb{R}\times \mathbb{R}\to \left]-\infty, +\infty\right] \colon (x,y) \mapsto  \begin{cases}
-\ln(xy-1)+1, &\text{if $xy>1$ and $y\in \mathbb{Q}$;}\\
-\ln(xy-1), &\text{if $xy>1$ and $y\notin \mathbb{Q}$;}\\
\pinf, &\text{otherwise.}
\end{cases}
\end{equation}
Let $(\varepsilon_n)_{n\in \mathbb{N}}$ be a sequence in 
$\mathopen]0, 1\mathclose[$ that strictly decreases to $0$. Let $M>0$ and $\alpha_n, \beta_n\in \mathopen]\frac{1}{\varepsilon_n}, \min\{\frac{1}{\varepsilon_{n+1}}, \frac{M+1}{\varepsilon_n}\}\mathclose[$ such that $\alpha_n<\beta_n$. Set $G_1:=\bigcup_{n\in \mathbb{N}}(\left[\alpha_n, \beta_n\right] \times \{\varepsilon_n\})$, $G_2:=\bigcup_{n\in \mathbb{N}} (\{\varepsilon_n\}\times \left[\alpha_n, \beta_n\right])$ and $G:=G_1 \cup G_2$. Then
\begin{align*}
\text{$G$ is $c$-path bounded.} 
\end{align*}
The detailed proof for $c$-path boundedness is presented in \cref{ex:5.20}.
\begin{figure}[h]
\centering
\begin{tikzpicture}
  \begin{axis}[
     width=10cm,
     height=10cm,
    axis lines=middle,
    xmin=0, xmax=5,
    ymin=0, ymax=5,
    samples=200,
    domain=0.17:5,
    xlabel={$x$},
    ylabel={$y$},
    xtick=\empty,
    ytick=\empty,
    axis line style={->},
    enlargelimits,
    clip=false
  ]
    \addplot [
      name path=lower,
      domain=0.17:5,
      draw=none
    ] {1/x};

    \path[name path=upper] (axis cs:0.17,5) -- (axis cs:5,5);

    \addplot [
      fill=green!20,
      draw=none
    ] fill between [
      of=upper and lower,
      soft clip={domain=0.17:6}
    ];

    \draw[thick, blue] (axis cs:1,1.3) -- (axis cs:1,2);
    \draw[dashed, thick, red] (axis cs:1,1.3) -- (axis cs:0,1.3);
    \node[anchor=east] at (axis cs:0, 1.3) {\small $\alpha_0$};
    \draw[dashed, thick, red] (axis cs:1,2) -- (axis cs:0,2);
    \node[anchor=east] at (axis cs:0, 2) {\small $\beta_0$};
    \draw[dashed, thick, red] (axis cs:1,1.3) -- (axis cs:1,0);
    \node[anchor=north] at (axis cs:1, 0) {\small $\varepsilon_0$};

    \draw[thick, blue] (axis cs:0.5,2.3) -- (axis cs:0.5,3);
    \draw[dashed, thick, red] (axis cs:0.5,2.3) -- (axis cs:0,2.3);
    \node[anchor=east] at (axis cs:0, 2.3) {\small $\alpha_1$};
    \draw[dashed, thick, red] (axis cs:0.5,3) -- (axis cs:0,3);
    \node[anchor=east] at (axis cs:0, 3) {\small $\beta_1$};
    \draw[dashed, thick, red] (axis cs:0.5,2.3) -- (axis cs:0.5,0);
    \node[anchor=north] at (axis cs:0.6, 0) {\small $\varepsilon_1$};
    
    \draw[thick, blue] (axis cs:0.33,3.3) -- (axis cs:0.33,4);
    \draw[dashed, thick, red] (axis cs:0.33,3.3) -- (axis cs:0,3.3);
    \node[anchor=east] at (axis cs:0, 3.3) {\small $\alpha_2$};
    \draw[dashed, thick, red] (axis cs:0.33,4) -- (axis cs:0,4);
    \node[anchor=east] at (axis cs:0, 4) {\small $\beta_2$};
    \draw[dashed, thick, red] (axis cs:0.33,3.3) -- (axis cs:0.33,0);
    \node[anchor=north] at (axis cs:0.33, 0) {\small $\varepsilon_2$};
    
    \draw[thick, blue] (axis cs:1.3,1) -- (axis cs:2,1);
    \draw[dashed, thick, red] (axis cs:1.3,1) -- (axis cs:1.3,0);
    \node[anchor=north] at (axis cs:1.3, 0) {\small $\alpha_0$};
    \draw[dashed, thick, red] (axis cs:2,1) -- (axis cs:2,0);
    \node[anchor=north] at (axis cs:2, 0) {\small $\beta_0$};
    \draw[dashed, thick, red] (axis cs:1.3,1) -- (axis cs:0,1);
    \node[anchor=east] at (axis cs:0, 1) {\small $\varepsilon_0$};

    \draw[thick, blue] (axis cs:2.3,0.5) -- (axis cs:3,0.5);
    \draw[dashed, thick, red] (axis cs:2.3,0.5) -- (axis cs:2.3,0);
    \node[anchor=north] at (axis cs:2.3, 0) {\small $\alpha_1$};
    \draw[dashed, thick, red] (axis cs:3,0.5) -- (axis cs:3,0);
    \node[anchor=north] at (axis cs:3, 0) {\small $\beta_1$};
    \draw[dashed, thick, red] (axis cs:2.3,0.5) -- (axis cs:0,0.5);
    \node[anchor=east] at (axis cs:0, 0.5) {\small $\varepsilon_1$};
    
    \draw[thick, blue] (axis cs:3.3,0.33) -- (axis cs:4,0.33);
    \draw[dashed, thick, red] (axis cs:3.3,0.33) -- (axis cs:3.3,0);
    \node[anchor=north] at (axis cs:3.3, 0) {\small $\alpha_2$};
    \draw[dashed, thick, red] (axis cs:4,0.33) -- (axis cs:4,0);
    \node[anchor=north] at (axis cs:4, 0) {\small $\beta_2$};
    \draw[dashed, thick, red] (axis cs:3.3,0.33) -- (axis cs:0,0.33);
    \node[anchor=east] at (axis cs:0, 0.33) {\small $\varepsilon_2$};

      \fill[black] (axis cs:4.5,0.33) circle[radius=1pt];
     \fill[black] (axis cs:4.6, 0.33) circle[radius=1pt];
     \fill[black] (axis cs:4.7, 0.33) circle[radius=1pt];
     
    \fill[black] (axis cs:0.33,4.5) circle[radius=1pt];
     \fill[black] (axis cs:0.33, 4.6) circle[radius=1pt];
     \fill[black] (axis cs:0.33, 4.7) circle[radius=1pt];

  \end{axis}
\end{tikzpicture}
\caption{{\color{black} Depicted is a visualization of \cref{ex:5.1}. The green area is the domain of $c$ in the first quadrant, and $G$ consists of the union of the blue segments.}
}
    \label{fig:1}
\end{figure}
However, none of the previously known results regarding $c$-potentiability is applicable:

\begin{enumerate}
\item $G$ is uncountable and has an infinite black hole $G_1$ \footnote{Please see \cref{r:5.2} below for details.}, so \cref{f:fix} cannot be applied.

\item Suppose that $\Omega$ is a nonempty subset of $G$ and for every $(x,y)\in G$, we have 
$$\sup_{(x_e,y_e)\in \Omega} F_G\big((x,y), (x_e,y_e)\big)\in \mathbb{R}.$$ Then there is a sequence $(k_n)_{\nnn}$ such that 
{\color{black} $k_n\to+\infty$ and }
$\{(x_{k_n}, y_{k_n})\}_{\nnn}\subseteq \Omega$ where $(x_{k_n}, y_{k_n})\in \left[\alpha_{k_n}, \beta_{k_n}\right]\times \{\varepsilon_{k_n}\}$. Fix $(x,y)\in \left[\alpha_0, \beta_0\right]\times \{\varepsilon_0\}$. 
Then 
\begin{align*}
F_G\big((x,y), (x_{k_n}, y_{k_n})\big)\geq c(x,y)-c(x_{k_n}, y)=\ln\Big(\frac{x_{k_n}y-1}{xy-1}\Big)\to +\infty.
\end{align*}
Hence, $\sup_{(x_e,y_e)\in \Omega} F_G((x,y), (x_e,y_e))=+\infty$ which contradicts our assumption. So the hypothesis of \cref{l:lemma1} is not satisfied, and \cref{l:lemma1} is not applicable. Moreover, the condition of \cref{t:exist1} is more restrictive than that of \cref{l:lemma1}, so \cref{t:exist1} is not applicable either. 

\item Suppose that $A$ is a nonempty subset of $G$ and for every $(x,y)\in G$, $\sup_{(x_s,y_s)\in A} F_G\big((x_s,y_s), (x,y)\big)\in \mathbb{R}$. Then there is a sequence $(k_n)_{\nnn}$ such that $\{(x_{k_n}, y_{k_n})\}_{\nnn}\subseteq A$ where $(x_{k_n}, y_{k_n})\in  \{\varepsilon_{k_n}\}\times \left[\alpha_{k_n}, \beta_{k_n}\right]$. Fix $(x,y)\in \{\varepsilon_0\}\times \left[\alpha_0, \beta_0\right]$. Then we have
\begin{equation}
\label{e:5.1e1}
F_G\big((x_{k_n},y_{k_n}), (x, y)\big)\geq c(x_{k_n},y_{k_n})-c(x, y_{k_n})=\ln\big(\frac{xy_{k_n}-1}{x_{k_n}y_{k_n}-1}\big).
\end{equation}
Note that $xy_{k_n}-1\geq x\alpha_{k_n}-1>\frac{x}{\varepsilon_{k_n}}-1\to +\infty$ and $x_{k_n}y_{k_n}-1\leq \varepsilon_{k_n} \beta_{k_n}-1\leq M$. Hence, $\ln((xy_{k_n}-1)/(x_{k_n}y_{k_n}-1))\to +\infty$ and  $F_G((x_{k_n},y_{k_n}), (x, y))\to +\infty$ by \cref{e:5.1e1}. So $\sup_{(x_s,y_s)\in A} F_G((x_s,y_s), (x,y))\in \mathbb{R}$ is false; 
thus, the hypothesis of \cref{l:lemma2} is not satisfied, and \cref{l:lemma2} is not applicable. Moreover, the condition of \cref{t:exist2} is more restrictive than that of \cref{l:lemma2}, so \cref{t:exist2} is not applicable either. 

\item The set $G$ is not strongly connected, so it is not ball chain connected by \cref{e:0713b} and not topologically connected by \cref{p:0713a}. Therefore, \cref{t:0713c} is not applicable. 

\item The cost $c$ is not continuous, so \cref{t:continuous} is not applicable. 
\end{enumerate}
\end{example}

\begin{remark}
\label{r:5.2}
Note the following for \cref{ex:5.1}:
\begin{enumerate}
\item If $(x,y)\in \{\varepsilon_n\}\times \left[\alpha_n, \beta_n\right]$ and $(u,v)\in \left[\alpha_m, \beta_m\right]\times \{\varepsilon_m\}$ then $xv=\varepsilon_n \varepsilon_m\leq 1$ and  $c(x,v)=+\infty$.

\item If $(x,y)\in \{\varepsilon_n\}\times \left[\alpha_n, \beta_n\right]$ and $(u,v)\in \{\varepsilon_m\}\times \left[\alpha_m, \beta_m\right]$ where $n>m$, then $xv=\varepsilon_n v<\frac{\varepsilon_n}{\varepsilon_{m+1}}\leq 1$ and $c(x,v)=+\infty$.

\item If $(x,y)\in \left[\alpha_n, \beta_n\right] \times \{\varepsilon_n\}$ and $(u,v)\in \left[\alpha_m, \beta_m\right]\times \{\varepsilon_m\}$ where $m>n$, then $xv=x\varepsilon_m <\frac{\varepsilon_m}{\varepsilon_{n+1}}\leq 1$ and $c(x,v)=+\infty$.
\end{enumerate}
 {\color{black} Hence, if $(x,y)$ is on a line segment to the left of the line segment on which $(u,v)$ lies, then there is no path from $(u,v)$ to $(x,y)$ in $G$. }
\end{remark}

\subsection{$c$-path bounded extension and ordering in $\mathbb{R}^2$}
In order to prove $c$-potentiability of a set $G$ in a scenario such as in \cref{ex:5.1}, we will employ
the notion of a $c$-path bounded extension:

\begin{defn}[$c$-path bounded extension]
\label{d:c-path bounded ex}
    {\color{black} Let $X, Y$ be two nonempty sets and $c: X\times Y\to \left]-\infty, +\infty\right]$ be a cost function with nonempty domain $D$.} Let $G$ be a $c$-path bounded subset of $D$. 
    We say that a subset $S$ of $D$ is a \emph{$c$-path bounded extension of $G$} if $S$ is $c$-path bounded and $G\subseteq S$.  
\end{defn}

{\color{black}
\begin{remark}
We focus our attention in this section on $\mathbb{R}\times \mathbb{R}$, however, we note that a $c$-path bounded extension is defined in a more general setting. The idea of this notion is that if we can find a larger $c$-path bounded set $S$, then there will be more paths from one point to another and some disconnected or unilaterally connected pairs may become strongly connected. 
In a sense, if a superset $S$ 
of $G$ is ``more potentiable" than the original set $G$, then $S$ is more likely to have a $c$-antiderivative $f$ whose restriction $f|_G$ is automatically a $c$-antiderivative of $G$ by \cref{r:easyRemark}.
\end{remark}
}

{\color{black} From now on, we let $\|\cdot\|$ be the Euclidean norm on $\mathbb{R}\times \mathbb{R}$.}  We start by introducing a binary relation on $\mathbb{R}\times \mathbb{R}$ via
\begin{align*}
    (x,y)\preccurlyeq^{\oplus} (u,v) \Leftrightarrow x\leq u \text{ and } y\leq v,
\end{align*}
and we shall write $(x,y)\prec ^{\oplus}(u,v)$ if $x<u$ and $y<v$. {\color{black} Moreover, we say that $(x,y)$ is the maximal (resp. minimal) element of a chain $(G, \preccurlyeq^\oplus)$ if $(x,y)\in G$ and $(u,v)\preccurlyeq^\oplus (x,y)$ (resp. $(x,y)\preccurlyeq^\oplus (u,v)$) for all $(u,v)\in G$. If $(G, \preccurlyeq^\oplus)$ has the maximal (resp. minimal) element, we denote it by $\max_{\preccurlyeq^\oplus}$ (resp. $\min_{\preccurlyeq^\oplus}$). 

Similarly, we define
\begin{align*}
    (x,y)\preccurlyeq^{\ominus} (u,v)\Leftrightarrow x\leq u \text{ and } y\geq v,
\end{align*}
and we write $(x,y)\prec^{\ominus}(u,v)$ if $x<u$ and $y>v$. Moreover, we say that $(x,y)$ is the maximal (resp. minimal) element of a chain $(G, \preccurlyeq^\ominus)$ if $(x,y)\in G$ and $(u,v)\preccurlyeq^\ominus (x,y)$ (resp. $(x,y)\preccurlyeq^\ominus (u,v)$) for all $(u,v)\in G$. If $(G, \preccurlyeq^\ominus)$ has the maximal (resp. minimal) element, we denote it by $\max_{\preccurlyeq^\ominus}$ (resp. $\min_{\preccurlyeq^\ominus}$). }

\begin{remark}
    \label{r:chain}
    Note that $\preccurlyeq^{\oplus}$ and $\preccurlyeq^{\ominus}$ are partial orders on $\mathbb{R}\times \mathbb{R}$. Moreover, 
    \begin{align*}
    &\lnot\big( (x,y)\prec^{\oplus}(u,v) \vee (u,v)\prec^{\oplus}(x,y)\big) \Leftrightarrow (x,y)\preccurlyeq^{\ominus}(u,v) \vee (u,v)\preccurlyeq^{\ominus}(x,y)\\
    & \lnot \big((x,y)\prec^\ominus (u,v)\vee (u,v)\prec^\ominus (x,y))\Leftrightarrow (x,y)\preccurlyeq^\oplus (u,v) \vee (u,v) \preccurlyeq^\oplus (x,y).
    \end{align*}
\end{remark}
{\color{black} $(\mathbb{R}\times \mathbb{R}, \preccurlyeq^\ominus)$ and $(\mathbb{R}\times \mathbb{R}, \preccurlyeq^\oplus)$ are order isomorphic with the isometry 
\begin{equation}
\label{e:order isomorphic}
    \theta: \mathbb{R}\times \mathbb{R}\to \mathbb{R}\times \mathbb{R}: (x,y)\mapsto (x, -y).
\end{equation}
Consequently, we will focus only on $\preccurlyeq^\ominus$ for the rest of this subsection because all the statements and proofs work similarly for $\preccurlyeq^\oplus$ by \cref{e:order isomorphic}.
}

The following remark and lemmas provide basic observations from the definition of $\preccurlyeq^\ominus$ which we will use to prove results later in this section. 
\begin{remark} The following assertions hold: 
    \label{r:cluster}
    \begin{enumerate}
        \item If $(G, \preccurlyeq^\ominus)$ is a chain, then $(\overline{G}, \preccurlyeq^\ominus)$ is a chain. 
        \item If $(x_1, y_1)\preccurlyeq^\ominus (a, b)\preccurlyeq^\ominus (x_2, y_2)$, then $\|(x_1, y_1)-(a,b)\|\leq \|(x_1, y_1)-(x_2,y_2)\|$ and $\|(x_2,y_2)-(a,b)\|\leq \|(x_2,y_2)-(x_1,y_1)\|$. 
    \end{enumerate}
\end{remark}

\begin{lemma}
    \label{l:keylemma}
    Assume that $(G, \preccurlyeq^\ominus)$ is a chain and $(x_1, y_1), (x_2,y_2)\in G$ satisfy 
    \begin{enumerate}
        \item $(x_1, y_1)\preccurlyeq^\ominus (x_2,y_2)$;
        \item and there exist $\delta_1, \delta_2>0$ such that $\mathcal{B}_{\delta_1}(x_1, y_1)\cap \mathcal{B}_{\delta_2}(x_2,y_2)\cap G=\varnothing$.
    \end{enumerate}
    Then $(a_1,b_1)\preccurlyeq^\ominus (a_2,b_2)$ for all $(a_1,b_1)\in \mathcal{B}_{\delta_1}(x_1,y_1)\cap G$, $(a_2,b_2)\in \mathcal{B}_{\delta_2}(x_2,y_2)\cap G$. 
\end{lemma}

\begin{proof}
By contradiction, suppose there are $(a_1,b_1)\in \mathcal{B}_{\delta_1}(x_1,y_1)\cap G$ and $(a_2,b_2)\in \mathcal{B}_{\delta_2}(x_2,y_2)\cap G$ such that $(a_2,b_2)\preccurlyeq^\ominus(a_1,b_1)$. 

\par \textit{Case 1:} Suppose $(x_1,y_1)\preccurlyeq^\ominus(a_2,b_2)\preccurlyeq^\ominus(a_1,b_1)$. Then  $\|(x_1, y_1)-(a_2,b_2)\|\leq \|(x_1, y_1)-(a_1, b_1)\|<\delta_1$ by \cref{r:cluster}.(ii). Hence, $(a_2,b_2)\in \mathcal{B}_{\delta_1}(x_1,y_1)\cap G$, which is absurd because $(a_2,b_2)\in \mathcal{B}_{\delta_2}(x_2,y_2)\cap G$ and $\mathcal{B}_{\delta_2}(x_2,y_2)\cap\mathcal{B}_{\delta_1}(x_1, y_1)\cap G=\varnothing$.
    
    \par \textit{Case 2:} Suppose $(a_2,b_2)\preccurlyeq^\ominus(x_1, y_1)\preccurlyeq^\ominus(a_1,b_1)$. Together with assumption (i), we know $(a_2,b_2)\preccurlyeq^\ominus(x_1, y_1)\preccurlyeq^\ominus(x_2,y_2)$. By \cref{r:cluster}.(ii), we have $\|(x_2,y_2)-(x_1,y_1)\|\leq \|(x_2,y_2)-(a_2,b_2)\|<\delta_2$, so $(x_1, y_1)\in \mathcal{B}_{\delta_2}(x_2,y_2)\cap G$, which is absurd.

  \par \textit{Case 3:} Suppose $(a_2,b_2)\preccurlyeq^\ominus(a_1,b_1)\preccurlyeq^\ominus(x_1, y_1)$. Then $(a_2,b_2)\preccurlyeq^\ominus(a_1,b_1)\preccurlyeq^\ominus(x_1, y_1)\preccurlyeq^\ominus (x_2, y_2)$. By \cref{r:cluster}.(ii), we know $\|(x_2, y_2)-(x_1, y_1)\|\leq \|(x_2, y_2)-(a_2, b_2)\|<\delta_2$, so $(x_1, y_1)\in \mathcal{B}_{\delta_2}(x_2,y_2)\cap G$ which is absurd. 
\end{proof}

\begin{lemma}
    \label{l:keylemma2}
    Assume that $(G, \preccurlyeq^\ominus)$ is a chain and $(x_1, y_1), (x_2,y_2)\in G$ satisfy $(x_1, y_1)\preccurlyeq^\ominus(x_2,y_2)$. Let $\delta>0$. Then
    
    \begin{enumerate}
        \item For any $(a_2,b_2)\in \mathcal{B}_\delta (x_2,y_2)\cap G$, there is a $(a_1,b_1)\in \mathcal{B}_\delta (x_1,y_1)\cap G$ such that $(a_1,b_1)\preccurlyeq^\ominus (a_2,b_2)$;
        \item For any $(a_1,b_1)\in \mathcal{B}_\delta (x_1,y_1)\cap G$, there is a $(a_2,b_2)\in \mathcal{B}_\delta (x_2,y_2)\cap G$ such that $(a_1,b_1)\preccurlyeq^\ominus (a_2,b_2)$. 
    \end{enumerate}
\end{lemma}
\begin{proof}
   (i): Take $(a_2, b_2)\in \mathcal{B}_\delta (x_2,y_2)\cap G$ and assume $(a_2, b_2)\preccurlyeq^\ominus(x_1,y_1)$ otherwise $(a_1, b_1):=(x_1, y_1)$ is the desired point. By \cref{r:cluster}.(ii), $\|(a_2,b_2)-(x_1,y_1)\|\leq \|(a_2, b_2)-(x_2,y_2)\|<\delta$, so $(a_2, b_2)\in \mathcal{B}_\delta(x_1,y_1)\cap G$. Therefore, $(a_1, b_1):=(a_2, b_2)$ is the desired point.  

    \par (ii): Take $(a_1, b_1)\in \mathcal{B}_\delta (x_1,y_1)\cap G$ and assume $(x_2, y_2)\preccurlyeq^\ominus(a_1,b_1)$ otherwise $(a_2, b_2):=(x_2,y_2)$ is the desired point. By \cref{r:cluster}.(ii), $\|(a_1,b_1)-(x_2,y_2)\|\leq \|(a_1, b_1)-(x_1,y_1)\|<\delta$, so $(a_1, b_1)\in \mathcal{B}_\delta(x_2,y_2)\cap G$. Therefore, $(a_2, b_2):=(a_1, b_1)$ is the desired point.  
\end{proof}

\subsection{Monotone costs in $\mathbb{R}^2$}
\begin{defn}[monotone costs]
\label{d:monocost}
    A cost $c:\mathbb{R}\times \mathbb{R}\to \left]-\infty,+\infty\right]$ is said to be $\ominus$-monotone if for all $(x_s,y_s), (x_e,y_e)\in D$, 
one has the implication 
{\color{black}
\begin{equation}
\label{e:2-mono}
[(x_s,y_s)\prec^\oplus (x_e, y_e) \vee (x_e,y_e)\prec^\oplus (x_s,y_s)]\Rightarrow c(x_s,y_s)-c(x_e,y_s)+c(x_e,y_e)-c(x_s,y_e)> 0.
\end{equation}
}
And $c$ is said to be $\oplus$-monotone if for all $(x_s,y_s), (x_e,y_e)\in D$, 
one has the implication 
{\color{black}
\begin{align*}
[(x_s,y_s)\prec^\ominus (x_e, y_e) \vee (x_e, y_e)\prec^\ominus(x_s,y_s)] \Rightarrow c(x_s,y_s)-c(x_e,y_s)+c(x_e, y_e)-c(x_s,y_e)> 0.
\end{align*}
}
\end{defn}

\begin{remark}
\label{r:separable}
Note that if $c$ is $\ominus$-monotone and $g,h\colon \mathbb{R}\to \mathbb{R}$ are any functions, then $\widetilde{c}(x,y):=c(x,y)+g(x)+h(y)$ is $\ominus$-monotone. Similarly, if $c$ is $\oplus$-monotone and $g,h \colon \mathbb{R}\to \mathbb{R}$ are any functions, then $\widetilde{c}(x,y):=c(x,y)+g(x)+h(y)$ is $\oplus$-monotone. 
\end{remark}

\begin{remark}
\label{r:smooth}
Suppose that $c \colon \mathbb{R}\times \mathbb{R} \to \left]-\infty, +\infty\right]$ has convex open domain and {\color{black}$\partial_x c$ is continuous on $D$.} Using calculus, one can show that $c$ is $\ominus$-monotone if and only if 
\begin{equation}
\label{e:5.6e1}
\begin{split}
&\quad \ \ \ [(x_s,y_s), (x_e,y_e)\in D \land (x_s,y_s)\prec^\oplus (x_e,y_e)]\\
&\Rightarrow [(x_s,y_e), (x_e,y_s)\in D \land \int_{x_s}^{x_e} \big(\partial_x c(x, y_e)-\partial_x c(x, y_s)\big) dx>0].
\end{split}
\end{equation}
Similarly, $c$ is $\oplus$-monotone if and only if 
\begin{equation}
\label{e:5.6e2}
\begin{split}
&\quad \ \ \ [(x_s,y_s), (x_e,y_e)\in D \land (x_s,y_s)\prec^\ominus (x_e,y_e)]\\
&\Rightarrow [(x_e,y_s), (x_s,y_e)\in D \land \int_{x_s}^{x_e} \big(\partial_x c(x, y_e)-\partial_x c(x, y_s)\big)dx>0].
\end{split}
\end{equation}
If we further assume that $c$ is {\color{black} continuously} twice differentiable, in particular, {\color{black} $\partial_{xy} c$ is continuous on $D$}, then $c$ is $\ominus$-monotone if and only if 
\begin{equation}
\label{e:5.6e3}
\begin{split}
 [(x_s,y_s), (x_e,y_e)\in D \land (x_s,y_s)\prec^\oplus (x_e,y_e)]
\Rightarrow [(x_s,y_e), (x_e,y_s)\in D \land \int_{x_s}^{x_e} \int_{y_s}^{y_e}\partial_{xy}c(x,y) dy\ dx>0];
\end{split}
\end{equation}
similarly, 
$c$ is $\oplus$-monotone if and only if
\begin{equation}
\label{e:5.6e4}
\begin{split}
 [(x_s,y_s), (x_e,y_e)\in D \land (x_s,y_s)\prec^\ominus (x_e,y_e)]
\Rightarrow [(x_s,y_e), (x_e,y_s)\in D \land \int_{x_s}^{x_e} \int_{y_e}^{y_s}\partial_{xy}c(x,y) dy\ dx<0].
\end{split}
\end{equation}
\end{remark}

\begin{remark}[monotonicity and the twist condition]
\label{r:twist}
{\color{black} Let the cost $c \colon \mathbb{R}\times \mathbb{R} \to \mathbb{R}$ be continuously differentiable with respect to the first variable, i.e.,  $\partial_x c$ is continuous on $\mathbb{R}\times \mathbb{R}$. In the context of optimal transport, $c$ is called 
\emph{twisted} if for every fixed $x\in \mathbb{R}$, the map 
\begin{align*}
y\mapsto \partial_{x} c(x, y) 
\end{align*}
is injective (see, for example \cite[Page~234]{Villani} and \cite[Definition~1.16]{Santambrogio}).}
This holds, in particular,  if 
\begin{enumerate}
\item $y \mapsto \partial_x c(x,y)$ 
is strictly increasing for all $x\in \mathbb{R}$
\end{enumerate}
or if 
\begin{enumerate}
\setcounter{enumi}{1}
\item $y \mapsto \partial_x c(x,y)$ 
is strictly decreasing for all $x\in \mathbb{R}$. 
\end{enumerate}
If (i) holds, then $\partial_x c(x, y_e)-\partial_x c(x,y_s)>0$ for $(x_s,y_s)\prec^\oplus (x_e,y_e)$ and $\forall x\in [x_s, x_e]$ hence \cref{e:5.6e1} holds and $c$ is $\ominus$-monotone. If (ii) holds, then $\partial_x c(x, y_e)-\partial_x c(x,y_s)>0$ for $(x_s,y_s)\prec^\ominus (x_e,y_e)$ and $\forall x\in [x_s, x_e]$ hence \cref{e:5.6e2} holds and $c$ is $\oplus$-monotone.
\end{remark}

\begin{example}[$\ominus$-monotone costs] \label{ex:5.11}\ 
\begin{enumerate}
\item Consider $c\colon \mathbb{R}\times \mathbb{R}\to \mathbb{R}\colon (x,y)\mapsto xy$. Then $\partial_{xy}c(x,y)=1$ for all $(x,y)\in \mathbb{R}\times \mathbb{R}$. Take $(x_s,y_s), (x_e,y_e)\in \mathbb{R}\times \mathbb{R}$ such that $(x_s,y_s)\prec^\oplus (x_e,y_e)$ and 
\begin{align*}
\int_{x_s}^{x_e}\int_{y_s}^{y_e} \partial_{xy}c(x,y) dy\ dx =(x_e-x_s)(y_e-y_s)>0.
\end{align*}
Hence, $c$ is $\ominus$-monotone by \cref{e:5.6e3}. 

\item Consider 
\begin{align*}
c\colon \mathbb{R}\times \mathbb{R}\to \left]-\infty, +\infty\right] \colon (x,y) \mapsto  \begin{cases}
-\ln(xy-1), &\text{if $0<\frac{1}{x}<y$;}\\
\pinf, &\text{otherwise.}
\end{cases}
\end{align*}
Then $D=\{(x,y)\in \mathbb{R}\times \mathbb{R}\ |\ 0<\frac{1}{x}<y\}$. For $(x_s, y_s), (x_e,y_e)\in D$ such that $(x_s, y_s)\prec^\oplus(x_e,y_e)$, we have $x_sy_s>1, x_ey_e>1$ and $x_s<x_e, y_s<y_e$, so $x_sy_e>x_sy_s>1$ and $x_ey_s>x_sy_s>1$. Hence, $(x_s,y_e), (x_e,y_s)\in D$. Moreover, $\partial_{xy}c(x,y)=\frac{1}{(xy-1)^2}>0$ so 
\begin{align*}
\int_{x_s}^{x_e}\int_{y_s}^{y_e} \partial_{xy}c(x,y) dy\ dx >0.
\end{align*}
Hence, $c$ is $\ominus$-monotone by \cref{e:5.6e3}.

\item Similar to (ii), the cost 
\begin{align*}
c\colon \mathbb{R}\times \mathbb{R}\to \left]-\infty, +\infty\right] \colon (x,y) \mapsto  \begin{cases}
-\ln(xy-1), &\text{if $y<\frac{1}{x}<0$;}\\
\pinf, &\text{otherwise.}
\end{cases}
\end{align*}
is $\ominus$-monotone. One can check this by mimicking the reasoning 
in (ii).
\end{enumerate}
\end{example}

\begin{example}[$\oplus$-monotone costs] \label{ex:5.12}\ 
\begin{enumerate}
\item Consider $c\colon \mathbb{R}\times \mathbb{R}\to \mathbb{R}\colon (x,y)\mapsto -xy$. Then $\partial_{xy}c(x,y)=-1$ for all $(x,y)\in \mathbb{R}\times \mathbb{R}$. Take $(x_s,y_s), (x_e,y_e)\in \mathbb{R}\times \mathbb{R}$ such that $(x_s,y_s)\prec^\ominus (x_e,y_e)$ and 
\begin{align*}
\int_{x_s}^{x_e}\int_{y_e}^{y_s} \partial_{xy}c(x,y) dy\ dx =-(x_e-x_s)(y_s-y_e)<0.
\end{align*}
Hence, $c$ is $\oplus$-monotone by \cref{e:5.6e4}. 

\item Consider 
\begin{align*}
c\colon \mathbb{R}\times \mathbb{R}\to \left]-\infty, +\infty\right] \colon (x,y) \mapsto  \begin{cases}
\frac{1}{y-x}, &\text{if $y>x$;}\\
\pinf, &\text{if $y\leq x$.}
\end{cases}
\end{align*}
Then $D=\{(x,y)\in \mathbb{R}\times \mathbb{R}\ |\ y>x\}$. For $(x_s, y_s), (x_e,y_e)\in D$ such that $(x_s, y_s)\prec^\ominus(x_e,y_e)$, we have $y_s>x_s, y_e>x_e$ and $x_s<x_e, y_e<y_s$, so $y_s>y_e>x_e>x_s$. Hence, $(x_s,y_e), (x_e,y_s)\in D$. Moreover, $\partial_{xy}c(x,y)=-\frac{2}{(y-x)^3}<0$ so 
\begin{align*}
\int_{x_s}^{x_e}\int_{y_e}^{y_s} \partial_{xy}c(x,y) dy\ dx <0.
\end{align*}
Hence, $c$ is $\oplus$-monotone by \cref{e:5.6e4}.

\item Now let $f: \mathbb{R} \to \left]-\infty, +\infty\right]$ be 
strictly convex and proper. Suppose that $f$ is differentiable on $\inte \dom(f)$. Consider 
\begin{align*}
c\colon \mathbb{R}\times \mathbb{R}\to \left]-\infty, +\infty\right] \colon (x,y) \mapsto  \begin{cases}
f(x)-f(y)-f'(y)(x-y), &\text{if $(x,y)\in \dom(f)\times \inte \dom(f)$;}\\
\pinf, &\text{otherwise,}
\end{cases}
\end{align*}
which is the Bregman cost \footnote{Recall \cref{ex:Bregman}.} on $\mathbb{R}\times\mathbb{R}$, and $D=\dom(f)\times \inte \dom(f)$.  Let $(x_s,y_s), (x_e,y_e)\in D$ such that $(x_s,y_s)\prec^\ominus (x_e,y_e)$. Then  $(x_s,y_e), (x_e,y_s)\in D$ as shown in \cref{ex:Bregman}. Moreover, $f$ is strictly convex, so  
\begin{align*}
c(x_s,y_s)-c(x_e,y_s)+c(x_e,y_e)-c(x_s,y_e)=\big(f'(y_s)-f'(y_e)\big)(x_e-x_s)>0.
\end{align*}
Therefore, $c$ is $\oplus$-monotone by \cref{d:monocost}.
\end{enumerate}
\end{example}

\subsection{$c$-path bounded extension for $\ominus$-monotone costs}
\label{sec:5.1}
In this subsection we provide results regarding $c$-path bounded extensions for $\ominus$-monotone costs.

\begin{lemma}
\label{l:3points-minus}
    Let $c: \mathbb{R}\times \mathbb{R}\to \left]-\infty, +\infty\right]$ be a cost that is $\ominus$-monotone and has convex domain $D$. 
    Suppose that $(G, \preccurlyeq^{\ominus})$ is a chain. For $(x_s, y_s), (x_m, y_m), (x_e, y_e)\in G$, assume that 
    \begin{enumerate}
        \item $c(x_{m}, y_s)<+\infty$, $c(x_e, y_m)<+\infty$, and 
        \item $\min_{\preccurlyeq^\ominus}\{(x_s, y_s), (x_e,y_e)\}\preccurlyeq^\ominus (x_m, y_m) \preccurlyeq^\ominus \max_{\preccurlyeq^\ominus}\{(x_s,y_s),(x_e, y_e)\}$ does not hold.
    \end{enumerate}
    Then $$c(x_s, y_s)-c(x_m, y_s)+c(x_m, y_m)-c(x_e, y_m)\leq c(x_s, y_s)-c(x_e, y_s).$$ 
\end{lemma}

\begin{proof}
    Without loss of generality, we assume $\min_{\preccurlyeq^\ominus}\{(x_s, y_s), (x_e,y_e)\}=(x_s,y_s)$ and $\max_{\preccurlyeq^\ominus}\{(x_s,y_s),(x_e, y_e)\}=(x_e,y_e)$. Note that (ii) implies $(x_m,y_m)$ is not in the box with vertices $(x_s, y_s), (x_e, y_s), (x_e, y_e), (x_s, y_e)$, which can be discussed by the following cases:\\

    \textbf{Case 1:} Suppose $x_m<x_s$ and $y_m=y_s$. Then we have 
    \begin{align*}
        &\quad \ \ c(x_s, y_s)-c(x_m, y_s)+c(x_m, y_m)-c(x_e, y_m)\\
        &=c(x_s, y_s)-c(x_m, y_s)+c(x_m, y_s)-c(x_e, y_s)\\
        &=c(x_s, y_s)-c(x_e, y_s). 
    \end{align*}

    \textbf{Case 2:} Suppose $x_m=x_e$ and $y_m<y_e$. Then we have 
    \begin{align*}
        &\quad \ \ c(x_s, y_s)-c(x_m, y_s)+c(x_m, y_m)-c(x_e, y_m)\\
        &=c(x_s, y_s)-c(x_e, y_s)+c(x_e, y_m)-c(x_e, y_m)\\
        &={\color{black}c(x_s, y_s)-c(x_e, y_s)}. 
    \end{align*}

    \textbf{Case 3:} Suppose $x_m=x_s$ and $y_m>y_s$. Then
    \begin{equation}
    \label{e:l4.2equ1}
    \begin{split}
        &\quad \ \ c(x_s, y_s)-c(x_m, y_s)+c(x_m, y_m)-c(x_e, y_m)\\
        &=c(x_s, y_s)-c(x_s, y_s)+c(x_s, y_m)-c(x_e, y_m)\\
        &=c(x_s, y_m)-c(x_e, y_m). 
    \end{split}
    \end{equation}
    By the assumption of this case, we have $(x_s, y_s)\preccurlyeq^\oplus (x_e, y_m)$. If $x_s=x_e$, then 
    \begin{align*}
        c(x_s, y_m)-c(x_e, y_m)=c(x_s, y_s)-c(x_e, y_s)=0. 
    \end{align*}
    If $x_s<x_e$, then $(x_s, y_s)\prec^\oplus (x_e, y_m)$. By (i), we know $(x_e, y_m)\in D$ and the monotonicity of the cost gives 
    \begin{equation}
    \label{e:l4.2equ2}
    \begin{split}
        &\quad \ \ c(x_s, y_s)-c(x_e, y_s)+c(x_e, y_m)-c(x_s, y_m)>0\\
        & \Rightarrow c(x_s, y_s)-c(x_e, y_s)>c(x_s, y_m)-c(x_e, y_m). 
    \end{split}
    \end{equation}
    Combining \cref{e:l4.2equ1} and \cref{e:l4.2equ2}, we have the desired result.\\

    \textbf{Case 4:} Suppose $x_m>x_e$ and $y_m=y_e$. Because $(x_s, y_s), (x_m, y_s)\in D$ and $D$ is convex, we know $(x_e, y_s)\in \dom(c)$ and obtain
    \begin{equation}
    \label{e:l4.2equ3}
    \begin{split}
        &\quad \ \ c(x_s, y_s)-c(x_m, y_s)+c(x_m, y_m)-c(x_e, y_m)\\
        &=c(x_s, y_s)-c(x_e, y_s)+c(x_e, y_s)-c(x_m, y_s)+c(x_m, y_e)-c(x_e, y_e).
    \end{split}
    \end{equation}
    By the assumption of this case, we know $(x_e, y_e)\preccurlyeq^\oplus(x_m, y_s)$. If $y_s=y_e$, then $y_s=y_m=y_e$ and 
    \begin{align*}
        &\quad \ \ c(x_s, y_s)-c(x_e, y_s)+c(x_e, y_s)-c(x_m, y_s)+c(x_m, y_e)-c(x_e, y_e)\\
        &=c(x_s, y_s)-c(x_e, y_s)+c(x_e, y_s)-c(x_m, y_s)+c(x_m, y_s)-c(x_e, y_s)\\
        &=c(x_s, y_s)-c(x_e, y_s).
    \end{align*}
    If $y_s>y_e$, then $(x_e, y_e)\prec^\oplus (x_m,y_s)$. Since the cost is $\ominus$-monotone, we have {\color{black} $c(x_e, y_e)-c(x_m, y_e)+c(x_m, y_s)-c(x_e, y_s)>0$. $(x_m, y_s)$ and $(x_e, y_s)$ are both in $D$, so $c(x_e, y_s)-c(x_m, y_s)<c(x_e,y_e)-c(x_m, y_e)$. Together with \cref{e:l4.2equ3}, we have
    \begin{align*}
        &\quad \ \ c(x_s, y_s)-c(x_m, y_s)+c(x_m, y_m)-c(x_e, y_m)\\
        &=c(x_s, y_s)-c(x_e, y_s)+c(x_e, y_s)-c(x_m, y_s)+c(x_m, y_e)-c(x_e, y_e)\\
        &<c(x_s, y_s)-c(x_e, y_s)+c(x_e, y_e)-c(x_m, y_e)+c(x_m, y_e)-c(x_e, y_e)\\
        &=c(x_s, y_s)-c(x_e, y_s).
    \end{align*}
}

    \textbf{Case 5:} Suppose $x_m<x_s$ and $y_m>y_s$. Because $(x_m, y_m), (x_e, y_m)\in D$ and $x_m<x_s\leq x_e$, we know $(x_s, y_m)\in D$ by convexity of $D$. Therefore, we have 
    \begin{equation}
        \label{e:l4.2equ4}
        \begin{split}
            &\quad \ \ c(x_s, y_s)-c(x_m, y_s)+c(x_m, y_m)-c(x_e, y_m)\\
        &=c(x_s, y_s)-c(x_m, y_s)+c(x_m, y_m)-c(x_s, y_m)+c(x_s, y_m)-c(x_e, y_m).
        \end{split}
    \end{equation}
    If $x_s=x_e$, then $c(x_s, y_s)-c(x_e, y_s)=0$ and (\ref{e:l4.2equ4}) becomes 
    \begin{align*}
        &\quad \ \ c(x_s, y_s)-c(x_m, y_s)+c(x_m, y_m)-c(x_s, y_m)+c(x_s, y_m)-c(x_e, y_m)\\
        &= c(x_s, y_s)-c(x_m, y_s)+c(x_m, y_m)-c(x_s, y_m)+c(x_s, y_m)-c(x_s, y_m)\\
        & = c(x_s, y_s)-c(x_m, y_s)+c(x_m, y_m)-c(x_s, y_m)\\
        &<0=c(x_s, y_s)-c(x_e, y_s)
    \end{align*}
    where the inequality is because $(x_m, y_s)\prec^\oplus(x_s, y_m)$ and $\ominus$-monotonicity of the cost. If $x_s<x_e$, then $(x_m, y_s)\prec^\oplus(x_s, y_m)$ and $(x_s, y_s)\prec^\oplus (x_e, y_m)$. By monotonicity of the cost, we have $c(x_m, y_s)-c(x_s, y_s)+c(x_s, y_m)-c(x_m, y_m)>0$ and $c(x_s, y_s)-c(x_e, y_s)+c(x_e, y_m)-c(x_s, y_m)>0$. Hence, (\ref{e:l4.2equ4}) yields 
    \begin{align*}
        &\quad \ \ c(x_s, y_s)-c(x_m, y_s)+c(x_m, y_m)-c(x_s, y_m)+c(x_s, y_m)-c(x_e, y_m)\\
        &<c(x_s, y_m)-c(x_e, y_m)\\
        &<c(x_s, y_s)-c(x_e, y_s).
    \end{align*}

    \textbf{Case 6:} Suppose $x_m>x_e$ and $y_m<y_e$. Because $(x_s, y_s), (x_m, y_s)\in D$ and $x_s\leq x_e < x_m$, by convexity of $D$, we know $(x_e, y_s)\in \dom(c)$. Hence,
    \begin{equation}
        \label{e:l4.2equ5}
        \begin{split}
             &\quad \ \ c(x_s, y_s)-c(x_m, y_s)+c(x_m, y_m)-c(x_e, y_m)\\
        &=c(x_s, y_s)-c(x_e, y_s)+c(x_e, y_s)-c(x_m, y_s)+c(x_m, y_m)-c(x_e, y_m)\\
        &<c(x_s, y_s)-c(x_e, y_s)
        \end{split}
    \end{equation}
    where the last inequality comes from $(x_e, y_m)\prec^\oplus(x_m, y_s)$ and monotonicity of the cost. 
\end{proof}

\begin{proposition}
\label{p:charat-minus}
    Suppose that the cost $c: \mathbb{R}\times \mathbb{R}\to \left]-\infty, +\infty\right]$ is $\ominus$-monotone and has convex domain $D$, 
    and that $G$ is a nonempty subset of $D$. Assume that 
    \begin{enumerate}
        \item $(G, \preccurlyeq^\ominus)$ is a chain,
        \item $\big((x_i,y_i)\big)_{i=1}^{N+1}$ is a finite sequence in $G$, where $N\geq 2$, and 
        \item For all $2\leq j\leq N$, {\color{black} we have the following:
        \begin{align*}
        \text{$\min_{\preccurlyeq^\ominus}\{(x_1, y_1), (x_{N+1}, y_{N+1})\}\preccurlyeq^\ominus (x_j, y_j)\preccurlyeq^\ominus\max_{\preccurlyeq^\ominus} \{(x_1, y_1), (x_{N+1}, y_{N+1})\}$ does not hold. }
        \end{align*}}
    \end{enumerate}
    Then 
    \begin{equation}
        \label{e:goal}
        \sum_{i=1}^{N} \big(c(x_i,y_i)-c(x_{i+1}, y_i)\big)\leq c(x_1,y_1)-c(x_{N+1}, y_1) .
    \end{equation}
\end{proposition}
\begin{proof}
    {\color{black}Without loss of generality, we assume that $(x_i, y_i)\neq (x_{i+1}, y_{i+1})$ for $\forall i\in \{1, \ldots, N\}$ because consecutive duplicate points in a path do not contribute to the summation in \cref{e:goal}.} Moreover, we assume that $c(x_{i+1}, y_i)<+\infty$ for all $i=1, \ldots, N$. Otherwise, \cref{e:goal} holds trivially. For $N=2$, \cref{e:goal} is true by \cref{l:3points-minus}. Assume \cref{e:goal} is true for $N\geq 2$. For $N+1$, assumption (iii) implies that there is $2\leq j\leq N+1$ such that {\color{black} $\min_{\preccurlyeq^\ominus}\{(x_{j-1}, y_{j-1}), (x_{j+1}, y_{j+1})\} \preccurlyeq^\ominus (x_j, y_j) \preccurlyeq^\ominus \max_{\preccurlyeq^\ominus}\{(x_{j-1}, y_{j-1}), (x_{j+1}, y_{j+1})\}$ does not hold} \footnote{{\color{black} Indeed, if $\min_{\preccurlyeq^\ominus}\{(x_{j-1}, y_{j-1}), (x_{j+1}, y_{j+1})\} \preccurlyeq^\ominus (x_j, y_j) \preccurlyeq^\ominus \max_{\preccurlyeq^\ominus}\{(x_{j-1}, y_{j-1}), (x_{j+1}, y_{j+1})\}$ for all $2\leq j \leq N+1$, then by the assumption that there is no consecutive duplicate points and a simple induction argument, we have $\min_{\preccurlyeq^\ominus}\{(x_1, y_1), (x_{N+2}, y_{N+2})\}\preccurlyeq^\ominus (x_j, y_j)\preccurlyeq^\ominus\max_{\preccurlyeq^\ominus} \{(x_1, y_1), (x_{N+2}, y_{N+2})\}$, which contradicts assumption (iii).}}. By \cref{l:3points-minus} and the induction assumption for $N$, we have 
    \begin{align*}
        &\quad \ \ \sum_{i=1}^{N+1} \big(c(x_i,y_i)-c(x_{i+1}, y_i)\big)\\
        &=\sum_{i=1}^{j-2} \big(c(x_i,y_i)-c(x_{i+1}, y_i)\big)+c(x_{j-1}, y_{j-1})-c(x_j, y_{j-1})+c(x_j, y_j)-c(x_{j+1}, y_{j})\\
        &\quad +\sum_{k=j+1}^{N+1} \big(c(x_k,y_k) -c(x_{k+1}, y_k)\big)\\
        &\leq \sum_{i=1}^{j-2} \big(c(x_i,y_i)-c(x_{i+1}, y_i)\big)+c(x_{j-1}, y_{j-1})-c(x_{j+1}, y_{j-1})+\sum_{k=j+1}^{N+1} \big(c(x_k,y_k) -c(x_{k+1}, y_k)\big) \tag{by \cref{l:3points-minus}}\\
        &\leq c(x_1, y_1)-c(x_{N+2}, y_1) \tag{by the induction assumption for $N$}. 
    \end{align*}
\end{proof}

\begin{theorem}
\label{c:charaofcyc}
    Let $c: \mathbb{R}\times \mathbb{R}\to \left]-\infty, +\infty\right]$ be a cost function that is $\ominus$-monotone and has convex domain $D$. 
    Suppose that $G$ is a nonempty subset of $D$. Then 
\begin{align*}
\text{$G$ is $c$-cyclically monotone} \Leftrightarrow \text{$(G, \preccurlyeq^\ominus)$ is a chain}.
\end{align*}
\end{theorem}
\begin{proof}
We only show the proof for $\ominus$-monotone cost, and $\oplus$-monotone case can be proved similarly. 

    "$\Rightarrow$": Suppose $G$ is $c$-cyclically monotone. Then it is monotone, i.e., 
    \begin{align*}
        c(x, y)-c(u,y)+c(u,v)-c(x,v)\leq 0
    \end{align*}
    for any $(x,y), (u,v)\in G$. Because $c$ is $\ominus$-monotone, we know $(x,y)\preccurlyeq^\ominus (u,v)$ or $(u,v)\preccurlyeq^\ominus(x,y)$. 

    \par "$\Leftarrow$": It is enough to show that for any $\big((x_i,y_i)\big)_{i=1}^{N+1}$ in $G$ with $(x_1, y_1)=(x_{N+1},y_{N+1})$
    \begin{equation}
        \label{e:cor4.8equ1}
        \sum_{i=1}^N \big(c(x_i,y_i)-c(x_{i+1}, y_i)\big)\leq 0.
    \end{equation}
    {\color{black} Let $\{i_k\}_{k=1}^M$ be the subset of indices $\{1, \ldots, N+1\}$ such that $(x_{i_k},y_{i_k})=(x_1,y_1)=(x_{N+1}, y_{N+1})$ and $i_k<i_{k+1}$ for $\forall k\in \{1, \ldots, M-1\}$.}  Note that if $i_{k+1}-i_k\geq 2$, then $\sum_{j=i_k}^{i_{k+1}-1}(c(x_j, y_j)-c(x_{j+1},y_j))\leq c(x_{i_k},y_{i_k})-c(x_{i_{k+1}},y_{i_k})=0$ by \cref{p:charat-minus}. Hence, 
    \begin{align*}
        \sum_{i=1}^{N} \big(c(x_i, y_i)-c(x_{i+1}, y_i)\big)= \sum_{k=1}^M \sum_{j=i_k}^{i_{k+1}-1}\big(c(x_j, y_j)-c(x_{j+1}, y_j)\big)\leq 0.
       \end{align*}
\end{proof}

\begin{defn}[chain extension]
\label{d:completion}
    Suppose $(G, \preccurlyeq^\ominus)$ is a chain. For $(x,y)\in G$, set 
    \begin{align*}
        A^G_{(x,y)}:=\{(u,v)\in G\ |\ (x,y)\preccurlyeq^\ominus(u,v)\}. 
    \end{align*}
    If $G\smallsetminus A^G_{(x,y)}\neq \varnothing$, i.e., $(x,y)$ is not the minimal element of $G$, we set 
    \begin{align*}
        \X^G_{(x,y)}:=\sup P_X(G\smallsetminus A^G_{(x,y)});\\
        \Y^G_{(x,y)}:=\inf P_Y(G\smallsetminus A^G_{(x,y)}).
    \end{align*}
    Otherwise, set 
    \begin{align*}
         \X^G_{(x,y)}:=x;\\
        \Y^G_{(x,y)}:=y.
    \end{align*}
    Moreover, set
    \begin{align*}
        \mathcal{I}^G_{(x,y)}:=\{(1-\lambda)(\X^G_{(x,y)}, \Y^G_{(x,y)})+\lambda (x,y)\ |\ \lambda\in [0,1]\}.
    \end{align*}
    The set $\com(G):=\cup_{(x,y)\in G}\mathcal{I}^G_{(x,y)}$ is said to be the \emph{chain extension} of $G$.  
\end{defn}

\begin{remark}
\label{r:completion}
    Note that if $(G, \preccurlyeq^\ominus)$ is a chain, then: 
    \begin{enumerate}
        \item $(\X^G_{(x,y)}, \Y^G_{(x,y)}) \in \overline{G}$ for every $(x,y)\in G$;
        \item $(u,v)\preccurlyeq^\ominus(\X^G_{(x,y)}, \Y^G_{(x,y)})$ for every $(x,y)\in G$ with $A^G_{(x,y)}\neq G$ and every $(u,v)\in G\smallsetminus A^G_{(x,y)}$;
        \item $(\mathcal{I}^G_{(x,y)}, \preccurlyeq^\ominus)$ is a chain for every $(x,y)\in G$;
        \item $(a,b)\in \mathcal{I}^G_{(x,y)}$ implies $(\X^G_{(x,y)}, \Y^G_{(x,y)})\preccurlyeq^\ominus(a,b) \preccurlyeq^\ominus(x,y)$;
        \item $(x_1,y_1)\preccurlyeq^\ominus (x_2,y_2)$ and $(x_1, y_1)\neq (x_2,y_2)$ imply $(a_1,b_1)\preccurlyeq^\ominus (a_2,b_2)$ for $(a_1,b_1)\in \mathcal{I}^G_{(x_1,y_1)}, (a_2,b_2)\in \mathcal{I}^G_{(x_2,y_2)}$. 
    \end{enumerate} 
\end{remark}

\begin{proposition}
    \label{p:completion}
    Let $(G, \preccurlyeq^\ominus)$ be a chain and $\com(G)$ be the chain extension of $G$. Then
    \begin{enumerate}
        \item $(\com(G), \preccurlyeq^\ominus)$ is a chain;
        \item If $G$ has the minimal element, then $\com(G)$ has the minimal element and $\min_{\preccurlyeq^\ominus}\com(G)=\min_{\preccurlyeq^\ominus}G$;
        \item If $G$ has the maximal element, then $\com(G)$ has the maximal element and $\max_{\preccurlyeq^\ominus}\com(G)=\max_{\preccurlyeq^\ominus}G$.
    \end{enumerate}
    If we further assume that $G$ is closed, then $\com(G)$ has the 
    following additional properties:
   \begin{enumerate}
   \setcounter{enumi}{3}
       \item $\com(G)$ is closed;
       \item $(\X^{\com(G)}_{(x,y)}, \Y^{\com(G)}_{(x,y)})=(x,y)$ for every $(x,y)\in \com(G)$.
   \end{enumerate}
\end{proposition}

\begin{proof}
    (i): Take $(a_1,b_1), (a_2,b_2)\in \com(G)$. Then $(a_1,b_1)\in \mathcal{I}^G_{(x_1, y_1)}$ and $(a_2, b_2)\in \mathcal{I}^G_{(x_2,y_2)}$ for some $(x_1, y_1), (x_2, y_2)\in G$. If $(x_1, y_1)=(x_2,y_2)$, then $(a_1,b_1)\preccurlyeq^\ominus(a_2,b_2)$ or $(a_2, b_2)\preccurlyeq^\ominus (a_1,b_1)$ because $(\mathcal{I}^G_{(x_1, y_1)}, \preccurlyeq^\ominus)$ is a chain\footnote{See \cref{r:completion}.(iii).}. If $(x_1, y_1)\neq (x_2,y_2)$, then we assume WLOG that $(x_1, y_1)\preccurlyeq^\ominus(x_2,y_2)$ because $(G, \preccurlyeq^\ominus)$ is a chain. By \cref{r:completion}.(v), we have $(a_1, b_1)\preccurlyeq^\ominus (a_2, b_2)$.

    \par (ii): Assume $G$ has the minimal element $\min_{\preccurlyeq^\ominus}G$. Take $(a,b)\in \com(G)$, i.e., $(a,b)\in \mathcal{I}^G_{(x,y)}$ for some $(x,y)\in G$. If $(x,y)=\min_{\preccurlyeq^\ominus} G$, then $G\smallsetminus A^G_{(x,y)}=\varnothing$ and $(\X^G_{(x,y)}, \Y^G_{(x,y)})=(x,y)=(a,b)$ \footnote{Recall \cref{d:completion}.}. If $(x,y)\neq \min_{\preccurlyeq^\ominus} G$ then $\min_{\preccurlyeq^\ominus}G \preccurlyeq^\ominus(\X^G_{(x,y)}, \Y^G_{(x,y)})\preccurlyeq^\ominus(a,b)\preccurlyeq^\ominus(x,y)$ by \cref{r:completion}.(ii)\&(iv). Therefore, $\min_{\preccurlyeq^\ominus}G\preccurlyeq^\ominus(a,b)$ for all $(a,b)\in \com(G)$. Because $\min_{\preccurlyeq^\ominus}G\in G\subseteq \com(G)$ and $(\com(G), \preccurlyeq^\ominus)$ is a chain, we know $\min_{\preccurlyeq^\ominus}\com(G)=\min_{\preccurlyeq^\ominus}G$.
    
    \par (iii): Assume $G$ has the maximal element $\max_{\preccurlyeq^\ominus}G$. Take $(a,b)\in \com(G)$, i.e., $(a,b)\in \mathcal{I}^G_{(x,y)}$ for some $(x,y)\in G$. By \cref{r:completion}.(iv), we know $(a,b)\preccurlyeq^\ominus(x,y)\preccurlyeq^\ominus\max_{\preccurlyeq^\ominus}G$. Because $(\com(G), \preccurlyeq^\ominus~)$ is a chain and $\max_{\preccurlyeq^\ominus}G\in G \subseteq \com(G)$, we know $\max_{\preccurlyeq^\ominus} \com(G)=\max_{\preccurlyeq^\ominus} G$. 

    \par (iv): Assume that $(\overline{x}, \overline{y})\in \overline{\com(G)}\smallsetminus G$. By closedness of $G$, $\mathcal{B}_{\delta}(\overline{x}, \overline{y})\cap G=\varnothing$ for some $\delta>0$. Since $(\overline{x}, \overline{y}) \in \overline{\com(G)}$, we can take $(u,v)\in \mathcal{B}_\delta (\overline{x}, \overline{y})\cap \com(G)$.  Then $(u,v)\in \mathcal{I}^G_{(x,y)}$ for some $(x,y)\in G$ and thus $(\X^G_{(x,y)},\Y^G_{(x,y)})\preccurlyeq^\ominus(u,v)\preccurlyeq^\ominus(x,y)$ \footnote{Recall \cref{r:completion}.(iv).}. Note that 
    \begin{align*}
    \{(\X^G_{(x,y)}, \Y^G_{(x,y)}), (u,v), (\overline{x}, \overline{y}), (x,y)\}\subseteq (\com(G)\cup \{(\overline{x}, \overline{y})\})\subseteq \overline{\com(G)}
    \end{align*}
    and $\overline{\com(G)}$ is a chain with respect to $\preccurlyeq^\ominus$ \footnote{See \cref{r:cluster}.(i).}, so $(\{(\X^G_{(x,y)}, \Y^G_{(x,y)}), (u,v), (\overline{x}, \overline{y}), (x,y)\}, \preccurlyeq^\ominus)$ is a chain. On the other hand, $(\X^G_{(x,y)}, \Y^G_{(x,y)})$ and $(x,y)$ are both in $G$ \footnote{Recall \cref{r:completion}.(i).}, so $\|(\overline{x}, \overline{y})-(\X^G_{(x,y)},\U^G_{(x,y)})\|\geq \delta$ and $\|(\overline{x},\overline{y})-(x,y)\|\geq \delta$ \footnote{Recall that $\mathcal{B}_\delta (\overline{x}, \overline{y})\cap G=\varnothing$.}. Recall that $(\X^G_{(x,y)},\Y^G_{(x,y)})\preccurlyeq^\ominus(u,v)\preccurlyeq^\ominus(x,y)$, so \cref{r:cluster}.(ii) implies 
    \begin{equation}
    \label{e:5.15e1}
    (\X^G_{(x,y)}, \Y^G_{(x,y)})\preccurlyeq^\ominus (\overline{x}, \overline{y})\preccurlyeq^\ominus (x,y) \footnote{Suppose the relation does not hold. Then there are two cases: (i) $(\overline{x}, \overline{y})\preccurlyeq^\ominus(\X^G_{(x,y)}, \Y^G_{(x,y)})\preccurlyeq^\ominus (x,y)$ or (ii) $(\X^G_{(x,y)}, \Y^G_{(x,y)})\preccurlyeq^\ominus (x,y)\preccurlyeq^\ominus (\overline{x}, \overline{y})$. If (i) holds, then $(\overline{x}, \overline{y})\preccurlyeq^\ominus(\X^G_{(x,y)}, \Y^G_{(x,y)})\preccurlyeq^\ominus(u,v)\preccurlyeq^\ominus (x,y)$ thus $\|(\overline{x}, \overline{y})-(\X^G_{(x,y)}, \Y^G_{(x,y)})\|\leq \|(\overline{x}, \overline{y})-(u,v)\|<\delta$ by \cref{r:cluster}.(ii),  which contradicts $\|(\overline{x}, \overline{y})-(\X^G_{(x,y)},\U^G_{(x,y)})\|\geq \delta$. If (ii) holds, then $(\X^G_{(x,y)}, \Y^G_{(x,y)})\preccurlyeq^\ominus(u,v)\preccurlyeq^\ominus (x,y)\preccurlyeq^\ominus(\overline{x}, \overline{y})$ thus $\|(\overline{x}, \overline{y})-(x,y)\|\leq \|(\overline{x}, \overline{y})-(u, v)\|<\delta$ by \cref{r:cluster}, which contradicts $\|(\overline{x},\overline{y})-(x,y)\|\geq \delta$.}.
    \end{equation}
Now, we claim that $\mathcal{B}_\delta(\overline{x}, \overline{y})\cap \com(G)\subseteq \mathcal{I}^G_{(x,y)}$. Indeed, suppose there is $(a,b)\in \mathcal{B}_\delta(\overline{x}, \overline{y})\cap \com(G)$ and $(a,b)\notin \mathcal{I}^G_{(x,y)}$. Then $(a,b)\in \mathcal{I}^G_{(x',y')}$ for some $(x',y')\in G\smallsetminus \{(x,y)\}$. 

\textbf{Case 1:} Assume $(x',y')\preccurlyeq^\ominus (x,y)$. Then together with \cref{e:5.15e1}, we have 
    \begin{equation}
    (a,b)\preccurlyeq^\ominus (x',y')\preccurlyeq^\ominus (\X^G_{(x,y)}, \Y^G_{(x,y)}) \preccurlyeq^\ominus (\overline{x}, \overline{y}) \preccurlyeq^\ominus (x,y) \footnote{Recall \cref{r:completion}.(iv)\&(v).}.
    \end{equation}
    By \cref{r:cluster}.(ii)
    , $\|(\overline{x}, \overline{y})-(\X^G_{(x,y)}, \Y^G_{(x,y)})\|\leq \|(a,b)-(\overline{x}, \overline{y})\|<\delta$, which implies $(\X^G_{(x,y)}, \Y^G_{(x,y)})\in \mathcal{B}_\delta(\overline{x}, \overline{y})$. However, it is absurd because $(\X^G_{(x,y)}, \Y^G_{(x,y)})\in G$ and $\mathcal{B}_\delta(\overline{x}, \overline{y})\cap G=\varnothing$. 

\textbf{Case 2:} Assume that $(x,y)\preccurlyeq^\ominus (x',y')$. Then together with \cref{e:5.15e1}, we have
\begin{equation}
(\X^G_{(x,y)}, \Y^G_{(x,y)}) \preccurlyeq^\ominus (\overline{x}, \overline{y}) \preccurlyeq^\ominus (x,y)\preccurlyeq^\ominus (a,b) \preccurlyeq^\ominus (x', y') \footnote{Recall 
\cref{r:completion}.(iv)\&(v).}.
\end{equation}
By \cref{r:cluster}.(ii), $\|(\overline{x}, \overline{y})-(x,y)\|\leq \|(\overline{x}, \overline{y}-(a,b)\|<\delta$, which implies $(x,y)\in \mathcal{B}_\delta(\overline{x}, \overline{y})$. However, it is absurd because $(x,y)\in G$ and $\mathcal{B}_\delta(\overline{x}, \overline{y})\cap G=\varnothing$. 

Therefore, our claim is correct, i.e., $\mathcal{B}_\delta (\overline{x}, \overline{y})\cap \com(G) \subseteq \mathcal{I}^G_{(x,y)}$. 
    
Now, take a sequence $\big((x_n,y_n)\big)_{n\in \mathbb{N}}$ in $(\mathcal{B}_\delta (\overline{x}, \overline{y})\cap \com(G))\subseteq \mathcal{I}^G_{(x,y)}$ such that $(x_n, y_n)\to (\overline{x}, \overline{y})$. Because $\mathcal{I}^G_{(x,y)}$ is closed, we know $(\overline{x}, \overline{y})\in \mathcal{I}^G_{(x,y)}\subseteq \com(G)$. 

    \par (v): Take $(u,v)\in \com(G)$. Then $(u,v)\in \mathcal{I}^G_{(x,y)}$ for some $(x,y)\in G$, i.e., $(u,v)=(1-\lambda) (\X^G_{(x,y)}, \Y^G_{(x,y)})+\lambda (x,y)$ for some $\lambda\in [0,1]$. 
    
    \par \textbf{Case 1: }Suppose $(\X^G_{(x,y)}, \Y^G_{(x,y)})=(x,y)$. Then $(\X^{\com(G)}_{(x,y)}, \Y^{\com(G)}_{(x,y)})=(\X^G_{(x,y)}, \Y^G_{(x,y)})=(x,y)=(u,v)$ because $(\X^G_{(x,y)}, \Y^G_{(x,y)})\preccurlyeq^\ominus(\X^{\com(G)}_{(x,y)}, \Y^{\com(G)}_{(x,y)})\preccurlyeq^\ominus(x,y)$. Hence, $(\X^{\com(G)}_{(u,v)}, \Y^{\com(G)}_{(u,v)})=(u,v)$.
    
    \par \textbf{Case 2: }Suppose $(\X^G_{(x,y)}, \Y^G_{(x,y)})\neq (x,y)$ and $(u,v)=(1-\lambda)(\X^G_{(x,y)}, \Y^G_{(x,y)})+\lambda (x,y)$ for some $\lambda \in \left]0, 1\right]$. Then 
     \begin{equation}
     \label{e:converge}
         (1-t)(\X^G_{(x,y)}, \Y^G_{(x,y)})+t(x,y)\to (u,v)
     \end{equation}
     for $t\in \left]0, \lambda\right[$ such that $t\to \lambda$. Note that $(1-t)(\X^G_{(x,y)}, \Y^G_{(x,y)})+t(x,y)\in \com(G) \smallsetminus A^{\com(G)}_{(u,v)}$ for $t\in \left] 0, \lambda\right[$, so \cref{e:converge} implies $(\X^{\com(G)}_{(u,v)}, \Y^{\com(G)}_{(u,v)})=(u,v)$. 
     
     \par \textbf{Case 3: }Suppose $(\X^G_{(x,y)}, \Y^G_{(x,y)})\neq (x,y)$ and $(u,v)=(\X^G_{(x,y)}, \Y^G_{(x,y)})$.  Then we set $(\tilde{x}, \tilde{y}):=(\X^G_{(x,y)}, \Y^G_{(x,y)})$ and we know $(u,v)\in \mathcal{I}^{G}_{(\tilde{x}, \tilde{y})}$ by the assumption of this case. If $(\X^G_{(\tilde{x},\tilde{y})}, \Y^G_{(\tilde{x},\tilde{y})})= (\tilde{x}, \tilde{y})$, then we deduce that $(\X^{\com(G)}_{(u,v)}, \Y^{\com(G)}_{(u,v)})=(u,v)$ by the same reasoning of \textbf{Case 1}. If $(\X^G_{(\tilde{x},\tilde{y})}, \Y^G_{(\tilde{x},\tilde{y})})\neq (\tilde{x}, \tilde{y})$, then $(u,v)\neq (\X^G_{(\tilde{x}, \tilde{y})}, \Y^G_{(\tilde{x}, \tilde{y})})$ and we deduce that $(\X^{\com(G)}_{(u,v)}, \Y^{\com(G)}_{(u,v)})=(u,v)$ by the same reasoning of \textbf{Case 2}.
\end{proof}

We are now ready for the main result of this section.
The main idea of its proof is simple: We ``fill in'' 
the original graph with line segments to obtain a better extension. 
In the context of \cref{ex:5.1}, we visualize 
in Figure~\ref{fig:2} 
the extension $\widetilde{G}$ of the original graph 
$G$ from \cref{ex:5.1}. See \cref{ex:5.20} for details. 
}

\begin{theorem}
    \label{t:cpathextension}
    Suppose that $c \colon \mathbb{R}\times \mathbb{R} \to \left]-\infty, +\infty\right]$ satisfies the following: 
    \begin{enumerate}
\item $D$ is a nonempty open convex proper subset of $\mathbb{R}\times \mathbb{R}$;
\item $c$ is $\ominus$-monotone;
\item $\forall (x,y)\in D$, there is an open ball $\mathcal{B}_{\delta}(x,y)$ such that $\sup_{(u,v)\in \mathcal{B}_{\delta}(x,y)}c(u,v)<+\infty$; 
\item $\forall (x,y)\in \bd(D)$,  $c(x_k,y_k)\to +\infty$ if $(x_k,y_k) \to (x,y)$.
\end{enumerate}
Moreover, suppose that $G$ is a nonempty subset of $D$ which is 
\begin{enumerate}
   \setcounter{enumi}{4}
       \item $c$-path bounded;
       \item for any partition $\{S_1, S_2\}$ of $G$, there are $(x,y)\in S_1, (u,v)\in S_2$ such that they are strongly connected or unilaterally connected \footnote{Recall \cref{d:connected}.}.
   \end{enumerate}
Then $G$ has a strongly connected $c$-path bounded extension, namely $\widetilde{G}:=\com(\overline{G})\cap D$, and $G$ is $c$-potentiable. 
\end{theorem}
\begin{proof}
    $G$ is $c$-path bounded, so $G$ is $c$-cyclically monotone by \cref{f:0707b}. By \cref{c:charaofcyc}, we know $(G, \preccurlyeq^\ominus)$ is a chain and thus $(\overline{G}, \preccurlyeq^\ominus)$ is a chain by \cref{r:cluster}.(i). Hence, $(\com(\overline{G}), \preccurlyeq^\ominus)$ is a chain by \cref{p:completion}.(i). We shall show that $\widetilde{G}:=\com(\overline{G}) \cap \dom(c)$ is a strongly connected $c$-path bounded extension of $G$ and $G$ is $c$-potentiable by the following steps: 

    \par \textbf{Step 1:}  \textbf{Show that if $(a,b)\in \overline{G}$ is neither the minimal nor the maximal element of $(\overline{G}, \preccurlyeq^\ominus)$, then $(a,b)\in D$}. By contradiction, suppose the opposite: There is a cluster point of $G$, say $(a,b)$ such that $(a,b)\notin D$ and $(a,b)$ is neither the minimal nor the maximal element of $(\overline{G}, \preccurlyeq^\ominus)$. Because $(a,b)$ is a cluster point of $G$, we know at least one of the following cases holds: 
    

     \par \textit{Case 1:} There is a sequence $\big((a_n,b_n)\big)_{n\in \mathbb{N}}$ in $A^{\overline{G}}_{(a,b)}\cap G$ such that $(a_n,b_n)\to (a,b)$.

     \par \textit{Case 2:} There is a sequence $\big((a_n,b_n)\big)_{n\in \mathbb{N}}$ in $(\overline{G}\smallsetminus A^{\overline{G}}_{(a,b)})\cap G$ such that $(a_n,b_n)\to (a,b)$.

     \par Assume \textit{Case 1} holds. Set $S_1:=A^{\overline{G}}_{(a,b)}\cap G\neq \varnothing$ and $S_2:=(\overline{G}\smallsetminus A^{\overline{G}}_{(a,b)})\cap G$. We will show that $\{S_1, S_2\}$ is a partition of $G$. Note that $S_1\cup S_2=(A^{\overline{G}}_{(a,b)}\cap G)\cup \big((\overline{G}\smallsetminus A^{\overline{G}}_{(a,b)})\cap G\big)=\big(A^{\overline{G}}_{(a,b)}\cup (\overline{G}\smallsetminus A^{\overline{G}}_{(a,b)})\big)\cap G=\overline{G}\cap G=G$. Moreover, because $(a,b)$ is not the minimal element of $(\overline{G}, \preccurlyeq^\ominus)$, there is $(u,v)\in (\overline{G}\smallsetminus A^{\overline{G}}_{(a,b)})$. Set $\delta:=\frac{\|(a,b)-(u,v)\|}{2}>0$ and because $(u,v)\in \overline{G}$, we know there is $(x,y)\in \mathcal{B}_\delta (u,v)\cap G$. Therefore, we have $(u,v)\preccurlyeq^\ominus(a,b)$ \footnote{Recall that $(u,v)\in \overline{G}\smallsetminus A^{\overline{G}}_{(a,b)}$.}, $\mathcal{B}_\delta(u,v) \cap \mathcal{B}_\delta(a,b)=\varnothing$ and $(x,y)\in \mathcal{B}_\delta (u,v) \cap G \subseteq \mathcal{B}_\delta (u,v)\cap \overline{G}$. On the other hand, $(a, b), (u,v), (x,y)\in \overline{G}$ which is a chain, so by \cref{l:keylemma}, we know $(x,y)\preccurlyeq^\ominus(a,b)$. It implies $S_2\neq \varnothing$, hence, $\{S_1, S_2\}$ is a partition of $G$. By assumption (vi), we know that there are $(x_1,y_1)\in S_1, (x_2,y_2)\in S_2$ such that $(x_2,y_1)\in D$ or $(x_1, y_2)\in D$. Without loss of generality, we assume $(x_2, y_1)\in D$. Because $(x_1, y_1)\neq (a,b)$ \footnote{Recall that $(a,b)\notin D$ and $(x_1, y_1)\in S_1 \subseteq G \subseteq D$.}, we can find $\varepsilon>0$ such that $\mathcal{B}_{\varepsilon}(x_1,y_1)\cap \mathcal{B}_\varepsilon (a,b)=\varnothing$. On the other hand, $(x_1,y_1)\in A^{\overline{G}}_{(a,b)}$ so $(a,b)\preccurlyeq^\ominus(x_1,y_1)$. By \cref{l:keylemma}, we know $(a_n,b_n)\preccurlyeq^\ominus(x_1,y_1)$ for $(a_n, b_n)\in \mathcal{B}_\varepsilon (a,b)$. Recall that  $(x_2, y_2)\in \overline{G}\smallsetminus A^{\overline{G}}_{(a,b)}$, so $(x_2,y_2)\preccurlyeq^\ominus (a,b)$ and thus $(x_2, y_2)\preccurlyeq^\ominus(a, b)\preccurlyeq^\ominus(a_n,b_n)\preccurlyeq^\ominus(x_1,y_1)$ for all $(a_n, b_n)\in \mathcal{B}_\varepsilon (a,b)$. Hence, $x_2\leq a\leq a_n \leq x_1$ and $y_1\leq b_n\leq b \leq y_2$. Because $(x_1, y_1), (x_2,y_1)\in D$ and $D$ is convex, $(a, y_1)\in D$. Similarly, because $(x_2, y_1), (x_2,y_2)\in D$ and $D$ is convex, we know $(x_2, b)\in D$. On the other hand, assumption (iii) implies there are $T_1>0$ and $T_2>0$ such that $\sup c(\mathcal{B}_{\varepsilon_1}(a, y_1))\leq T_1$ for some $\varepsilon_1>0$ and $\sup c(\mathcal{B}_{\varepsilon_2}(x_2, b))\leq T_2$ for some $\varepsilon_2>0$, which implies $c(a_n, y_1)\leq T_1$ for $|a_n-a|<\varepsilon_1$ and $c(x_2, b_n)\leq T_2$ for $|b_n-b|<\varepsilon_2$. Now, set $r=\min \{\varepsilon, \varepsilon_1, \varepsilon_2\}$ and for all $(a_n, b_n)\in \mathcal{B}_r (a,b)$, we have 
     \begin{equation}
     \label{e:5.19e1}
         c(x_1, y_1)+c(a_n,b_n)-T_1-T_2 \leq c(x_1, y_1)-c(a_n, y_1)+c(a_n,b_n)-c(x_2,b_n).
     \end{equation}
     Because $(a_n, b_n)\to (a,b)\in \bd(D)$, we know $c(a_n,b_n)\to +\infty$. Hence, the right-hand side of \cref{e:5.19e1} tends to $+\infty$. It implies $F_G\big((x_1, y_1), (x_2,y_2)\big)=+\infty$ because $((x_1,y_1), (a_n,b_n), (x_2,y_2))\in P^G_{(x_1,y_1)\to (x_2,y_2)}$. Therefore, $G$ is not $c$-path bounded, which is absurd.
     
     \par Assume \textit{Case 2} holds. Set $S_1:=A^{\overline{G}}_{(a,b)}\cap G$ and $S_2:=(\overline{G}\smallsetminus A^{\overline{G}}_{(a,b)})\cap G \neq \varnothing$. We will show that $\{S_1, S_2\}$ is a partition of $G$. Note that $S_1\cup S_2=(A^{\overline{G}}_{(a,b)}\cap G)\cup \big((\overline{G}\smallsetminus A^{\overline{G}}_{(a,b)})\cap G\big)=\big(A^{\overline{G}}_{(a,b)}\cup (\overline{G}\smallsetminus A^{\overline{G}}_{(a,b)})\big)\cap G=\overline{G}\cap G=G$. Moreover, because $(a,b)$ is not the maximal element of $(\overline{G}, \preccurlyeq^\ominus)$, there is $(u,v)\in A^{\overline{G}}_{(a,b)}$ such that $(u,v)\neq (a,b)$. Set $\delta:=\frac{\|(a,b)-(u,v)\|}{2}>0$ and because $(u,v)\in \overline{G}$, we know there is $(x,y)\in \mathcal{B}_\delta (u,v)\cap G$. Therefore, we have $(a,b)\preccurlyeq^\ominus(u,v)$, $\mathcal{B}_\delta(a,b) \cap \mathcal{B}_\delta(u,v)=\varnothing$ and $(x,y)\in \mathcal{B}_\delta (u,v) \cap G \subseteq \mathcal{B}_\delta (u,v)\cap \overline{G}$. On the other hand, $(a, b), (u,v), (x,y)\in \overline{G}$ which is a chain, so by \cref{l:keylemma}, we know $(a,b)\preccurlyeq^\ominus(x,y)$. It implies $S_1\neq \varnothing$, hence, $\{S_1, S_2\}$ is a partition of $G$. By assumption (vi), there are $(x_1,y_1)\in S_1, (x_2,y_2)\in S_2$ such that $(x_2,y_1)\in D$ or $(x_1, y_2)\in D$. Without loss of generality, we assume $(x_2, y_1)\in D$. Because $(x_2, y_2)\neq (a,b)$ \footnote{Recall that $(a,b)\notin D$ and {\color{black} $(x_2, y_2)\in S_2 \subseteq G \subseteq D$.}}, we can find $\varepsilon>0$ such that $\mathcal{B}_{\varepsilon}(x_2,y_2)\cap \mathcal{B}_\varepsilon (a,b)=\varnothing$. On the other hand, $(x_2,y_2)\in \overline{G} \smallsetminus A^{\overline{G}}_{(a,b)}$ so $(x_2,y_2)\preccurlyeq^\ominus(a,b)$. By \cref{l:keylemma}, we know $(x_2,y_2)\preccurlyeq^\ominus(a_n,b_n)$ for $(a_n, b_n)\in \mathcal{B}_\varepsilon (a,b)$. Recall that  $(x_1, y_1)\in A^{\overline{G}}_{(a,b)}$, so $(a,b)\preccurlyeq^\ominus (x_1,y_1)$ and thus $(x_2, y_2)\preccurlyeq^\ominus(a_n, b_n)\preccurlyeq^\ominus(a,b)\preccurlyeq^\ominus(x_1,y_1)$ for all $(a_n, b_n)\in \mathcal{B}_\varepsilon (a,b)$. Hence, $x_2\leq a_n\leq a \leq x_1$ and $y_1\leq b\leq b_n \leq y_2$. Because $(x_1, y_1), (x_2,y_1)\in D$ and $D$ is convex, $(a, y_1)\in D$. Similarly, because $(x_2, y_1), (x_2,y_2)\in D$ and $D$ is convex, we know $(x_2, b)\in D$. On the other hand, assumption (iii) implies there are $T_1>0$ and $T_2>0$ such that $\sup c(\mathcal{B}_{\varepsilon_1}(a, y_1))\leq T_1$ for some $\varepsilon_1>0$ and $\sup c(\mathcal{B}_{\varepsilon_2}(x_2, b))\leq T_2$ for some $\varepsilon_2>0$, which implies $c(a_n, y_1)\leq T_1$ for $|a_n-a|<\varepsilon_1$ and $c(x_2, b_n)\leq T_2$ for $|b_n-b|<\varepsilon_2$. Now, set $r=\min \{\varepsilon, \varepsilon_1, \varepsilon_2\}$ and for all $(a_n, b_n)\in \mathcal{B}_r (a,b)$, we have 
     \begin{equation}
     \label{e:Ae1}
         c(x_1, y_1)+c(a_n,b_n)-T_1-T_2 \leq c(x_1, y_1)-c(a_n, y_1)+c(a_n,b_n)-c(x_2,b_n).
     \end{equation}
     Because $(a_n, b_n)\to (a,b)\in \bd(D)$, we know $c(a_n,b_n)\to +\infty$. Hence, the right-hand side of \cref{e:Ae1} tends to $+\infty$. It implies $F_G\big((x_1, y_1), (x_2,y_2)\big)=+\infty$ because $((x_1,y_1), (a_n,b_n), (x_2,y_2))\in P^G_{(x_1,y_1)\to (x_2,y_2)}$. Therefore, $G$ is not $c$-path bounded, which is absurd.
     
     \par Summing up, we conclude that if $(a,b)\in \overline{G}$ is neither the minimal nor the maximal element of $(\overline{G}, \preccurlyeq^\ominus)$, then $(a,b)\in D$.


     \par \textbf{Step 2:} By \cref{p:completion}.(i), we know $(\com(\overline{G}), \preccurlyeq^\ominus)$ is a chain. \textbf{We want to show $(a,b)\in D$ for all $(a,b)\in \com(\overline{G})$ that is neither the minimal nor the maximal element of $(\com(\overline{G}), \preccurlyeq^\ominus)$}. Let $(a,b)\in \com(\overline{G})$ and assume $(a,b)$ is neither the minimal nor the maximal element of $\com(\overline{G})$. Then $(a,b)\in \mathcal{I}^{\overline{G}}_{(x,y)}$ for some $(x,y)\in \overline{G}$. Moreover, because $\overline{G}$ is closed, $(\X^{\overline{G}}_{(x,y)}, \Y^{\overline{G}}_{(x,y)})\in \overline{G}\subseteq \overline{D}$ by \cref{r:completion}.(i). Now, we discuss by following cases:  

    \par \textit{Case 1:} Suppose $(\X^{\overline{G}}_{(x,y)}, \Y^{\overline{G}}_{(x,y)})$ is the minimal element of $(\overline{G}, \preccurlyeq^\ominus)$ and $(x,y)$ is not the maximal element of $(\overline{G}, \preccurlyeq^\ominus)$. By \cref{p:completion}.(ii), $(\X^{\overline{G}}_{(x,y)},\Y^{\overline{G}}_{(x,y)})=\min_{\preccurlyeq^\ominus} \com(\overline{G})=\min_{\preccurlyeq^\ominus} \overline{G}$. Recall that $(a,b)\neq \min_{\preccurlyeq^\ominus} \com(\overline{G})$ and $(a,b)\preccurlyeq^\ominus (x,y)$, so $(x,y)\neq \min_{\preccurlyeq^\ominus}\com(\overline{G})=\min_{\preccurlyeq^\ominus} \overline{G}$. Hence, $(x,y)$ is neither the minimal nor the maximal element of $\overline{G}$. By \textbf{Step 1}, $(x,y)\in D=\inte (D)$ \footnote{Recall that $D$ is open by assumption (i).} and we have $(a,b)\in D$ by convexity of $D$.  
     
     \par \textit{Case 2:} Suppose $(\X^{\overline{G}}_{(x,y)}, \Y^{\overline{G}}_{(x,y)})$ is not the minimal element of $(\overline{G}, \preccurlyeq^\ominus)$ and $(x,y)$ is the maximal element of $(\overline{G}, \preccurlyeq^\ominus)$. By \cref{p:completion}.(iii), $(x,y)=\max_{\preccurlyeq^\ominus} \com(\overline{G})=\max_{\preccurlyeq^\ominus} \overline{G}$. Recall that $(a,b)\neq \max_{\preccurlyeq^\ominus}\com(\overline{G})=\max_{\preccurlyeq^\ominus} \overline{G}$ and $(\X^{\overline{G}}_{(x,y)}, \Y^{\overline{G}}_{(x,y)})\preccurlyeq^\ominus (a,b)$, so $(\X^{\overline{G}}_{(x,y)}, \Y^{\overline{G}}_{(x,y)})\neq \max_{\preccurlyeq^\ominus} \overline{G}$. Hence, $(\X^{\overline{G}}_{(x,y)}, \Y^{\overline{G}}_{(x,y)})$ is neither the minimal nor the maximal element of $(\overline{G}, \preccurlyeq^\ominus)$ and $(\X^{\overline{G}}_{(x,y)}, \Y^{\overline{G}}_{(x,y)})\in D=\inte (D)$ by \textbf{Step 1}. Moreover, $(a,b)\neq (x,y)$ since $(a,b)$ is not the maximal element of $\com(\overline{G})$. By convexity of $D$, we know $(a,b)\in D$. 
     
     \par \textit{Case 3:} Suppose $(\X^{\overline{G}}_{(x,y)}, \Y^{\overline{G}}_{(x,y)})$ is not the minimal element of $(\overline{G}, \preccurlyeq^\ominus)$ and $(x,y)$ is not the maximal element of $(\overline{G}, \preccurlyeq^\ominus)$. Then both of them are in $D$, so $(a,b)\in D$ by convexity of $D$. 
     
    \par \textit{Case 4:} Suppose  $(\X^{\overline{G}}_{(x,y)}, \Y^{\overline{G}}_{(x,y)})$ and $(x,y)$ are the minimal element and the maximal element of $(\overline{G}, \preccurlyeq^\ominus)$ respectively. We may assume $(\X^{\overline{G}}_{(x,y)}, \Y^{\overline{G}}_{(x,y)})\neq (x,y)$ otherwise $\overline{G}=\{(\X^{\overline{G}}_{(x,y)}, \Y^{\overline{G}}_{(x,y)})\}=\{(x,y)\}=\{(a,b)\}=G\subseteq D$. If $(\X^{\overline{G}}_{(x,y)}, \Y^{\overline{G}}_{(x,y)})\neq (x,y)$ and they are the minimal and the maximal elements of $(\overline{G}, \preccurlyeq^\ominus)$, then $\overline{G}=\{(\X^{\overline{G}}_{(x,y)}, \Y^{\overline{G}}_{(x,y)}), (x,y)\}$ \footnote{Suppose there is another point $(\widetilde{x}, \widetilde{y})$ between $(\X^{\overline{G}}_{(x,y)}, \Y^{\overline{G}}_{(x,y)})$ and $(x,y)$. Since $A^{\overline{G}}_{(x,y)}=\{(x,y)\}\neq \overline{G}$, we know from \cref{r:completion}.(ii) that $(\widetilde{x}, \widetilde{y})\preccurlyeq^\ominus (\X^{\overline{G}}_{(x,y)}, \Y^{\overline{G}}_{(x,y)})$, which contradicts the assumption that $(\X^{\overline{G}}_{(x,y)}, \Y^{\overline{G}}_{(x,y)})$ is the minimal element of $(\overline{G}, \preccurlyeq^\ominus)$.} and thus $\overline{G}=G$. Therefore, $(a,b)\in D$ by convexity of $D$. 

     \par \textbf{Step 3:} \textbf{Show that any two points in $\com(\overline{G})\cap D$ are ball chain connected}. Take $(x_s, y_s), (x_e,y_e)\in \com(\overline{G}) \cap D$. Without loss of generality, assume $(x_s, y_s)\preccurlyeq^\ominus (x_e,y_e)$ and $(x_s, y_s)\neq (x_e, y_e)$. Set 
     \begin{align*}
S:=\menge{(x,y)\in \com(\overline{G})}{(x_s, y_s)\preccurlyeq^\ominus(x,y)\preccurlyeq^\ominus(x_e,y_e)}
\end{align*}
which is the intersection of the compact rectangle $[x_s, x_e]\times [y_e, y_s]$ and a closed set $\com(\overline{G})$ \footnote{Recall \cref{p:completion}.(iv).}, hence $S$ is compact. Moreover, if $(x,y)\in S\smallsetminus\{(x_s, y_s), (x_e, y_e)\}$, then $(x,y)$ is neither the minimal nor the maximal element of $\com(\overline{G})$ so $(x,y)\in D$ by \textbf{Step 2}. Therefore, $S$ is compact and $S\subseteq D$. Since $\mathbb{R}^2\smallsetminus D$ is closed, we know 
     \begin{align*}
         d:=\inf_{\substack{(x,y)\in S \\ (u,v)\in \mathbb{R}^2\smallsetminus D}} \|(x,y)-(u,v)\|>0. 
     \end{align*}
     By compactness of $S$, we have a finite open cover of $S$
     \begin{align*}
         \menge{\mathcal{B}_{d}(x_n, y_n)}{(x_n, y_n)\in S, n\in \{1, \ldots, N\}}. 
     \end{align*}
     Moreover, we may assume $(x_1, y_1)=(x_s,y_s), (x_N, y_N)=(x_e,y_e)$ and $(x_n, y_n)\preccurlyeq^\ominus (x_{n+1}, y_{n+1})$ for all $n\in \{1, \ldots, N-1\}$ because $(S, \preccurlyeq^\ominus)$ is a chain.  We claim $\mathcal{B}_d(x_n, y_n)\cap \mathcal{B}_d(x_{n+1}, y_{n+1})\cap S \neq \varnothing$ for all $n\in \{1, \ldots, N-1\}$. By contradiction, suppose there is $n\in \{1, \ldots, N-1\}$ such that $\mathcal{B}_d(x_n,y_n) \cap \mathcal{B}_d(x_{n+1}, y_{n+1})\cap S=\varnothing$. Set 
     \begin{align*}
         &u^-:=\sup P_X(\mathcal{B}_d(x_n, y_n) \cap S)\\
         &v^-:=\inf P_Y(\mathcal{B}_d(x_n, y_n) \cap S)\\
         &u^+:=\inf P_X(\mathcal{B}_d(x_{n+1}, y_{n+1}) \cap S)\\
         &v^+:=\sup P_Y(\mathcal{B}_d(x_{n+1}, y_{n+1}) \cap S). 
     \end{align*}
     Then $(u^-, v^-), (u^+, v^+)$ are in $\overline{\mathcal{B}_d(x_n,y_n)\cap S}$ and $\overline{\mathcal{B}_d(x_{n+1}, y_{n+1})\cap S}$ respectively and thus both in $S$.  By \cref{l:keylemma2}, we know 
\begin{equation}
\label{e:revision1}
\big(\forall k\leq n\big) \ \big(\forall (u,v)\in \mathcal{B}_d(x_k, y_k)\cap S\big) \ \big(\exists (u_n,v_n)\in  \mathcal{B}_d(x_n, y_n)\cap S\big) \quad (u,v)\preccurlyeq^\ominus (u_n,v_n)\preccurlyeq^\ominus (u^-,v^-)
\end{equation}
\begin{equation}
\label{e:revision2}
\begin{split}
&\big(\forall k\geq n+1\big) \ \big(\forall (u,v)\in \mathcal{B}_d(x_k, y_k)\cap S\big) \ \big(\exists (u_{n+1},v_{n+1})\in  \mathcal{B}_d(x_{n+1}, y_{n+1})\cap S\big)\\
&(u^+,v^+)\preccurlyeq^\ominus (u_{n+1},v_{n+1})\preccurlyeq^\ominus (u,v).
\end{split}
\end{equation}
This implies $(u,v)\preccurlyeq^\ominus(u^-,v^-)$ if $(u,v)\in \mathcal{B}_d(x_k, y_k)\cap S$ for $k\leq n$, and $(u^+, v^+)\preccurlyeq^\ominus (\alpha, \beta)$ if $(\alpha,\beta)\in \mathcal{B}_d(x_k, y_k)\cap S$ for $k\geq n+1$. Recall that we assume $\mathcal{B}_d(x_n, y_n)\cap \mathcal{B}_d(x_{n+1}, y_{n+1})\cap S=\varnothing$, so $(u,v)\preccurlyeq^\ominus(\alpha, \beta)$ for all {\color{black}$(u,v)\in \mathcal{B}_d(x_n, y_n)\cap S, (\alpha, \beta)\in \mathcal{B}_d(x_{n+1}, y_{n+1})\cap S$} by \cref{l:keylemma} and thus $(u^-, v^-)\preccurlyeq^\ominus(u^+, v^+)$. On the other hand, $(\com(\overline{G})\smallsetminus A^{\com(\overline{G})}_{(u^+,v^+)})\subseteq \menge{(u,v)\in \com(\overline{G})}{(u,v)\preccurlyeq^\ominus(u^-,v^-)}$ \footnote{Take $(u,v)\in \com(\overline{G})\smallsetminus A^{\com(\overline{G})}_{(u^+,v^+)}$. If $(u,v)\notin S$, then $(u,v)\preccurlyeq^\ominus (x_1,y_1)\preccurlyeq^\ominus(u^-,v^-)$. If $(u,v)\in S$, then $(u,v)\in \mathcal{B}_d (x_k,y_k)\cap S$ for some $k\in \{1, \ldots, n\}$ thus $(u,v)\preccurlyeq^\ominus(u^-, v^-)$.}, so $(\X^{\com(\overline{G})}_{(u^+,v^+)}, \Y^{\com(\overline{G})}_{(u^+,v^+)})\preccurlyeq^\ominus (u^-,v^-)$. By \cref{p:completion}.(v), $(\X^{\com(\overline{G})}_{(u^+,v^+)}, \Y^{\com(\overline{G})}_{(u^+,v^+)})=(u^-, v^-)=(u^+, v^+)$. {\color{black} Set $(\gamma, \omega):=(u^-,v^-)=(u^+, v^+)$. Since $(u^-, v^-), (u^+, v^+)$ are in $\overline{\mathcal{B}_d(x_n, y_n)\cap S}$ and $\overline{\mathcal{B}_d(x_{n+1}, y_{n+1})\cap S}$ respectively, we have $\|(u^-, v^-)-(x_n, y_n)\|\leq d$ and $\|(u^+, v^+)-(x_{n+1}, y_{n+1})\|\leq d$. Hence,
     \begin{align*}
         \|(\gamma, \omega)-(x_n, y_n)\|&\leq d, \\
         \|(\gamma, \omega)-(x_{n+1}, y_{n+1})\|&\leq d.
     \end{align*} 
     
     \textit{Case 1:} Suppose that $\|(\gamma, \omega)-(x_n, y_n)\|<d$ and $\|(\gamma, \omega)-(x_{n+1}, y_{n+1})\|<d$. Then $(\gamma, \omega)\in \mathcal{B}_d(x_n,y_n)\cap \mathcal{B}_d(x_{n+1}, y_{n+1})$. Moreover, $(\gamma, \omega)\in S$, so $\mathcal{B}_d(x_n, y_n)\cap \mathcal{B}_d(x_{n+1}, y_{n+1})\cap S\neq \varnothing$. This contradicts the assumption $\mathcal{B}_d(x_n, y_n)\cap \mathcal{B}_d(x_{n+1}, y_{n+1})\cap S= \varnothing$. 

     \textit{Case 2:} Suppose that $\|(\gamma, \omega)-(x_n, y_n)\|<d$ and $\|(\gamma, \omega)-(x_{n+1}, y_{n+1})\|=d$. Then $(\gamma, \omega)\in \mathcal{B}_d(x_n, y_n)$ and thus $\mathcal{B}_\delta(\gamma, \omega)\subseteq \mathcal{B}_d(x_n, y_n)$ for some $\delta>0$. Since $(\gamma, \omega)\in \overline{\mathcal{B}_d(x_{n+1}, y_{n+1})\cap S}$, there exists $(\widetilde{\gamma}, \widetilde{\omega})\in \mathcal{B}_d(x_{n+1}, y_{n+1})\cap S$ such that $\|(\gamma, \omega)-(\widetilde{\gamma}, \widetilde{\omega})\|<\delta$, i.e., $(\widetilde{\gamma}, \widetilde{\omega})\in \mathcal{B}_\delta(\gamma, \omega)$. This implies that $(\widetilde{\gamma}, \widetilde{\omega})\in \mathcal{B}_\delta(\gamma, \omega)\cap \mathcal{B}_d(x_{n+1}, y_{n+1})\cap S\subseteq \mathcal{B}_d(x_n, y_n)\cap \mathcal{B}_d(x_{n+1}, y_{n+1})\cap S$. Hence, $\mathcal{B}_d(x_n, y_n)\cap \mathcal{B}_d(x_{n+1}, y_{n+1})\cap S\neq \varnothing$, which contradicts the assumption $\mathcal{B}_d(x_n, y_n)\cap \mathcal{B}_d(x_{n+1}, y_{n+1})\cap S= \varnothing$. 

     \textit{Case 3:} Suppose that $\|(\gamma, \omega)-(x_n, y_n)\|=d$ and $\|(\gamma, \omega)-(x_{n+1}, y_{n+1})\|<d$. Then $(\gamma, \omega)\in \mathcal{B}_d(x_{n+1}, y_{n+1})$ and thus $\mathcal{B}_\delta(\gamma, \omega)\subseteq \mathcal{B}_d(x_{n+1}, y_{n+1})$ for some $\delta>0$.
    Since $(\gamma, \omega)\in \overline{\mathcal{B}_d(x_{n}, y_{n})\cap S}$ , there exists $(\widetilde{\gamma}, \widetilde{\omega})\in \mathcal{B}_d(x_{n}, y_{n})\cap S$ such that $\|(\gamma, \omega)-(\widetilde{\gamma}, \widetilde{\omega})\|<\delta$, i.e., $(\widetilde{\gamma}, \widetilde{\omega})\in \mathcal{B}_\delta(\gamma, \omega)$. This implies that $(\widetilde{\gamma}, \widetilde{\omega})\in \mathcal{B}_\delta(\gamma, \omega)\cap \mathcal{B}_d(x_{n}, y_{n})\cap S\subseteq \mathcal{B}_d(x_n, y_n)\cap \mathcal{B}_d(x_{n+1}, y_{n+1})\cap S$, which contradicts the assumption $\mathcal{B}_d(x_n, y_n)\cap \mathcal{B}_d(x_{n+1}, y_{n+1})\cap S= \varnothing$. 

    \textit{Case 4:} Suppose that $\|(\gamma, \omega)-(x_n, y_n)\|=d$ and $\|(\gamma, \omega)-(x_{n+1}, y_{n+1})\|=d$. 
Let $k\in\{1,\ldots,N\}$, and suppose first 
that $k\leq n$.
Then $(x_k, y_k)\preccurlyeq^\ominus (x_n, y_n)$ 
    and $(x_n, y_n)\preccurlyeq^\ominus (u^-, v^-)=(\gamma, \omega)$ \footnote{Recall \cref{e:revision1}.}. Hence, $(x_k, y_k)\preccurlyeq^\ominus (x_n,y_n)\preccurlyeq^\ominus (\gamma, \omega)$. By \cref{r:cluster}.(ii),  $d=\|(x_n,y_n)-(\gamma, \omega)\|\leq \|(x_k, y_k)-(\gamma, \omega)\|$. Therefore, $(\gamma, \omega)\notin \mathcal{B}_d(x_k, y_k)$ for $\forall k\leq n$. 
Now suppose that $k\geq n+1$. 
Then  $(\gamma, \omega)=(u^+, v^+)\preccurlyeq^\ominus(x_{n+1}, y_{n+1})\preccurlyeq^\ominus(x_k, y_k)$ \footnote{Recall \cref{e:revision2}.}. By \cref{r:cluster}.(ii), $d=\|(x_{n+1}, y_{n+1})-(\gamma, \omega)\|\leq \|(x_k, y_k)-(\gamma, \omega)\|$. Hence, $(\gamma, \omega)\notin \mathcal{B}_d(x_k, y_k)$ for $\forall k\geq n+1$. Summing up, we have shown that
    \begin{align*}
        (\forall k\in \{1, \ldots, N\})\quad (\gamma, \omega)\notin \mathcal{B}_d(x_k, y_k).
    \end{align*}
    However, $(\gamma, \omega)\in S$. This contradicts the fact that $\{\mathcal{B}_d(x_n, y_n)\ |\ (x_n, y_n)\in S, n\in \{1, \ldots, N\}\}$ is a cover of $S$. 
     }
     
     \par Therefore, the assumption $\mathcal{B}_d(x_n,y_n) \cap \mathcal{B}_d(x_{n+1}, y_{n+1})\cap S=\varnothing$ is wrong and thus $\mathcal{B}_d(x_n,y_n)\cap \mathcal{B}_d(x_{n+1},y_{n+1})\cap \com(\overline{G})\cap D \neq \varnothing$ \footnote{Recall that $S\subseteq \com(\overline{G})\cap D$.}. Recall that $\mathcal{B}_d (x_n, y_n)\subseteq D$, $\mathcal{B}_d(x_{n+1}, y_{n+1})\subseteq D$, so $\mathcal{B}_d (x_n, y_n)\subseteq \mathcal{B}(x_n,y_n)$ and $\mathcal{B}_d (x_{n+1}, y_{n+1})\subseteq \mathcal{B}(x_{n+1}, y_{n+1})$ by \cref{d:maxball}. Hence, $\mathcal{B}(x_n,y_n) \cap \mathcal{B}(x_{n+1},y_{n+1})\cap \com(\overline{G}) \cap D\neq \varnothing$. It implies $(x_s,y_s), (x_e,y_e)$ are ball chain connected. Since $(x_s,y_s), (x_e, y_e)$ are arbitrary in $\com(\overline{G})\cap D$, we conclude that any two points in $\com(\overline{G})\cap D$ are ball chain connected. 
     \par \textbf{Step 4:} Combining the conclusion in \textbf{Step 3} and \cref{e:0713b}, we deduce that $\com(\overline{G})\cap D$ is strongly connected. Moreover, $(\com(\overline{G})\cap D, \preccurlyeq^\ominus)$ is a chain, so it is $c$-cyclically monotone by \cref{c:charaofcyc}. Applying \cref{f:0707c}, $\com(\overline{G})\cap D$ is $c$-path bounded. Hence, $\com(\overline{G})\cap D$ is a strongly connected $c$-path bounded extension of $G$. 
     \par \textbf{Step 5:} We know that $\widetilde{G}=\com(\overline{G})\cap D$ is $c$-path bounded and strongly connected, so it is $c$-potentiable, i.e., it has a $c$-antiderivative $f$ by \cref{f:fix} or \cref{f:0707c}. Since $G\subseteq \widetilde{G}$, $f|_G$ is a $c$-antiderivative of $G$. Hence, $G$ is $c$-potentiable.  
\end{proof}

Recall that at the beginning of the section, we have shown that none of the results from \cref{s:intro} to \cref{s:metric} is applicable to prove the $c$-potentiability in \cref{ex:5.1}. Now, with \cref{t:cpathextension}, we can prove the existence of $c$-potential as follows: 

\begin{example}[\cref{ex:5.1} revisited]
\label{ex:5.20}
Let $\widetilde{c}$ and $G$ be the cost and the set in \cref{ex:5.1}. Set $D_1:=\{(x,y)\in \mathbb{R}\times \mathbb{R}\ |\ 0<\frac{1}{x}<y\}$, $D_2:=\{(x,y)\in \mathbb{R}\times \mathbb{R}\ |\ y<\frac{1}{x}<0\}$. Note that $D=D_1\cup D_2$ and $(x,y), (u,v)$ are disconnected for all $(x,y)\in D_1$, $(u,v)\in D_2$. Next, set 
\begin{align*}
\widetilde{c}_1\colon X\times Y \to \left]-\infty, +\infty\right]\colon (x,y)\mapsto\begin{cases}
\widetilde{c}(x,y), &\text{if $(x,y)\in D_1$;}\\
\pinf, &\text{otherwise,}
\end{cases}
\end{align*}
\begin{align*}
\widetilde{c}_2\colon X\times Y \to \left]-\infty, +\infty\right]\colon (x,y)\mapsto\begin{cases}
\widetilde{c}(x,y), &\text{if $(x,y)\in D_2$;}\\
\pinf, &\text{otherwise.}
\end{cases}
\end{align*}
Our goal is to prove $\widetilde{c}$-potentiability of $G$, i.e., $G$ has a $\widetilde{c}$-antiderivative. To do this, we shall first show that $G\cap D_1$ has a $\widetilde{c}_1$-antiderivative using \cref{t:cpathextension}. Let us start by checking that $\widetilde{c}_1$ satisfies items (i)--(iv) of \cref{t:cpathextension}:
\begin{itemize}
\item [(i)] The domain of $\widetilde{c}_1$ is $D_1=\{(x,y)\in \mathbb{R}\times \mathbb{R}\ |\ 0<\frac{1}{x}<y\}$ which is obviously a nonempty open convex proper subset of $\mathbb{R}\times \mathbb{R}$. 

\item [(ii)] Note that $\widetilde{c}_1(x,y)=c_1(x,y)+\chi_{\mathbb{Q}}(y)$ where $\chi_{\mathbb{Q}}$ is the Dirichlet function \footnote{$\chi_\mathbb{Q}(x)=1$ if $x\in \mathbb{Q}$ and $\chi_\mathbb{Q}(x)=0$ if $x\notin \mathbb{Q}$.} and $c_1$ is given in \cref{ex:5.11}.(ii) where we showed that it is $\ominus$-monotone. Hence, $\widetilde{c}_1$ is $\ominus$-monotone by \cref{r:separable} and $\widetilde{c}_1$ satisfies condition (ii) in \cref{t:cpathextension}.

\item [(iii)] For $(x,y)\in D_1$, there is $\mathcal{B}_{\delta}(x,y)$ such that $\sup_{(u,v)\in \mathcal{B}_\delta (x,y)}c_1 (u,v)<+\infty$ because $c_1$ is continuous at $(x,y)$. On the other hand, $\chi_{\mathbb{Q}}$ is bounded, so $\sup_{(u,v)\in \mathcal{B}_\delta (x,y)}\widetilde{c}_1 (u,v)<+\infty$.

\item [(iv)] Note that $c_1(x_n,y_n)\to +\infty$ as $(x_n,y_n)\to (x,y)\in \bd(D_1)$. Since $\chi_{\mathbb{Q}}$ is bounded below, $\widetilde{c}_1(x_n,y_n)\to +\infty$. 
\end{itemize}
Now we show that $G\cap D_1$ satisfies conditions (v) and (vi) in \cref{t:cpathextension}: 
\begin{itemize}
\item [(v)] $G$ is $\widetilde{c}_1$-path bounded (see \hyperref[s:appB]{Appendix}).  

\item [(vi)] Let $\{S_1, S_2\}$ be a partition of $G$. We want to show that there are $(x,y)\in S_1$, $(u,v)\in S_2$ such that there is a path from $(x,y)$ to $(u,v)$ or there is a path from $(u,v)$ to $(x,y)$. For each $\nnn$, take $(x_n,y_n)\in \{\varepsilon_n\}\times [\alpha_n, \beta_n]$. Then $\{(x_n,y_n)\}_{\nnn}\subseteq S_1 \cup S_2$, so at least one of $S_1$ and $S_2$ contains infinitely many elements of $\{(x_n,y_n)\}_\nnn$. Without loss of generality, assume $S_1$ is the set. Fix $(x,y)\in S_2$. 
\begin{enumerate}
\item [Case 1:] Suppose $(x, y)\in [\alpha_m, \beta_m]\times \{\varepsilon_m\}$. Then there is $(x_n,y_n)\in (\{\varepsilon_n\}\times [\alpha_n, \beta_n])\cap S_1$ such that $(x_n, y_n)\prec^\ominus(x,y)$. Hence, $(x, y_n)\in D_1$ and $((x_n,y_n), (x,y))$ is a path from $(x_n, y_n)$ to $(x,y)$. 

\item [Case 2:] Suppose $(x, y)\in \{\varepsilon_m\}\times [\alpha_m, \beta_m]$. Then there is $(x_n,y_n)\in (\{\varepsilon_n\}\times [\alpha_n, \beta_n])\cap S_1$ such that $n>m$. Hence, $(x, y_n)\in D_1$ and $((x_n,y_n), (x,y))$ is a path from $(x_n, y_n)$ to $(x,y)$. 
\end{enumerate}

\end{itemize}
Therefore, all the assumptions in \cref{t:cpathextension} hold and we conclude that $G\cap D_1$ has a strongly connected $\widetilde{c}_1$-path bounded extension ({\color{black} visualized in Figure~\ref{fig:2}}) and a $\widetilde{c}_1$-potential, i.e., 
\begin{align*}
\text{$G\cap D_1$ has a $\widetilde{c}_1$-antiderivative.}
\end{align*}
Moreover, $G\cap D_2=\varnothing$ and $(x,y), (u,v)$ are disconnected for all $(x,y)\in D_1$, $(u,v)\in D_2$. Finally, $G$ is $\widetilde{c}$-potentiable, i.e., 
\begin{align*}
\text{$G$ has a $\widetilde{c}$-antiderivative}
\end{align*}
by \cref{e:2.3}.
\end{example}

\begin{figure}[h]
\centering
\begin{tikzpicture}
  \begin{axis}[
     width=10cm,
     height=10cm,
    axis lines=middle,
    xmin=0, xmax=5,
    ymin=0, ymax=5,
    samples=200,
    domain=0.17:5,
    xlabel={$x$},
    ylabel={$y$},
    xtick=\empty,
    ytick=\empty,
    axis line style={->},
    enlargelimits,
    clip=false
  ]
    \addplot [
      name path=lower,
      domain=0.17:5,
      draw=none
    ] {1/x};

    \path[name path=upper] (axis cs:0.17,5) -- (axis cs:5,5);

    \addplot [
      fill=green!20,
      draw=none
    ] fill between [
      of=upper and lower,
      soft clip={domain=0.17:6}
    ];

    \draw[thick, blue] (axis cs:1,1.3) -- (axis cs:1,2);
    \draw[dashed, thick, red] (axis cs:1,1.3) -- (axis cs:0,1.3);
    \node[anchor=east] at (axis cs:0, 1.3) {\small $\alpha_0$};
    \draw[dashed, thick, red] (axis cs:1,2) -- (axis cs:0,2);
    \node[anchor=east] at (axis cs:0, 2) {\small $\beta_0$};
    \draw[dashed, thick, red] (axis cs:1,1.3) -- (axis cs:1,0);
    \node[anchor=north] at (axis cs:1, 0) {\small $\varepsilon_0$};

\draw[thick, blue] (axis cs:1,2) -- (axis cs:0.5, 2.3);
    \draw[thick, blue] (axis cs:0.5,2.3) -- (axis cs:0.5,3);
    \draw[dashed, thick, red] (axis cs:0.5,2.3) -- (axis cs:0,2.3);
    \node[anchor=east] at (axis cs:0, 2.3) {\small $\alpha_1$};
    \draw[dashed, thick, red] (axis cs:0.5,3) -- (axis cs:0,3);
    \node[anchor=east] at (axis cs:0, 3) {\small $\beta_1$};
    \draw[dashed, thick, red] (axis cs:0.5,2.3) -- (axis cs:0.5,0);
    \node[anchor=north] at (axis cs:0.6, 0) {\small $\varepsilon_1$};

\draw[thick, blue] (axis cs:0.5, 3) -- (axis cs:0.33, 3.3);
    \draw[thick, blue] (axis cs:0.33,3.3) -- (axis cs:0.33,4);
    \draw[dashed, thick, red] (axis cs:0.33,3.3) -- (axis cs:0,3.3);
    \node[anchor=east] at (axis cs:0, 3.3) {\small $\alpha_2$};
    \draw[dashed, thick, red] (axis cs:0.33,4) -- (axis cs:0,4);
    \node[anchor=east] at (axis cs:0, 4) {\small $\beta_2$};
    \draw[dashed, thick, red] (axis cs:0.33,3.3) -- (axis cs:0.33,0);
    \node[anchor=north] at (axis cs:0.33, 0) {\small $\varepsilon_2$};
    
 \draw[thick, blue] (axis cs:1,1.3) -- (axis cs:1.3,1);
    \draw[thick, blue] (axis cs:1.3,1) -- (axis cs:2,1);
    \draw[dashed, thick, red] (axis cs:1.3,1) -- (axis cs:1.3,0);
    \node[anchor=north] at (axis cs:1.3, 0) {\small $\alpha_0$};
    \draw[dashed, thick, red] (axis cs:2,1) -- (axis cs:2,0);
    \node[anchor=north] at (axis cs:2, 0) {\small $\beta_0$};
    \draw[dashed, thick, red] (axis cs:1.3,1) -- (axis cs:0,1);
    \node[anchor=east] at (axis cs:0, 1) {\small $\varepsilon_0$};
    
\draw[thick, blue] (axis cs:2.3,0.5) -- (axis cs:2,1);
    \draw[thick, blue] (axis cs:2.3,0.5) -- (axis cs:3,0.5);
    \draw[dashed, thick, red] (axis cs:2.3,0.5) -- (axis cs:2.3,0);
    \node[anchor=north] at (axis cs:2.3, 0) {\small $\alpha_1$};
    \draw[dashed, thick, red] (axis cs:3,0.5) -- (axis cs:3,0);
    \node[anchor=north] at (axis cs:3, 0) {\small $\beta_1$};
    \draw[dashed, thick, red] (axis cs:2.3,0.5) -- (axis cs:0,0.5);
    \node[anchor=east] at (axis cs:0, 0.5) {\small $\varepsilon_1$};

\draw[thick, blue] (axis cs:3.3,0.33) -- (axis cs:3, 0.5);
    \draw[thick, blue] (axis cs:3.3,0.33) -- (axis cs:4,0.33);
    \draw[dashed, thick, red] (axis cs:3.3,0.33) -- (axis cs:3.3,0);
    \node[anchor=north] at (axis cs:3.3, 0) {\small $\alpha_2$};
    \draw[dashed, thick, red] (axis cs:4,0.33) -- (axis cs:4,0);
    \node[anchor=north] at (axis cs:4, 0) {\small $\beta_2$};
    \draw[dashed, thick, red] (axis cs:3.3,0.33) -- (axis cs:0,0.33);
    \node[anchor=east] at (axis cs:0, 0.33) {\small $\varepsilon_2$};

      \fill[black] (axis cs:4.5,0.33) circle[radius=1pt];
     \fill[black] (axis cs:4.6, 0.33) circle[radius=1pt];
     \fill[black] (axis cs:4.7, 0.33) circle[radius=1pt];
     
    \fill[black] (axis cs:0.33,4.5) circle[radius=1pt];
     \fill[black] (axis cs:0.33, 4.6) circle[radius=1pt];
     \fill[black] (axis cs:0.33, 4.7) circle[radius=1pt];

  \end{axis}
\end{tikzpicture}
\caption{{\scriptsize 
{\color{black} The $\widetilde{c}_1$-path bounded extension 
$\widetilde{G}$ of $G$ from Figure~\ref{fig:1} and \cref{ex:5.1}. See 
\cref{ex:5.20} for details.}}}
\label{fig:2}
\end{figure}

{\color{black}
\begin{remark}
    To avoid repetition, we omit the analogous discussion 
    of  $c$-path bounded extensions for $\oplus$-monotone costs.
\end{remark}
}

{\color{black}
\section{Future work}
We conclude by highlighting several directions for future research.

First, uniqueness of $c$-potentials deserves attention and additional study. Although uniqueness does not hold in general, as demonstrated by \cref{ex:not unique}, it is natural to seek conditions on the cost $c$ and the set $G$ that guarantee uniqueness. A classical result in convex analysis\footnote{See, for example, \cite[Theorem~22.18 and Proposition~22.19]{BC2017}.} states that if $X=Y=\HH$ is a Hilbert space, and if $G$ is maximally cyclically monotone with respect to the cost $-\langle \cdot,\cdot\rangle$, then the potential of $G$ is unique up to a constant. This raises the natural question of whether this result can be extended to more general spaces and costs.

A second natural question is how to derive explicit formulas for $c$-potentials. We have constructed $c$-potentials using maximal inner variations, which depend only on the cost. Given an explicit formula for the cost, we therefore ask whether the corresponding maximal inner variation can be formulated explicitly as a function on $G\times G$. Such a formula will facilitate computations and visualizations of $c$-potentials.

Third, it is worthwhile to explore $c$-path bounded extensions in higher dimensions. 
The technique in \cref{s:extension} applies to monotone costs on $\mathbb{R}\times\mathbb{R}$, where it relies crucially on the order structure of the real line. Since no analogous natural total order exists in higher dimensions, it remains unclear how the $c$-path bounded extension technique can be adapted to $\mathbb{R}^d\times\mathbb{R}^d$ for $d\geq 2$.
}

\section*{Acknowledgments}
\small
The authors are grateful to Kasia Wyczesany for helpful email discussions
regarding non-traditional costs {\color{black} and to two anonymous referees for constructive suggestions which greatly improved this work over the original version \cite{BBGv1}.} 
The research of SB was partially supported by a collaboration grant for mathematicians of the Simons Foundation. The research of HHB was partially supported by a Discovery Grant
of the Natural Sciences and Engineering Research Council of
Canada.

\appendix
\section*{Appendix}
\phantomsection
\label{s:appB}
We shall show that $G$ in \cref{ex:5.20} is $\widetilde{c}_1$-path bounded. Take $(x_s,y_s), (x_e,y_e)\in G$ and suppose that $\big((x_k, y_k)\big)_{k=1}^{N+1}\in P^G_{(x_s,y_s)\to (x_e,y_e)} \neq \varnothing$ \footnote{Otherwise it is trivial that $\sum_{k=1}^N \widetilde{c}_1(x_k,y_k)-\widetilde{c}_1(x_{k+1}, y_k)$ has an upper bound.}. We need to find an upper bound for $\sum_{k=1}^N (\widetilde{c}_1(x_k,y_k)-\widetilde{c}_1(x_{k+1}, y_k))$, which only depends on $(x_s,y_s), (x_e,y_e)$. 

\begin{enumerate}
\item [Case 1:] Suppose $(x_s,y_s), (x_e,y_e)\in G_2$. Then for each $k\in \{1, \ldots, N\}$, there are $m,n\in \mathbb{N}$ with $m\leq n$ such that $(x_k, y_k)\in \{\varepsilon_n\}\times [\alpha_n, \beta_n]$ and $(x_{k+1}, y_{k+1})\in \{\varepsilon_m\}\times [\alpha_m, \beta_m]$. Let $n(k)$ be the number that $(x_k, y_k)\in \{\varepsilon_{n(k)}\}\times [\alpha_{n(k)}, \beta_{n(k)}]$ and we know $n(k+1)\leq n(k)$ for each $k\in \{1, \ldots, N\}$. Note that for $k\in \{1, \ldots, N\}$, $(\varepsilon_{n(k)}, \alpha_{n(k)}), (\varepsilon_{n(k+1)}, y_k)\in D_1$ and $(\varepsilon_{n(k)}, \alpha_{n(k)})\preccurlyeq^\oplus (\varepsilon_{n(k+1)}, y_k)$, so $\ominus$-monotonicity of $\widetilde{c}_1$ gives 
\begin{align*}
\widetilde{c}_1(x_k,y_k)-\widetilde{c}_1(x_{k+1}, y_k)=\widetilde{c}_1(\varepsilon_{n(k)}, y_k)-\widetilde{c}_1(\varepsilon_{n(k+1)}, y_k)\leq \widetilde{c}_1(\varepsilon_{n(k)}, \alpha_{n(k)})-\widetilde{c}_1(\varepsilon_{n(k+1)}, \alpha_{n(k)}).
\end{align*}
Hence, 
\begin{equation}
\label{e:5.20e1}
\sum_{k=1}^N \big(\widetilde{c}_1(x_k, y_k)-\widetilde{c}_1(x_{k+1}, y_k)\big)\leq \sum _{k=1}^N\big(\widetilde{c}_1(\varepsilon_{n(k)}, \alpha_{n(k)})-\widetilde{c}_1(\varepsilon_{n(k+1)}, \alpha_{n(k)})\big).
\end{equation}
Set $I:=\{k\in \{1, \ldots, N\}\ |\ n(k)-n(k+1)\geq 1\}$. Note that $|I|\leq n(1)-n(N+1)$ and $\widetilde{c}_1(\varepsilon_{n(k)}, \alpha_{n(k)})-\widetilde{c}_1(\varepsilon_{n(k+1)}, \alpha_{n(k)})=0$ if $k\notin I$. Therefore, \cref{e:5.20e1} becomes 
\begin{equation}
\label{e:5.20e2}
\begin{split}
\sum_{k=1}^N \big(\widetilde{c}_1(x_k, y_k)-\widetilde{c}_1(x_{k+1}, y_k)\big)&\leq \sum _{k=1}^N\big(\widetilde{c}_1(\varepsilon_{n(k)}, \alpha_{n(k)})-\widetilde{c}_1(\varepsilon_{n(k+1)}, \alpha_{n(k)})\big)\\
&=\sum_{k\in I}\big(\widetilde{c}_1(\varepsilon_{n(k)}, \alpha_{n(k)})-\widetilde{c}_1(\varepsilon_{n(k+1)}, \alpha_{n(k)})\big).
\end{split}
\end{equation}
Moreover, fix $k\in I$ and for $i\in \{0, \ldots, n(k)-n(k+1)-1\}$, we have $(\varepsilon_{n(k)-i}, \alpha_{n(k)-i}), (\varepsilon_{n(k)-i-1}, \alpha_{n(k)})\in D_1$ and $(\varepsilon_{n(k)-i}, \alpha_{n(k)-i})\preccurlyeq^\oplus (\varepsilon_{n(k)-i-1}, \alpha_{n(k)})$, so $\ominus$-monotonicity of $\widetilde{c}_1$ yields
\begin{equation}
\label{e:5.20e3}
\begin{split}
\widetilde{c}_1(\varepsilon_{n(k)}, \alpha_{n(k)})-\widetilde{c}_1(\varepsilon_{n(k+1)}, \alpha_{n(k)})&=\sum_{i=0}^{n(k)-n(k+1)-1}\big(\widetilde{c}_1(\varepsilon_{n(k)-i}, \alpha_{n(k)})-\widetilde{c}_1(\varepsilon_{n(k)-i-1}, \alpha_{n(k)})\big)\\
&\leq \sum_{i=0}^{n(k)-n(k+1)-1}\big(\widetilde{c}_1(\varepsilon_{n(k)-i}, \alpha_{n(k)-i})-\widetilde{c}_1(\varepsilon_{n(k)-i-1}, \alpha_{n(k)-i})\big).
\end{split}
\end{equation}
Combining \cref{e:5.20e2} and \cref{e:5.20e3}, we get 
\begin{align*}
\sum_{k=1}^N \big(\widetilde{c}_1(x_k, y_k)-\widetilde{c}_1(x_{k+1}, y_k)\big)&\leq \sum_{k\in I} \sum_{i=0}^{n(k)-n(k+1)-1}\big(\widetilde{c}_1(\varepsilon_{n(k)-i}, \alpha_{n(k)-i})-\widetilde{c}_1(\varepsilon_{n(k)-i-1}, \alpha_{n(k)-i})\big)\\
&=\sum_{j=0}^{n(1)-n(N+1)-1}\big(\widetilde{c}_1(\varepsilon_{n(1)-j}, \alpha_{n(1)-j})-\widetilde{c}_1(\varepsilon_{n(1)-j-1}, \alpha_{n(1)-j})\big)
<+\infty.
\end{align*}
Note that the right-hand side only depends on $n(1)$ and $n(N+1)$ which are determined by $(x_s, y_s)$ and $(x_e, y_e)$. Therefore, we have $\widetilde{c}_1$-path boundedness. 

\item [Case 2:] Suppose that $(x_s,y_s), (x_e,y_e)\in G_1$. Then for each $k\in \{1, \ldots, N\}$, there are $m,n\in \mathbb{N}$ with $n\leq m$ such that $(x_k,y_k)\in [\alpha_n, \beta_n]\times \{\varepsilon_n\}$ and $(x_{k+1}, y_{k+1})\in [\alpha_m, \beta_m]\times \{\varepsilon_m\}$. Let $n(k)$ be the number that $(x_k,y_k)\in [\alpha_{n(k)}, \beta_{n(k)}]\times \{\varepsilon_{n(k)}\}$. Set $\{k_i\}_{i=1}^M$ \footnote{{\color{black} We assume that $M\geq 1$, i.e., $\{k_i\}_{i=1}^M\neq \varnothing$. In fact, the upper bound we obtain under the assumption $M\geq 1$ also works for $M=0$. Indeed, if $M=0$, then $n(1)=n(2)=\cdots=n(N+1)$ and $\widetilde{c}_1(\alpha_{n(1)}, \varepsilon_{n(1)}) -\widetilde{c}_1(\beta_{n(1)}, \varepsilon_{n(1)})+\sum_{n=n(1)}^{n(N+1)-1} \big(\widetilde{c}_1(\alpha_n, \varepsilon_n)-\widetilde{c}_1(\alpha_{n+1}, \varepsilon_n)\big)
+2\sum_{n=n(1)}^{n(N+1)} \big(\widetilde{c}_1(\alpha_{n}, \varepsilon_{n})-\widetilde{c}_1(\beta_{n}, \varepsilon_{n})\big)=3(\widetilde{c}_1(\alpha_{n(1)}, \varepsilon_{n(1)})-\widetilde{c}_1(\beta_{n(1)}, \varepsilon_{n(1)}))\geq \widetilde{c}_1(\alpha_{n(1)}, \varepsilon_{n(1)})-\widetilde{c}_1(\beta_{n(1)}, \varepsilon_{n(1)})\geq \widetilde{c}_1(x_1, \varepsilon_{n(1)})-\widetilde{c}_1(x_{N+1}, \varepsilon_{n(1)})=\sum_{k=1}^N(\widetilde{c}_1(x_k, \varepsilon_{n(1)})-\widetilde{c}_1(x_{k+1}, \varepsilon_{n(1)}))=\sum_{k=1}^N(\widetilde{c}_1(x_k, y_k)-\widetilde{c}_1(x_{k+1}, y_{k}))$.} } be the subset of the indices $\{1, \ldots, N+1\}$ such that $n(k_{i+1})-n(k_{i})\geq 1$ and {\color{black} set $k_{M+1}:=N+1$}, i.e., for each $i\in \{1, \ldots, M\}$, $n(k_{i}+1)=n(k)$ for all $k\in \{k_{i}+1, \ldots, k_{i+1}\}$. Hence, 
\begin{align*}
&\quad \ \ \sum_{k=1}^N \big(\widetilde{c}_1(x_k, y_k)-\widetilde{c}_1(x_{k+1}, y_k)\big)\\
&=\sum_{k=1}^N \big(\widetilde{c}_1(x_k, \varepsilon_{n(k)})-\widetilde{c}_1(x_{k+1}, \varepsilon_{n(k)})\big)\\
&=\sum_{k=1}^{k_1-1} \big(\widetilde{c}_1(x_k, \varepsilon_{n(k)}) -\widetilde{c}_1(x_{k+1}, \varepsilon_{n(k)})\big)+ \sum_{i=1}^M \sum_{k=k_i}^{k_{i+1}-1}\big(\widetilde{c}_1(x_k, \varepsilon_{n(k)})-\widetilde{c}_1(x_{k+1}, \varepsilon_{n(k)})\big)\\
&=\sum_{k=1}^{k_1-1} \big(\widetilde{c}_1(x_k, \varepsilon_{n(k)}) -\widetilde{c}_1(x_{k+1}, \varepsilon_{n(k)})\big)\\
&\quad +\sum_{i=1}^M \big(\widetilde{c}_1(x_{k_i}, \varepsilon_{n(k_i)})-\widetilde{c}_1(x_{k_i+1}, \varepsilon_{n(k_i)})+\sum_{k=k_i+1}^{k_{i+1}-1}\widetilde{c}_1(x_k, \varepsilon_{n(k)})-\widetilde{c}_1(x_{k+1}, \varepsilon_{n(k)})\big)\\
&=\sum_{k=1}^{k_1-1} \big(\widetilde{c}_1(x_k, \varepsilon_{n(k)}) -\widetilde{c}_1(x_{k+1}, \varepsilon_{n(k)})\big)+\sum_{i=1}^M \big(\widetilde{c}_1(x_{k_i}, \varepsilon_{n(k_i)})-\widetilde{c}_1(x_{k_i+1}, \varepsilon_{n(k_i)})\big)\\
&\quad +\sum_{i=1}^M\sum_{k=k_i+1}^{k_{i+1}-1}\big(\widetilde{c}_1(x_k, \varepsilon_{n(k)})-\widetilde{c}_1(x_{k+1}, \varepsilon_{n(k)})\big)\\
&= \sum_{k=1}^{k_1-1} \big(\widetilde{c}_1(x_k, \varepsilon_{n(k)}) -\widetilde{c}_1(x_{k+1}, \varepsilon_{n(k)})\big)+\sum_{i=1}^M \big(\widetilde{c}_1(x_{k_i}, \varepsilon_{n(k_i)})-\widetilde{c}_1(x_{k_i+1}, \varepsilon_{n(k_i)})\big)\\
&\quad +\sum_{i=1}^M\sum_{k=k_i+1}^{k_{i+1}-1}\big(\widetilde{c}_1(x_k, \varepsilon_{n(k_i+1)})-\widetilde{c}_1(x_{k+1}, \varepsilon_{n(k_i+1)})\big)\\
&=\sum_{k=1}^{k_1-1} \big(\widetilde{c}_1(x_k, \varepsilon_{n(k)}) -\widetilde{c}_1(x_{k+1}, \varepsilon_{n(k)})\big)+\sum_{i=1}^M \big(\widetilde{c}_1(x_{k_i}, \varepsilon_{n(k_i)})-\widetilde{c}_1(x_{k_i+1}, \varepsilon_{n(k_i)})\big)\\
&\quad +\sum_{i=1}^M\big(\widetilde{c}_1(x_{k_i+1}, \varepsilon_{n(k_i+1)})-\widetilde{c}_1(x_{k_{i+1}}, \varepsilon_{n(k_i+1)})\big)\\
&\leq \sum_{k=1}^{k_1-1} \big(\widetilde{c}_1(x_k, \varepsilon_{n(k)}) -\widetilde{c}_1(x_{k+1}, \varepsilon_{n(k)})\big)+\sum_{i=1}^M \big(\widetilde{c}_1(x_{k_i}, \varepsilon_{n(k_i)})-\widetilde{c}_1(x_{k_i+1}, \varepsilon_{n(k_i)})\big)\\
&\quad +\sum_{i=1}^M\big(\widetilde{c}_1(\alpha_{n(k_i+1)}, \varepsilon_{n(k_i+1)})-\widetilde{c}_1(\beta_{n(k_i+1)}, \varepsilon_{n(k_i+1)})\big)\\
&\leq \sum_{k=1}^{k_1-1} \big(\widetilde{c}_1(x_k, \varepsilon_{n(k)}) -\widetilde{c}_1(x_{k+1}, \varepsilon_{n(k)})\big)+\sum_{i=1}^M \big(\widetilde{c}_1(\alpha_{n(k_i)}, \varepsilon_{n(k_i)})-\widetilde{c}_1(\beta_{n(k_i+1)}, \varepsilon_{n(k_i)})\big)\\
&\quad +\sum_{i=1}^M\big(\widetilde{c}_1(\alpha_{n(k_i+1)}, \varepsilon_{n(k_i+1)})-\widetilde{c}_1(\beta_{n(k_i+1)}, \varepsilon_{n(k_i+1)})\big)\\
&=\sum_{k=1}^{k_1-1} \big(\widetilde{c}_1(x_k, \varepsilon_{n(k)}) -\widetilde{c}_1(x_{k+1}, \varepsilon_{n(k)})\big)\\
&\quad +\sum_{i=1}^M \big(\widetilde{c}_1(\alpha_{n(k_i)}, \varepsilon_{n(k_i)})+\sum_{n=n(k_i)+1}^{n(k_i+1)}\big(-\widetilde{c}_1(\alpha_{n}, \varepsilon_{n(k_i)})+\widetilde{c}_1(\alpha_{n}, \varepsilon_{n(k_i)})\big)-\widetilde{c}_1(\beta_{n(k_i+1)}, \varepsilon_{n(k_i)})\big)\\
&\quad +\sum_{i=1}^M\big(\widetilde{c}_1(\alpha_{n(k_i+1)}, \varepsilon_{n(k_i+1)})-\widetilde{c}_1(\beta_{n(k_i+1)}, \varepsilon_{n(k_i+1)})\big).
\end{align*}
Note that $(\alpha_{n}, \varepsilon_n)\prec^\oplus (\alpha_{n+1}, \varepsilon_{n(k_i)})$ for $n\in \{n(k_i)+1, \ldots, n(k_{i+1})\}$, and $(\alpha_{n(k_i+1)}, \varepsilon_{n(k_i+1)})\prec^\oplus (\beta_{n(k_i+1)}, \varepsilon_{n(k_i)})$ so by $\ominus$-monotonicity of $\widetilde{c}_1$, we have 
\begin{equation}
\label{e:appCe1}
\begin{split}
\widetilde{c}_1(\alpha_n, \varepsilon_{n(k_i)})- \widetilde{c}_1(\alpha_{n+1}, \varepsilon_{n(k_i)})&\leq \widetilde{c}_1(\alpha_n, \varepsilon_n)-\widetilde{c}_1(\alpha_{n+1}, \varepsilon_n)\\
\widetilde{c}_1(\alpha_{n(k_i+1)}, \varepsilon_{n(k_i)})- \widetilde{c}_1(\beta_{n(k_i+1)}, \varepsilon_{n(k_i)})&\leq \widetilde{c}_1(\alpha_{n(k_i+1)}, \varepsilon_{n(k_i+1)})-\widetilde{c}_1(\beta_{n(k_i+1)}, \varepsilon_{n(k_i+1)}).
\end{split}
\end{equation}
Hence, 
\begin{align*}
&\sum_{k=1}^{k_1-1} \big(\widetilde{c}_1(x_k, \varepsilon_{n(k)}) -\widetilde{c}_1(x_{k+1}, \varepsilon_{n(k)})\big)\\
&\quad +\sum_{i=1}^M \big(\widetilde{c}_1(\alpha_{n(k_i)}, \varepsilon_{n(k_i)})+\sum_{n=n(k_i)+1}^{n(k_i+1)}\big(-\widetilde{c}_1(\alpha_{n}, \varepsilon_{n(k_i)})+\widetilde{c}_1(\alpha_{n}, \varepsilon_{n(k_i)})\big)-\widetilde{c}_1(\beta_{n(k_i+1)}, \varepsilon_{n(k_i)})\big)\\
&\quad +\sum_{i=1}^M\big(\widetilde{c}_1(\alpha_{n(k_i+1)}, \varepsilon_{n(k_i+1)})-\widetilde{c}_1(\beta_{n(k_i+1)}, \varepsilon_{n(k_i+1)})\big)\\
&\leq \sum_{k=1}^{k_1-1} \big(\widetilde{c}_1(x_k, \varepsilon_{n(k)}) -\widetilde{c}_1(x_{k+1}, \varepsilon_{n(k)})\big)\\
&\quad +\sum_{i=1}^M\big(\sum_{n=n(k_i)}^{n(k_i+1)-1}\big(\widetilde{c}_1(\alpha_{n}, \varepsilon_{n})-\widetilde{c}_1(\alpha_{n+1}, \varepsilon_{n})\big)+\widetilde{c}_1(\alpha_{n(k_i+1)}, \varepsilon_{n(k_i+1)})-\widetilde{c}_1(\beta_{n(k_i+1)}, \varepsilon_{n(k_i+1)})\big)\\
&\quad +\sum_{i=1}^M\big(\widetilde{c}_1(\alpha_{n(k_i+1)}, \varepsilon_{n(k_i+1)})-\widetilde{c}_1(\beta_{n(k_i+1)}, \varepsilon_{n(k_i+1)})\big) \tag{by \cref{e:appCe1}}\\
&= \widetilde{c}_1(x_1, \varepsilon_{n(1)}) -\widetilde{c}_1(x_{k_1}, \varepsilon_{n(1)})+\sum_{i=1}^M\sum_{n=n(k_i)}^{n(k_i+1)-1}\big(\widetilde{c}_1(\alpha_{n}, \varepsilon_{n})-\widetilde{c}_1(\alpha_{n+1}, \varepsilon_{n})\big)\\
&\quad +2\sum_{i=1}^M \big(\widetilde{c}_1(\alpha_{n(k_i+1)}, \varepsilon_{n(k_i+1)})-\widetilde{c}_1(\beta_{n(k_i+1)}, \varepsilon_{n(k_i+1)})\big)\\
&\leq \widetilde{c}_1(\alpha_{n(1)}, \varepsilon_{n(1)}) -\widetilde{c}_1(\beta_{n(1)}, \varepsilon_{n(1)})+\sum_{n=n(1)}^{n(N+1)-1} \big(\widetilde{c}_1(\alpha_n, \varepsilon_n)-\widetilde{c}_1(\alpha_{n+1}, \varepsilon_n)\big)\\
&\quad +2\sum_{n=n(1)}^{n(N+1)} \big(\widetilde{c}_1(\alpha_{n}, \varepsilon_{n})-\widetilde{c}_1(\beta_{n}, \varepsilon_{n})\big). 
\end{align*}
Note that the last inequality only depends on $n(1)$ and $n(N+1)$ which are determined by $(x_s, y_s), (x_e,y_e)$, so we have found a uniform upper bound and thus $\widetilde{c}_1$-path bounded. 

\item [Case 3:] Suppose that $(x_s, y_s)\in \big(\{\varepsilon_n\}\times [\alpha_n, \beta_n]\big)\subseteq G_2$ and $(x_e,y_e)\in \big([\alpha_m, \beta_m]\times \{\varepsilon_m\}\big)\subseteq G_1$. Then there is $j\in \{1, \ldots, N\}$ such that $(x_k, y_k)\in G_2$ for $k\in \{1, \ldots, j\}$ and $(x_k, y_k)\in G_1$ for $k\in \{j+1, \ldots, N+1\}$. Therefore, 
\begin{equation}
\label{e:5.20e4}
\begin{split}
&\quad \sum_{k=1}^{N} \big(\widetilde{c}_1(x_k, y_k)-\widetilde{c}_1(x_{k+1}, y_k)\big)\\
&=\sum_{k=1}^{j-1} \big(\widetilde{c}_1(x_k, y_k)-\widetilde{c}_1(x_{k+1}, y_k)\big)+\big(\widetilde{c}_1(x_j, y_j)-\widetilde{c}_1(x_{j+1}, y_j)\big)+\sum_{k=j+1}^N \big(\widetilde{c}_1(x_{k}, y_k)-\widetilde{c}_1(x_{k+1}, y_k)\big).
\end{split}
\end{equation}
Let $n(k)$ be the number such that $(x_k, y_k)\in \{\varepsilon_{n(k)}\}\times [\alpha_{n(k)}, \beta_{n(k)}]$ or $(x_k, y_k)\in [\alpha_{n(k)}, \beta_{n(k)}]\times \{\varepsilon_{n(k)}\}$. By Case 1, we have 
\begin{equation}
\label{e:5.20e5}
\begin{split}
\sum_{k=1}^{j-1}\big(\widetilde{c}_1(x_k, y_k)-\widetilde{c}_1(x_{k+1}, y_k)\big)&\leq\sum_{i=0}^{n(1)-n(j)-1}\big(\widetilde{c}_1(\varepsilon_{n(1)-i}, \alpha_{n(1)-i})-\widetilde{c}_1(\varepsilon_{n(1)-i-1}, \alpha_{n(1)-i})\big)\\
&\leq \sum_{i=0}^{n(1)-1}\big(\widetilde{c}_1(\varepsilon_{n(1)-i}, \alpha_{n(1)-i})-\widetilde{c}_1(\varepsilon_{n(1)-i-1}, \alpha_{n(1)-i})\big).
\end{split}
\end{equation}
By Case 2, we have
\begin{equation}
\label{e:5.20e6}
\begin{split}
&\quad \sum_{k=j+1}^N \big(\widetilde{c}_1(x_k,y_k)-\widetilde{c}_1(x_{k+1}, y_k)\big)\\
&\leq \widetilde{c}_1(\alpha_{n(j+1)}, \varepsilon_{n(j+1)}) -\widetilde{c}_1(\beta_{n(j+1)}, \varepsilon_{n(j+1)})+\sum_{n=n(j+1)}^{n(N+1)-1} \big(\widetilde{c}_1(\alpha_n, \varepsilon_n)-\widetilde{c}_1(\alpha_{n+1}, \varepsilon_n)\big)\\
&\quad +2\sum_{n=n(j+1)}^{n(N+1)} \widetilde{c}_1(\alpha_{n}, \varepsilon_{n})-\widetilde{c}_1(\beta_{n}, \varepsilon_{n}).
\end{split}
\end{equation}
Finally, $(\varepsilon_{n(j)}, \alpha_{n(j)})\preccurlyeq^\oplus (x_{j+1}, y_j)$, so $\ominus$-monotonicity of $\widetilde{c}_1$ yields
\begin{equation}
\label{e:Appendix1}
\begin{split}
 &\quad \widetilde{c}_1(x_j, y_j)-\widetilde{c}_1(x_{j+1}, y_j)\\
 &=\widetilde{c}_1(\varepsilon_{n(j)}, y_j)-\widetilde{c}_1(x_{j+1}, y_j)\\
 &\leq \widetilde{c}_1(\varepsilon_{n(j)}, \alpha_{n(j)})-\widetilde{c}_1(x_{j+1}, \alpha_{n(j)})\\
 &=\widetilde{c}_1(\varepsilon_{n(j)}, \alpha_{n(j)})+\sum_{i=1}^{n(j)}\big(-\widetilde{c}_1(\varepsilon_{n(j)-i}, \alpha_{n(j)})+\widetilde{c}_1(\varepsilon_{n(j)-i}, \alpha_{n(j)})\big)-\widetilde{c}_1(x_{j+1}, \alpha_{n(j)})\\
 &=\sum_{i=0}^{n(j)-1}\big(\widetilde{c}_1(\varepsilon_{n(j)-i}, \alpha_{n(j)})-\widetilde{c}_1(\varepsilon_{n(j)-i-1}, \alpha_{n(j)})\big)+\widetilde{c}_1(\varepsilon_0, \alpha_{n(j)})-\widetilde{c}_1(x_{j+1}, \alpha_{n(j)})\\
 &\leq \sum_{i=0}^{n(j)-1}\big(\widetilde{c}_1(\varepsilon_{n(j)-i}, \alpha_{n(j)})-\widetilde{c}_1(\varepsilon_{n(j)-i-1}, \alpha_{n(j)})\big)+\widetilde{c}_1(\varepsilon_0, \alpha_{n(j)})-\widetilde{c}_1(\beta_{n(j+1)}, \alpha_{n(j)}).
 \end{split}
\end{equation}
Note that $(\varepsilon_{n(j)-i}, \alpha_{n(j)-i})\preccurlyeq^\oplus (\varepsilon_{n(j)-i-1}, \alpha_{n(j)})$ for $i\in \{0, \ldots, n(j)-1\}$ \footnote{{\color{black}If $n(j)=0$, then $\sum_{i=0}^{n(j)-1}(\widetilde{c}_1(\varepsilon_{n(j)-i}, \alpha_{n(j)})-\widetilde{c}_1(\varepsilon_{n(j)-i-1}, \alpha_{n(j)}))=0$ and we only need to consider $\widetilde{c}_1(\varepsilon_0, \alpha_{n(j)})-\widetilde{c}_1(\beta_{n(j+1)},\alpha_{n(j)})$.}}, and $(\varepsilon_0, \alpha_0)\preccurlyeq^\oplus(\beta_{n(j+1)}, \alpha_{n(j)})$. By $\ominus$-monotonicity again 
\begin{equation}
\label{e:Appendix2}
\begin{split}
&\quad \sum_{i=0}^{n(j)-1}\big(\widetilde{c}_1(\varepsilon_{n(j)-i}, \alpha_{n(j)})-\widetilde{c}_1(\varepsilon_{n(j)-i-1}, \alpha_{n(j)})\big)+\widetilde{c}_1(\varepsilon_0, \alpha_{n(j)})-\widetilde{c}_1(\beta_{n(j+1)}, \alpha_{n(j)})\\
&\leq \sum_{i=0}^{n(j)-1}\big(\widetilde{c}_1(\varepsilon_{n(j)-i}, \alpha_{n(j)-i})-\widetilde{c}_1(\varepsilon_{n(j)-i-1}, \alpha_{n(j)-i})\big)+\widetilde{c}_1(\varepsilon_0, \alpha_{0})-\widetilde{c}_1(\beta_{n(j+1)}, \alpha_{0})\\
&=\sum_{i=0}^{n(j)-1}\big(\widetilde{c}_1(\varepsilon_{n(j)-i}, \alpha_{n(j)-i})-\widetilde{c}_1(\varepsilon_{n(j)-i-1}, \alpha_{n(j)-i})\big)\\
&\quad +\widetilde{c}_1(\varepsilon_0, \alpha_{0})+\sum_{n=0}^{n(j+1)}\big(-\widetilde{c}_1(\alpha_n, \alpha_{0})+\widetilde{c}_1(\alpha_n, \alpha_0)\big)-\widetilde{c}_1(\beta_{n(j+1)}, \alpha_0)\\
&=\sum_{i=0}^{n(j)-1}\big(\widetilde{c}_1(\varepsilon_{n(j)-i}, \alpha_{n(j)-i})-\widetilde{c}_1(\varepsilon_{n(j)-i-1}, \alpha_{n(j)-i})\big)\\
&\quad +\widetilde{c}_1(\varepsilon_0, \alpha_0)-\widetilde{c}_1(\alpha_0, \alpha_0)+\sum_{n=0}^{n(j+1)-1}\big(\widetilde{c}_1(\alpha_n,\alpha_0)-\widetilde{c}_1(\alpha_{n+1}, \alpha_0)\big)+\widetilde{c}_1(\alpha_{n(j+1)}, \alpha_0)-\widetilde{c}_1(\beta_{n(j+1)}, \alpha_0).
\end{split}
\end{equation}
Note that $(\alpha_n, \varepsilon_n)\preccurlyeq^\oplus (\alpha_{n+1}, \alpha_0)$ for $n\in \{0, \ldots, n(j+1)-1\}$ \footnote{{\color{black}If $n(j+1)=0$, then $\sum_{n=0}^{n(j+1)-1}\big(\widetilde{c}_1(\alpha_n,\alpha_0)-\widetilde{c}_1(\alpha_{n+1}, \alpha_0)\big)=0$.}}, and $(\alpha_{n(j+1)}, \varepsilon_{n(j+1)})\preccurlyeq^\oplus (\beta_{n(j+1)}, \alpha_0)$. $\ominus$-monotonicity implies 
\begin{equation}
\label{e:5.20e7}
\begin{split}
&\sum_{i=0}^{n(j)-1}\big(\widetilde{c}_1(\varepsilon_{n(j)-i}, \alpha_{n(j)-i})-\widetilde{c}_1(\varepsilon_{n(j)-i-1}, \alpha_{n(j)-i})\big)\\
&\quad +\widetilde{c}_1(\varepsilon_0, \alpha_0)-\widetilde{c}_1(\alpha_0, \alpha_0)+\sum_{n=0}^{n(j+1)-1}\big(\widetilde{c}_1(\alpha_n,\alpha_0)-\widetilde{c}_1(\alpha_{n+1}, \alpha_0)\big)+\widetilde{c}_1(\alpha_{n(j+1)}, \alpha_0)-\widetilde{c}_1(\beta_{n(j+1)}, \alpha_0)\\
&\leq \sum_{i=0}^{n(j)-1}\big(\widetilde{c}_1(\varepsilon_{n(j)-i}, \alpha_{n(j)-i})-\widetilde{c}_1(\varepsilon_{n(j)-i-1}, \alpha_{n(j)-i})\big)\\
&\quad +\widetilde{c}_1(\varepsilon_0, \alpha_0)-\widetilde{c}_1(\alpha_0, \alpha_0)+\sum_{n=0}^{n(j+1)-1}\big(\widetilde{c}_1(\alpha_n,\varepsilon_n)-\widetilde{c}_1(\alpha_{n+1}, \varepsilon_n)\big)+\widetilde{c}_1(\alpha_{n(j+1)}, \varepsilon_{n(j+1)})-\widetilde{c}_1(\beta_{n(j+1)}, \varepsilon_{n(j+1)})\\
&\leq \sum_{i=0}^{n(1)-n(j)-1}\big(\widetilde{c}_1(\varepsilon_{n(1)-i}, \alpha_{n(1)-i})-\widetilde{c}_1(\varepsilon_{n(1)-i-1}, \alpha_{n(1)-i})\big)+\sum_{i=0}^{n(j)-1}\big(\widetilde{c}_1(\varepsilon_{n(j)-i}, \alpha_{n(j)-i})-\widetilde{c}_1(\varepsilon_{n(j)-i-1}, \alpha_{n(j)-i})\big)\\
&\quad +\widetilde{c}_1(\varepsilon_0, \alpha_0)-\widetilde{c}_1(\alpha_0, \alpha_0)+\sum_{n=0}^{n(j+1)-1}\big(\widetilde{c}_1(\alpha_n,\varepsilon_n)-\widetilde{c}_1(\alpha_{n+1}, \varepsilon_n)\big)+\widetilde{c}_1(\alpha_{n(j+1)}, \varepsilon_{n(j+1)})-\widetilde{c}_1(\beta_{n(j+1)}, \varepsilon_{n(j+1)})\\
&=\sum_{i=0}^{n(1)-1}\big(\widetilde{c}_1(\varepsilon_{n(1)-i}, \alpha_{n(1)-i})-\widetilde{c}_1(\varepsilon_{n(1)-i-1}, \alpha_{n(1)-i})\big)\\
&\quad +\widetilde{c}_1(\varepsilon_0, \alpha_0)-\widetilde{c}_1(\alpha_0, \alpha_0)+\sum_{n=0}^{n(j+1)-1}\big(\widetilde{c}_1(\alpha_n,\varepsilon_n)-\widetilde{c}_1(\alpha_{n+1}, \varepsilon_n)\big)+\widetilde{c}_1(\alpha_{n(j+1)}, \varepsilon_{n(j+1)})-\widetilde{c}_1(\beta_{n(j+1)}, \varepsilon_{n(j+1)}). 
\end{split}
\end{equation}
Combining \cref{e:5.20e4}, \cref{e:5.20e5}, \cref{e:5.20e6}, \cref{e:Appendix1}, \cref{e:Appendix2} and \cref{e:5.20e7}, we have 
\begin{align*}
&\quad \sum_{k=1}^{N} \big(\widetilde{c}_1(x_k, y_k)-\widetilde{c}_1(x_{k+1}, y_k)\big)\\
&\leq \sum_{i=0}^{n(1)-1}\big(\widetilde{c}_1(\varepsilon_{n(1)-i}, \alpha_{n(1)-i})-\widetilde{c}_1(\varepsilon_{n(1)-i-1}, \alpha_{n(1)-i})\big)\\
&\quad + \widetilde{c}_1(\alpha_{n(j+1)}, \varepsilon_{n(j+1)}) -\widetilde{c}_1(\beta_{n(j+1)}, \varepsilon_{n(j+1)})+\sum_{n=n(j+1)}^{n(N+1)-1} \big(\widetilde{c}_1(\alpha_n, \varepsilon_n)-\widetilde{c}_1(\alpha_{n+1}, \varepsilon_n)\big)\\
&\quad +2\sum_{n=n(j+1)}^{n(N+1)} \big(\widetilde{c}_1(\alpha_{n}, \varepsilon_{n})-\widetilde{c}_1(\beta_{n}, \varepsilon_{n})\big)+\sum_{i=0}^{n(1)-1}\big(\widetilde{c}_1(\varepsilon_{n(1)-i}, \alpha_{n(1)-i})-\widetilde{c}_1(\varepsilon_{n(1)-i-1}, \alpha_{n(1)-i})\big)\\
&\quad +\widetilde{c}_1(\varepsilon_0, \alpha_0)-\widetilde{c}_1(\alpha_0, \alpha_0)+\sum_{n=0}^{n(j+1)-1}\big(\widetilde{c}_1(\alpha_n,\varepsilon_n)-\widetilde{c}_1(\alpha_{n+1}, \varepsilon_n)\big)\\
&\quad +\widetilde{c}_1(\alpha_{n(j+1)}, \varepsilon_{n(j+1)})-\widetilde{c}_1(\beta_{n(j+1)}, \varepsilon_{n(j+1)})\\
& = \widetilde{c}_1(\varepsilon_0, \alpha_0)-\widetilde{c}_1(\alpha_0, \alpha_0)+2\sum_{i=0}^{n(1)-1}\big(\widetilde{c}_1(\varepsilon_{n(1)-i}, \alpha_{n(1)-i})-\widetilde{c}_1(\varepsilon_{n(1)-i-1}, \alpha_{n(1)-i})\big)\\
&\quad +2\big(\widetilde{c}_1(\alpha_{n(j+1)}, \varepsilon_{n(j+1)}) -\widetilde{c}_1(\beta_{n(j+1)}, \varepsilon_{n(j+1)})\big)\\
&\quad +2\sum_{n=n(j+1)}^{n(N+1)} \big(\widetilde{c}_1(\alpha_{n}, \varepsilon_{n})-\widetilde{c}_1(\beta_{n}, \varepsilon_{n})\big)+\sum_{n=0}^{n(N+1)-1}\big(\widetilde{c}_1(\alpha_n, \varepsilon_n)-\widetilde{c}_1(\alpha_{n+1}, \varepsilon_n)\big)\\
&\leq \widetilde{c}_1(\varepsilon_0, \alpha_0)-\widetilde{c}_1(\alpha_0, \alpha_0)+2\sum_{i=0}^{n(1)-1}\big(\widetilde{c}_1(\varepsilon_{n(1)-i}, \alpha_{n(1)-i})-\widetilde{c}_1(\varepsilon_{n(1)-i-1}, \alpha_{n(1)-i})\big)\\
& \quad +4\sum_{n=n(j+1)}^{n(N+1)} \big(\widetilde{c}_1(\alpha_{n}, \varepsilon_{n})-\widetilde{c}_1(\beta_{n}, \varepsilon_{n})\big)+\sum_{n=0}^{n(N+1)-1}\big(\widetilde{c}_1(\alpha_n, \varepsilon_n)-\widetilde{c}_1(\alpha_{n+1}, \varepsilon_n)\big)\\
& \leq \widetilde{c}_1(\varepsilon_0, \alpha_0)-\widetilde{c}_1(\alpha_0, \alpha_0)+2\sum_{i=0}^{n(1)-1}\big(\widetilde{c}_1(\varepsilon_{n(1)-i}, \alpha_{n(1)-i})-\widetilde{c}_1(\varepsilon_{n(1)-i-1}, \alpha_{n(1)-i})\big)\\
& \quad +4\sum_{n=0}^{n(N+1)} \big(\widetilde{c}_1(\alpha_{n}, \varepsilon_{n})-\widetilde{c}_1(\beta_{n}, \varepsilon_{n})\big)+\sum_{n=0}^{n(N+1)-1}\big(\widetilde{c}_1(\alpha_n, \varepsilon_n)-\widetilde{c}_1(\alpha_{n+1}, \varepsilon_n)\big).
\end{align*}
The upper bound only depends on $n(1)$ and $n(N+1)$, which are determined by $(x_s, y_s), (x_e,y_e)$. 
\item [Case 4:] Suppose $(x_s, y_s)\in G_1$ and $(x_e, y_e)\in G_2$. Then $P^G_{(x_s,y_s)\to (x_e,y_e)}=\varnothing$ by \cref{ex:5.1}.(i). Hence, this case does not satisfy the overall assumption. 
\end{enumerate}

\end{document}